\documentclass[10pt]{amsart}
\usepackage{amscd, amssymb, amsfonts, amsmath, amstext, amsthm}

\newcounter{nomer}[section]
\newcommand{\n}{\par\refstepcounter{nomer}
{\arabic{section}.\arabic{nomer}.\ }}

\begin{document}

\title{Classification of semisimple Hopf algebras of dimension $16$.}
\author{Yevgenia Kashina}
\address{\hskip-\parindent
        Yevgenia Kashina\\
        Mathematical Sciences Research Institute\\
        1000 Centennial Drive\\
        Berkeley, CA 94720-5070}
\email{kashina@msri.org}

\thanks{The author was supported in part by NSF  
grants DMS-9701755 and  DMS 98-02086}
\newtheorem{theorem}{Theorem}[section]
\newtheorem{proposition}[theorem]{Proposition}
\newtheorem{lemma}[theorem]{Lemma}
\theoremstyle{definition}
\newtheorem{definition}[theorem]{Definition}
\theoremstyle{remark}
\newtheorem{remark}{Remark}[section]
\numberwithin{equation}{section}
\newcommand{\Char}{\mathop{ \rm char}\nolimits}

\begin{abstract}In this paper we  completely classify nontrivial semisimple Hopf 
algebras of dimension $16$. We also compute all the possible structures of 
the Grothendieck ring of
semisimple non-commutative Hopf algebras of dimension $16$. 
Moreover, we prove that non-commutative
semisimple Hopf algebras of dimension $p^{n}$, $p$ is prime, cannot have a
cyclic group of grouplikes.
\end{abstract}
\maketitle

\allowdisplaybreaks

\section{Introduction.}

Recently various classification results were obtained for finite-dimensional
semisimple Hopf algebras over an algebraically closed field of
characteristic $0$. The smallest dimension, for which the question was still
open, was $16$. In this paper we  completely classify all nontrivial 
(i.e. non-commutative and non-cocommutative) Hopf 
algebras of dimension $16$. Moreover, we consider all possible structures
of Grothendieck rings $K_{0}\left( H\right) $ for semisimple non-commutative
Hopf algebras of dimension $16$.

Let $H$ be a non-commutative semisimple Hopf algebra of dimension $16$ over
an algebraically closed field $k$ of characteristic $0$. Then irreducible
representations of $H$ of degree $1$ are exactly the grouplike elements of $
H^{\ast }$. Let $\mathbf{G}\left( H^{\ast }\right) $ denote the group of
grouplikes of $H^{\ast }$, then $k\mathbf{G}\left( H^{\ast }\right) $ is a
subHopfalgebra of $H^{\ast }$ and thus, by Nichols-Zoeller Theorem 
\cite{NZ}, $\left| \mathbf{G}\left( H^{\ast }
\right) \right| =\dim
k\mathbf{G}\left( H^{\ast }\right) $ divides $\dim H^{\ast}=
\dim H=16$.
Therefore by the Artin-Wedderburn Theorem, as an algebra $H$ is isomorphic to
either 
\begin{eqnarray}
&&k^{\left( 8\right) }\oplus M_{2}\left( k\right) \oplus M_{2}\left(
k\right) \text{ or}  \label{case1} \\
&&k^{\left( 4\right) }\oplus M_{2}\left( k\right) \oplus M_{2}\left(
k\right) \oplus M_{2}\left( k\right)  \label{case2}
\end{eqnarray}
$\dim Z\left( H\right) $ equals the number of summands in the
Artin-Wedderburn decomposition of $H$, thus in the case $\left( \ref{case1}
\right) $ $\dim Z\left( H\right) =10$ and $\left| \mathbf{G}\left( H^{\ast
}\right) \right| =8$ and in the case $\left( \ref{case2}\right) $ $\dim
Z\left( H\right) =7$ and $\left| \mathbf{G}\left( H^{\ast }\right) \right|
=4 $.

Our first result, which will be proved in the beginning of Section \ref{sec3}
, is the following:
\begin{theorem}
\label{th1}Let $H$ be a semisimple Hopf algebra of dimension $p^{n}$ over an
algebraically closed field $k$ of characteristic $0$. If $H\ncong kC_{p^{n}}$
then $\mathbf{G}\left( H\right) $ is not cyclic.
\end{theorem}

Our main result will be proved in Section \ref{summary}:
\begin{theorem}
\label{th2}Let $k$ be an algebraically closed field of characteristic $0$.
Then there are exactly $16$ nonisomorphic nontrivial semisimple Hopf
algebras of dimension $16$.
\begin{enumerate}
\item \  $\mathbf{G}\left( H\right) $ is abelian of order $8$ if and only if
$\mathbf{G}\left( H^{\ast }\right)$ is abelian of order $8$. There are $%
11$ Hopf algebras with such a group of grouplikes.
\item \  If $H$ has a nonabelian group of grouplikes then $\mathbf{G}\left(
H\right) =D_{8}$ and $\mathbf{G}\left( H^{\ast }\right) =C_{2}\times 
C_{2}$. There are $2$ Hopf algebras with a nonabelian group of
grouplikes.
\item \  There are $3$ Hopf algebras with $\mathbf{G}\left( H\right) =C%
_{2}\times C_{2}$. Two of them are dual to the Hopf algebras with a
nonabelian group of grouplikes and one of them is selfdual.
\end{enumerate}
\end{theorem}
\begin{remark}
A part of Theorem \ref{th2}, saying that if  $H$ has a nonabelian group of 
grouplikes then $\mathbf{G}\left(H\right) =D_{8}$ and 
$\mathbf{G}\left( H^{\ast }\right) =C_{2}\times 
C_{2}$ can also be obtained as a corollary to a theorem of Natale 
\cite{Na}, and Proposition \ref{th3}. This theorem states 
that if $\mathbf{G}\left( H\right) $ is nonabelian then $H^{\ast }$ has $4$ 
central grouplikes.
\end{remark} 
One method of constructing a new Hopf algebra from a known one $H$ is to
twist the comultiplication of $H$ by a $2$-pseudo-cocycle $\Omega \in
H\otimes H$ (or a $2$-cocycle $J\in H\otimes H$). The new Hopf algebra is
denoted $H_{\Omega }$ (or $H_{J}$). The next theorem summarizes the results of Sections \ref{sec5} and \ref{sec6}:
\begin{theorem}
\label{th4} Let $H$ be a semisimple Hopf algebra of dimension $16$ over an
algebraically closed field $k$ of characteristic $0$. Then there are
exactly $7$ possible structures of the Grothendieck ring 
$K_{0}\left(H\right)$. Moreover
\begin{enumerate}
\item \  $\mathbf{G}\left(H^{\ast }\right)$ is abelian if and only if the 
Grothendieck ring of $H$ is commutative. Then 
\begin{enumerate}
\item \  If $\left| \mathbf{G}\left( H^{\ast }\right) \right| =8$, as algebras
$K_{0}\left( H\right) \otimes _{\mathbb{Z}}k\cong k^{\left( 10\right) }$.
\item \  If $\left| \mathbf{G}\left( H^{\ast }\right) \right| =4$, as algebras
$K_{0}\left( H\right) \otimes _{\mathbb{Z}}k\cong k^{\left( 7\right) }$.
\item \ 
$K_{0}\left( H\right) =K_{0}\left( kG\right) $, where $G$ is one of the $9$
nonabelian groups of order $16$ (although only $6$ of those $K_{0}$-rings are
distinct).
\item \  $H$ is a twisting with a $2$-pseudo-cocycle of some
group algebra.
\end{enumerate}
\item \  If  
$K_{0}\left( H\right)$ is not commutative then
\begin{enumerate}
\item \ As algebras $K_{0}\left( H\right)\otimes _{\mathbb{Z}}k
\cong k^{\left(6\right) }\oplus M_{2}\left( k\right) $.
\item\  $H$ is not a twisting of a group algebra.
\item\  There is only one possible structure of the $K_{0}$-ring.
\item\  All Hopf algebras with non-commutative $K_{0}$-ring
are twistings of each other. 
\end{enumerate}
\end{enumerate}
\end{theorem}
\begin{remark}
By Theorem \ref{th2} there are only $2$ Hopf algebras with nonabelian 
$\mathbf{G}\left(H^{\ast }\right)$.
\end{remark}

We summarize the distinct non-commutative, non-cocommutative
semisimple Hopf algebras of dimension $16$ in the following table. We try
to distinguish nonisomorphic examples of Hopf algebras using the groups $%
\mathbf{G}\left( H\right) $ and $\mathbf{G}\left( H^{\ast }\right) $ and the
Grothendieck rings $K_{0}\left( H\right) $ (defined in Section \ref{sec2}).
Here we consider twistings of group algebras $kG$, where $G$ is a nonabelian
group of order $16$. There are exactly nine such groups, described in \cite{B}
(see Section \ref{sec4}). The twistings appearing here are explained in
Section \ref{seccoc}. The coproduct $\#^{\alpha }$ is explained
in Section \ref{seccoprod}. $H_{8}$ denotes the unique nontrivial
semisimple Hopf algebra of dimension $8$ (see \cite{KP} and \cite{Ma1}). 
\begin{table}[ht]
\caption{ } 
\begin{tabular}{|l|l|l|l|l|l|}
\hline
No. & Example & $\mathbf{G}\left( H\right) $ & $\mathbf{G}\left( H^{\ast
}\right) $ & $K_{0}\left( H\right) $ & Notes \\ \hline
$1$ & $H_{d:-1,1}\cong H_{8}\otimes kC_{2}$ & $\left( C_{2}\right) ^{3}$ & $%
\left( C_{2}\right) ^{3}$ & $K_{\ref{sec5}.\ref{5.1}}=K_{0}\left( D_{8}\times
C_{2}\right) $ & 
not triangular
\cr \hline
$2$ & $H_{d:1,1}\cong k\left( D_{8}\times C_{2}\right) _{J}$ & $\left(
C_{2}\right) ^{3}$ & $\left( C_{2}\right) ^{3}$ & $K_{\ref{sec5}.\ref{5.1}}=K_{0}\left(
D_{8}\times C_{2}\right) $ & triangular \cr \hline
$3$ & $\left( H_{c:\sigma _{1}}\right) ^{\ast }$ & $\left( C_{2}\right)
^{3} $ & $C_{2}\times C_{4}$ & $K_{\ref{sec5}.\ref{5.3}}=K_{0}\left( G_{5}\right) $ & 
\cr \hline
$4$ & $\left( H_{b:1}\right) ^{\ast }$ & $\left( C_{2}\right) ^{3}$ & $%
C_{2}\times C_{4}$ & $K_{\ref{sec5}.\ref{5.3}}=K_{0}\left( G_{5}\right) $ &  \cr \hline
$5$ & $H_{c:\sigma _{1}}$ & $C_{2}\times C_{4}$ & $%
\left( C_{2}\right) ^{3}$ & $K_{\ref{sec5}.\ref{5.2}}=K_{0}\left( G_{7}\right) $ &  \cr 
\hline
$6$ & $H_{b:1}$ & $C_{2}\times C_{4}$ & $\left( C_{2}\right) ^{3}$ & $K_{%
\ref{sec5}.\ref{5.1}}=K_{0}\left( D_{8}\times C_{2}\right) $ &  \cr \hline
$7$ & $H_{c:\sigma _{0}}$ & $C_{2}\times C_{4}$ & $C_{2}\times C_{4}$ & $K_{%
\ref{sec5}.\ref{5.4}}=K_{0}\left( G_{1}\right) $ & \cr \hline
$8$ & $H_{a:1}$ & $C_{2}\times C_{4}$ & $C_{2}\times C_{4}$ & $K_{\ref{sec5}.\ref{5.3}%
}=K_{0}\left( G_{5}\right) $ &  \cr \hline
$9$ & $H_{a:y}$ & $C_{2}\times C_{4}$ & $C_{2}\times C_{4}$ & $K_{\ref{sec5}.\ref{5.3}%
}=K_{0}\left( G_{5}\right) $ &  \cr \hline
$10$ & $H_{b:y}$ & $C_{2}\times C_{4}$ & $C_{2}\times C_{4}$ & $K_{\ref{sec5}.\ref
{5.3}}=K_{0}\left( G_{5}\right) $ &  \cr \hline
$11$ & $H_{b:x^2y}$ & $C_{2}\times C_{4}$ & $C_{2}\times C_{4}$ & $K_{\ref{sec5}.\ref
{5.3}}=K_{0}\left( G_{5}\right) $ &  \cr \hline
$12$ & $ H_{C:1}\cong \left( kD_{16}\right) _{J}$ & $D_{8}$ & 
$C_{2}\times C_{2}$ & $K_{\ref{sec6}.\ref
{6.3}}=K_{0}\left( D_{16}\right)$ & triangular \cr \hline
$13$ & $H_{E}\cong \left( kG_{2}\right) _{J}$ & $D_{8}$ & $C_{2}\times C_{2}$ 
& $K_{\ref{sec6}.\ref
{6.4}}=K_{0}\left( G_{2}\right)$ & triangular \cr \hline
$14$ & $ H_{B:1}\cong\left( \left( kD_{16}\right) _{J}\right) ^{\ast }$ & $%
C_{2}\times C_{2}$ & $D_{8}$ & $K_{\ref{sec5}.\ref
{5.5}}$  &  \cr \hline
$15$ & $H_{B:X}\cong\left(\left( kG_{2}\right) _{J}\right) 
^{\ast }$ & $C_{2}\times C_{2}$ & $D_{8}$ & $K_{\ref{sec5}.\ref
{5.5}}$ &$ kQ_{8}\#^{\alpha }kC_{2}$\cr \hline
$16$ & $H_{C:\sigma _{1}}$ & $%
C_{2}\times C_{2}$ & $C_{2}\times C_{2}$ & $K_{\ref{sec6}.\ref
{6.3}}=K_{0}\left( D_{16}\right)$  & not triangular \cr \hline
\end{tabular}
\end{table}
\begin{remark}
$H_{C:\sigma _{1}}$ is not triangular for the following reasons. If it were 
triangular then by \cite[Theorem 2.1]{EG} it would be 
equal to a twisting with a $2$-cocycle of a group algebra $kG$. Then by 
\cite[Theorem 4.1]{N1}
$K_{0}\left( H_{C:\sigma _{1}}\right)=K_{0}\left( kG\right)$ and therefore
$H_{C:\sigma _{1}}$ would be a twisting of $kD_{16}$ or of $kQ_{16}$. 
But by \cite[Theorem 4.1]{Ma7}, $kQ_{16}$ doesn't have nontrivial cocycle 
twistings and $ H_{C:1}\cong \left( kD_{16}\right) _{J}$ is the only cocycle 
twisting of $kD_{16}$. 
\end{remark}
\begin{remark}
The following Hopf algebras are selfdual: $H_{d:-1,1}\cong H_{8}\otimes kC_{2}$
(since $H_{8}$ is selfdual), $H_{c:\sigma _{0}}$ (since comparing $K_{0}$-rings
we see that $H_{c:\sigma _{0}}\cong A_{3}^{+} \cong 
\left( A_{3}^{+}\right) ^{\ast }$, described in \cite{K1} and \cite{K2}),
$H_{d:1,1}\cong k\left( D_{8}\times C_{2}\right) 
_{J}$ and $H_{C:\sigma _{1}}$ (since there is no other choice for the dual).
\end{remark}
\section{Preliminaries.} \label{sec2}

First we will need the following definition, which was introduced in \cite{S}%
.
\begin{definition}
Let $K_{0}\left( H\right) ^{+}$ denote the abelian semigroup of all
equivalence classes of representations of $H$ with the addition given by a
direct sum. Then its enveloping group $K_{0}\left( H\right) $ has the
structure of an ordered ring with involution $^{\ast }$ and is called the
Grothendieck ring.
\end{definition}
In \cite{NR} the structure of $K_{0}\left( H\right) $ was described for
comodules; it was then translated into the language of modules in \cite{N2}.
The multiplication in this ring is defined as follows: let $\left[ \pi _{1}%
\right] $ and $\left[ \pi _{2}\right] $ denote the classes of
representations equivalent to $\pi _{1}$ and $\pi _{2}$, then $\left[ \pi
_{1}\right] \bullet \left[ \pi _{2}\right] $ is the class of the
representation $\left( \pi _{1}\otimes \pi _{2}\right) \circ \Delta $; the
unit of this ring is the class $\left[ \varepsilon \right] $ and $\left[ \pi %
\right] ^{\ast }$ is the equivalence class of the dual representation $%
^{t}\left( \pi \circ S\right) $ defined by $\left\langle ^{t}\left( \pi
\circ S\left( h\right) \right) \left( f\right) ,v\right\rangle =\left\langle
f,\left( \pi \circ S\left( h\right) \right) \left( v\right) \right\rangle $.
The equivalence classes of irreducible representations of $H$ form a basis
of $K_{0}\left( H\right) $ and are called \textit{basic elements}. If $\left[
\pi _{1}\right] ,\ldots ,\left[ \pi _{d}\right] $ are the basic elements
then $\left[ \rho \right] =\sum_{i=1}^{d}\deg \pi _{i}\left[ \pi _{i}\right] 
$ is called the \textit{marked element}. For basic elements $x$ and $y$ we
write 
\begin{equation*}
x\bullet y=\sum_{z\text{ - basic}}m\left( z,x\bullet y\right) z
\end{equation*}
where $m\left( z,x\bullet y\right) $ are non-negative integers. Then the
following properties are true (see \cite{NR} and \cite{N2}): 
\begin{eqnarray}
m\left( z,x\bullet y\right) &=&m\left( x^{\ast },y\bullet z^{\ast }\right) \\
m\left( 1,x\bullet y^{\ast }\right) &=&\delta _{x,y} \\
\sum m\left( z,x\bullet y\right) \deg \left( z\right) &=&\deg \left(
x\bullet y\right)
\end{eqnarray}
For simplicity of notation we will write $\pi $ instead of $\left[ \pi %
\right] $ for elements of $K_{0}\left( H\right) $. We will denote the degree 
$2$ irreducible representations of $H$ by $\pi _{i}$ and the degree $1$
irreducible representations of $H$ (i.e. elements of $\mathbf{G}\left(
H^{\ast }\right) $ or multiplicative characters of $H$) by $\chi _{i}$. We
denote the generators of $\mathbf{G}\left( H^{\ast }\right) $ by $\chi $, 
$\varphi $ and $\psi $. If $H=kG$ then $\mathbf{G}\left( \left( kG\right)
^{\ast }\right) $ is the group of multiplicative characters of $G$.

The following proposition can be also obtained as a corollary to 
\cite[Proposition 2.4]{Ma7}:
\begin{proposition}\label{com_quot}
Let $H$ be a nontrivial semisimple Hopf algebra of dimension $16$. 
Assume that there exists an element $\chi \in \mathbf{G}\left( H^{\ast }\right)
\cap Z\left( H^{\ast }\right)$ of order $2$ such that $\chi \bullet \pi =\pi$
for every $2$-dimensional representation $\pi $ of $H$. Then $H^{\ast }$ has
a group algebra of dimension $8$ as a quotient.
\end{proposition}
\proof
Write $G=\mathbf{G}\left( H^{\ast }\right)$. Dualizing formulas 
(\ref{case1}) and (\ref{case2}) we get that as coalgebras
\begin{eqnarray*}
H^{\ast }&=&kG \oplus E_{1} \oplus E_{2} \qquad \ \qquad\text{if } 
\left|\mathbf{G}\left( H^{\ast }\right) \right|=8 \ \text{ or}\\
H^{\ast }&=&kG \oplus E_{1}\oplus E_{2}\oplus 
E_{3} \qquad \text{if } \left| \mathbf{G}\left( H^{\ast }\right) \right|=4
\end{eqnarray*}
where $E_{i}$ are simple subcoalgebras of dimension $4$ and $\chi E_{i}=E_{i}$.
$\left( \chi -1\right)H^{\ast }$ is a normal Hopf ideal of $H^{\ast }$.
Then $L=H^{\ast }/ \left( \chi -1\right)H^{\ast }$ is a Hopf algebra of 
dimension $8$. Consider the projection $p:H^{\ast } \to L$. 
Since $\chi E_{i}=E_{i}$, $p\left(E_{i}\right)=E_{i}/ \left( \chi -1\right)
E_{i}$. Therefore
\begin{eqnarray*}
L&=&k\left( G/<\chi > \right) \oplus p\left(E_{1}\right) \oplus 
p\left(E_{2}\right) \qquad\qquad 
\qquad\text{if }\left|\mathbf{G}\left( H^{\ast }\right)\right| =8\ \text{ or}\\
L&=&k\left( G/<\chi > \right)\oplus p\left(E_{1}\right) \oplus 
p\left(E_{2}\right) \oplus  p\left(E_{3}\right)
\qquad \text{if } \left| \mathbf{G}\left( H^{\ast }\right) \right|=4
\end{eqnarray*}
$p\left(E_{i}\right)$ are cosemisimple coalgebras of dimension $2$, therefore 
each of them is spanned by two grouplikes. Thus $L$ is spanned by $8$ 
grouplikes and $L$ is a group algebra.
\endproof
We will also need the notion of a twisting of a Hopf algebra (see \cite{D}, 
\cite{V}, \cite{N1}):

\begin{definition}
The \textit{twisting} $H_{\Omega }$ of a Hopf algebra $H$ is a Hopf algebra
with the same algebra structure and counit and with comultiplication and
antipode given by 
\begin{eqnarray*}
\Delta _{\Omega }\left( h\right) =\Omega \Delta \left( h\right) \Omega ^{-1}
\\
S_{\Omega }\left( h\right) =uS\left( h\right) u^{-1}
\end{eqnarray*}
for all $h\in H$, where $\Omega \in H\otimes H$ and $u\in H$ are invertible
elements.
\end{definition}

The new comultiplication $\Delta _{\Omega }$ is coassociative if and only if 
$\Omega $ is a $2$-\textit{pseudo-cocycle}, that is $\partial _{2}\left(
\Omega \right) $ lies in the centralizer of $\left( \Delta \otimes id\right)
\Delta \left( H\right) $ in $H\otimes H\otimes H$, where 
\begin{equation*}
\partial _{2}\left( \Omega \right) =\left( id\otimes \Delta \right) \left(
\Omega ^{-1}\right) \left( 1\otimes \Omega ^{-1}\right) \left( \Omega
\otimes 1\right) \left( \Delta \otimes id\right) \left( \Omega \right) .
\end{equation*}
$\Omega $ is called a $2$-\textit{cocycle} if $\partial _{2}\left( \Omega
\right) =1\otimes 1\otimes 1$ and in this case we will denote it by $J$.

\begin{remark}
By \cite[Theorem 4.1]{N1} $K_{0}\left( H\right) \cong K_{0}\left( H_{\Omega
}\right) $ as ordered rings with marked elements, and thus $\mathbf{G}\left(
H^{\ast }\right) \cong $ $\mathbf{G}\left( \left( H_{\Omega }\right) ^{\ast
}\right) $.
\end{remark}

\section{Hopf algebras of dimension $16$ with a commutative
subHopfalgebra of dimension $8$.} \label{sec3}

We apply the methods used by Masuoka in \cite{Ma1}, \cite{Ma2} and \cite{Ma4}. 
Let $H$ be a nontrivial semisimple Hopf algebra of dimension $16$ with a
subHopfalgebra $K=\left( kG\right) ^{\ast }$ of dimension $8$. 
Since $K$ is a subHopfalgebra of index $2$, by \cite[Proposition 2]{Ma6}
or \cite[Theorem 2.1.1]{Na1}
$K$ is normal in $H$ and thus we have an exact sequence of Hopf algebras 
\begin{equation}
K\stackrel{i}{\hookrightarrow }H\stackrel{\pi }{\twoheadrightarrow }F
\label{ext1}
\end{equation}
where $F=k\left\langle t\right\rangle \cong kC_{2}$ and $K=\left( kG\right)
^{\ast }$, which is cleft by \cite{Sch1} or \cite{MaD}. Such a sequence is 
called an extension of $F$ by $K$ and was first studied by Kac in 
\cite{Ka}. The construction of extensions from cohomological data was done in
\cite{Mj} and \cite{AD}. $K$ is commutative and $F$ is cocommutative and thus 
$\left( F,K\right) $ form an abelian matched pair of Hopf algebras and $\left(
G,\left\langle t\right\rangle \right) $ form an abelian matched pair of
groups (see \cite{T}, \cite{H}, \cite{Si}, \cite[Section 1]{Ma2}). Therefore 
$H$ becomes a bicrossed product $K\#_{\sigma }^{\theta }F$ with an action $%
\rightharpoonup :F\otimes K\rightarrow K$, a coaction $\rho :F\rightarrow
F\otimes K$, a cocycle $\sigma :F\otimes F\rightarrow K$ and a dual cocycle $%
\theta :F\rightarrow K\otimes K$. $G$ is a normal subgroup of the group $%
G\times \left\langle t\right\rangle $, arising from a matched pair $\left(
G,\left\langle t\right\rangle \right) $, since $G$ has index $2$ in $G\times
\left\langle t\right\rangle $. Thus $\rho $ is trivial and the action by $t$
is a Hopf algebra automorphism of $K$ (see \cite[Section 1]{Ma2}).
$\rightharpoonup $ is a nontrivial action on $K$, since otherwise $H\cong
K^{t}\left[ C_{2}\right] $ as an algebra, and thus $H$ is commutative. 

Let $v=\sigma \left( t,t\right) \in K.$ Then by the properties of the
cocycle $v$ is a unit and 
\begin{equation}
t\rightharpoonup v=v  \label{coc1}
\end{equation}
Multiplication in $H$ gives us 
\begin{eqnarray}
\overline{t}^{2} &=&v  \label{coc2} \\
\overline{t}c &=&\left( t\rightharpoonup c\right) \overline{t}
\label{coc3}
\end{eqnarray}
where $\overline{t}=1\#t$ and $c \in K$.

Moreover, if a unit $v\in K$ satisfies $\left( \ref{coc1}\right) $, $\left( 
\ref{coc2}\right) $ and $\left( \ref{coc3}\right) $, we can define a cocycle 
$\sigma $ by $\sigma \left( 1,1\right) =\sigma \left( 1,t\right) =\sigma
\left( t,1\right) =1$ and $\sigma \left( t,t\right) =v.$

We proceed by considering the possible $G$,
namely $C_{8}$, $C_{4}\times C_{2}$, $C_{2}\times C_{2}\times C_{2}$,
$D_{8}$ and $Q_{8}$.
Theorem \ref{th1} says that the first case cannot appear.

\begin{proof}[\bf{Theorem \ref{th1}}]
Let's prove the
statement by induction on $n$. When $n=2$, by \cite[Theorem 2]{Ma3} $H$ is a 
group algebra
and if $H\ncong kC_{p^{2}}$ then $H\cong k\left( C_{p}\times
C_{p}\right) $ and $\mathbf{G}\left( H\right) \cong C_{p}\times C_{p}$.

Now assume the statement is true for $n=m$. Consider $H$ of dimension $%
p^{m+1}$. $\dim \left( H^{\ast }\right) =p^{m+1}$ and thus, by \cite[Theorem 1]
{Ma3},
there exists a central grouplike of order $p$ in $H^{\ast }$ and therefore $%
H^{\ast }$ contains a normal subHopfalgebra $K\cong kC_{p}$. Thus we get a
short exact sequence of Hopf algebras 
\begin{equation}
K\overset{i}{\hookrightarrow }H^{\ast }\overset{\pi }{\twoheadrightarrow }F
\label{shortdual}
\end{equation}
where $F=H^{\ast }/K^{+}H^{\ast }$. Dualizing $\left( \ref{shortdual}\right) 
$ we get another short exact sequence of Hopf algebras: 
\begin{equation}
F^{\ast }\overset{\pi ^{\ast }}{\hookrightarrow }H\overset{i^{\ast }}{%
\twoheadrightarrow }K^{\ast }  \label{short}
\end{equation}
where $K^{\ast }\cong K\cong kC_{p}$ and $\dim F^{\ast }=\dim F=p^{m}$. Thus
we get $\mathbf{G}\left( F^{\ast }\right) \subseteq \mathbf{G}\left(
H\right) $ and $\mathbf{G}\left( F^{\ast }\right) $ is not cyclic unless $%
F^{\ast }\cong kC_{p^{m}}$. In the first case we are done since it implies
that $\mathbf{G}\left( H\right) $ is not cyclic. In the second case, since $%
K $ is normal in $H^{\ast }$, $H^{\ast }$ is isomorphic as an algebra to a
twisted group ring $K^{t}\left[ F\right] $ where $F\cong F^{\ast }\cong
kC_{p^{m}}$. It is easy to show that, since$\ F$ is a group algebra of a
cyclic group, $K^{t}\left[ F\right] $ is commutative. Thus $H$ is
cocommutative and the only possible $H$ with a cyclic group of grouplikes is $%
kC_{p^{m+1}}$. 
\end{proof}

\subsection{Case of $\mathbf{G}\left( H\right) =C_{4}\times C_{2}$.\label%
{3.1}}
\ \newline
We will show that there are at most $7$ possible Hopf algebras of this kind.
Let $H$ be a nontrivial semisimple Hopf algebra of dimension $16$ with a
subHopfalgebra $K=k\left( C_{4}\times C_{2}\right) ^{\ast } 
\cong k\left( C_{4}\times C_{2}\right) $. Then
$\mathbf{G}\left( H\right) =G \cong C_{4}\times C_{2}$.
 
Let $G=\left\langle x\right\rangle \times \left\langle
y\right\rangle $ with $\left| x\right| =4$ and $\left| y\right| =2$. Then
the dual basis of $K\cong K^{\ast }$ is given by 
\begin{equation*}
e_{pq}=1/8\left( 1+i^{p}x+i^{2p}x^{2}+i^{3p}x^{3}\right) \left( 1+\left(
-1\right) ^{q}y\right) ,\text{\qquad }p=0,1,2,3;q=0,1
\end{equation*}
Then 
\begin{eqnarray*}
\Delta _{H}\left( e_{pq}\right) &=&\Delta _{K}\left( e_{pq}\right) =\sum 
_{\substack{ p_{1}+p_{2}\equiv p\mod 4  \\ %
q_{1}+q_{2}\equiv q\mod 2 }}e_{p_{1}q_{1}}\otimes
e_{p_{2}q_{2}} \\
\Delta _{H}\left( \overline{t}\right) &=&\theta \left( t\right) \overline{t}%
\otimes \overline{t}
\end{eqnarray*}
where $\overline{t}=1\#t.$ Dualizing $\left( \ref{ext1}\right) $ we get
another extension 
\begin{equation*}
F^{\ast }\overset{\pi ^{\ast }}{\hookrightarrow }H^{\ast }\overset{i^{\ast }%
}{\twoheadrightarrow }K^{\ast }
\end{equation*}
and as in \cite[2.4]{Ma1}, \cite[2.11]{Ma2} or \cite[2.1]{Ma5}, since $k$ is
algebraically closed, there exist units $\overline{x}$ and $\overline{y}\in
H^{\ast }$, such that $\overline{x}^{4}=\overline{y}^{2}=1_{H^{\ast }}$, $%
\left\langle e_{pq},\overline{x}^{i}\overline{y}^{j}\right\rangle =\delta
_{ip}\delta _{jq}$ and $\alpha =\overline{y}^{-1}\overline{x}^{-1}\overline{y%
}\overline{x}\in F^{\ast }=k\left\{ e_{0},e_{1}\right\} $, where $\left\{
e_{r}\right\} $ is a dual basis of $\left\{ t^{r}\right\} $, $r=0,1$. $%
\varepsilon \left( \alpha \right) =\varepsilon \left( \overline{y}^{-1}%
\overline{x}^{-1}\overline{y}\overline{x}\right) =1$ and therefore $\alpha
=e_{0}+\xi e_{1}$. The right action $\rho ^{\ast }:F^{\ast }\otimes K^{\ast
}\rightarrow F^{\ast }$ is trivial, thus $F^{\ast }$ lies in the center of $%
H^{\ast }$. 
\begin{equation*}
\overline{x}=\overline{y}^{2}\overline{x}=\overline{y}\overline{x}\overline{y%
}\alpha =\overline{x}\overline{y}\alpha \overline{y}\alpha =\overline{x}%
\overline{y}^{2}\alpha ^{2}=\overline{x}\alpha ^{2}
\end{equation*}
Thus $\alpha ^{2}=1$ and therefore $\xi =\pm 1$. 
\begin{equation*}
\left\langle \Delta _{H}\left( \overline{t}\right) ,\overline{x}^{i}%
\overline{y}^{j}e_{k}\otimes \overline{x}^{p}\overline{y}^{q}e_{r}\right%
\rangle =\left\langle \overline{t},\overline{x}^{i}\overline{y}^{j}e_{k}%
\overline{x}^{p}\overline{y}^{q}e_{r}\right\rangle =\delta _{kr}\left\langle 
\overline{t},\overline{x}^{i+p}\overline{y}^{j+q}\alpha
^{jp}e_{k}\right\rangle =\xi ^{jp}\delta _{k1}\delta _{r1}
\end{equation*}
On the other hand 
\begin{equation*}
\left\langle \Delta _{H}\left( \overline{t}\right) ,\overline{x}^{i}%
\overline{y}^{j}e_{k}\otimes \overline{x}^{p}\overline{y}^{q}e_{r}\right%
\rangle =\left\langle \theta \left( t\right) \overline{t}\otimes \overline{t}%
,\overline{x}^{i}\overline{y}^{j}e_{k}\otimes \overline{x}^{p}\overline{y}%
^{q}e_{r}\right\rangle =\left\langle \theta \left( t\right) ,\overline{x}^{i}%
\overline{y}^{j}\otimes \overline{x}^{p}\overline{y}^{q}\right\rangle \delta
_{k1}\delta _{r1}
\end{equation*}
Therefore 
\begin{equation*}
\theta \left( t\right) =\sum_{ijpq}\xi ^{jp}e_{ij}\otimes e_{pq}
\end{equation*}
and since $H$ should be non-cocommutative, $\theta \left(
t,t\right) $ is nontrivial, and thus $\xi =-1$ and 
\begin{equation*}
\theta \left( t\right) =\sum_{ijpq}\left( -1\right) ^{jp}e_{ij}\otimes
e_{pq}=
{\frac12}
\left( \left( 1+y\right) \otimes 1+\left( 1-y\right) \otimes x^{2}\right)
\end{equation*}

Write $v=\sigma
\left( t,t\right) =\sum c_{i,j}e_{i,j}$ then $c_{0,0}=\varepsilon \left(
v\right) =1$ and $c_{i,j}\neq 0,$ since $v$ is a unit, and 
\begin{equation*}
\Delta _{H}\left( \overline{t}^{2}\right) =\Delta _{H}\left( v\right)
=\Delta _{K}\left( \sum c_{i,j}e_{i,j}\right) =\sum
c_{i+p,j+q}e_{i,j}\otimes e_{p,q}
\end{equation*}
On the other hand, if we write 
\begin{equation*}
t\rightharpoonup e_{p,q} =e_{\alpha _{1}\left( p,q\right) ,\alpha _{2}\left(
p,q\right) }
\end{equation*}
\begin{eqnarray*}
\Delta \left( \overline{t}\right)\Delta \left( \overline{t}\right) &=&
\sum_{ijpq}\left( -1\right) ^{jp}e_{ij} \overline{t}\otimes e_{pq} \overline{t}
\sum_{ijpq}\left( -1\right) ^{jp}e_{ij} \overline{t}\otimes e_{pq}
\overline{t}\\
&=& \sum_{ijpq}\left( -1\right) ^{jp}e_{ij}\otimes e_{pq}
\sum_{ijpq}\left( -1\right) ^{jp}
e_{\alpha _{1}\left( i,j\right) ,\alpha _{2}\left(
i,j\right) } \overline{t}^{2}\otimes 
e_{\alpha _{1}\left( p,q\right) ,\alpha _{2}\left( p,q\right) }
\overline{t}^{2}\\
&=& \sum_{ijpq}\left( -1\right) ^{jp}e_{ij}\otimes e_{pq}
\sum_{ijpq}\left( -1\right) ^{\alpha _{2}\left( i,j\right)
\alpha _{1}\left( p,q\right) }e_{ij}\overline{t}^{2}\otimes 
 e_{pq}\overline{t}^{2}\\ 
&=&
\sum_{ijpq}\left( -1\right) ^{jp + \alpha _{2}\left( i,j\right)
\alpha _{1}\left( p,q\right) }c_{ij}c_{pq}e_{ij}\otimes e_{pq}
\end{eqnarray*}
Thus for $H$ to be a bialgebra we should have
\begin{equation}
c_{i+p,j+q}=\left( -1\right) ^{jp + \alpha _{2}\left( i,j\right)
\alpha _{1}\left( p,q\right) }c_{i,j}c_{p,q}  \label{coc4}
\end{equation}

Action by $t$ is a Hopf algebra map and therefore $t\rightharpoonup G=G$
and $f_{t}:G\rightarrow G$ defined by $f_{t}\left( g\right)
=t\rightharpoonup g$ is a group automorphism of order $2$. There are three
possibilities for such an automorphism; we consider them below: 

\subsubsection*{Case a)}
The action is given by 
\begin{eqnarray*}
t &\rightharpoonup &x=xy \\
t &\rightharpoonup &y=y
\end{eqnarray*}
Then $t\rightharpoonup e_{i,j}=e_{i+2j,j}$. Write $v=\sigma \left(
t,t\right) =\sum c_{i,j}e_{i,j}$ . By $\left( \ref{coc1}\right) $ and
$\left( \ref{coc4}\right) $ 
\begin{eqnarray}
c_{i,j}&=&c_{i+2j,j}  \label{cond3}\\
c_{i+p,j+q}&=&c_{i,j}c_{p,q}  \label{cond4}
\end{eqnarray}
Conditions $\left( \ref{cond3}\right) $ and $\left( \ref{cond4}\right) $
imply that $c_{1,0}=\left( -1\right) ^k$ and 
$c_{0,1}=\left( -1\right) ^l$ for $k,l=0,1$ and 
\begin{eqnarray*}
\sigma \left( t,t\right) =\sum _{p,q}\left( -1\right)^{kp+lq}e_{p,q}=
\sum \left( -1\right)^{kp}e_{p,q}\sum \left( -1\right)^{lq}e_{p,q}
=x^{2k}y^{l} \qquad k,l=0,1
\end{eqnarray*}
For $k,l=0,1$ let $H_{k,l}$ be the Hopf algebras with the structures
described above with cocycles $\sigma _{k,l}\left( t,t\right) =x^{2k}y^{l}$.
Define 
\begin{eqnarray*}
f:H_{k,l} &\rightarrow &H_{k+1,l+1}\qquad \text{by} \\
f\left( e_{r,s}\right) &=&e_{r,s} \\
f\left( \overline{t}\right) &=&x\overline{t}
\end{eqnarray*}
and extend it multiplicatively to $f\left( e_{r,s}\overline{t}\right) $.
Then $f$ is a trivial group homomorphism on $\mathbf{G}\left( H_{k,l}\right) $
 and 
\begin{eqnarray*}
f\left( \overline{t}\right) f\left( \overline{t}\right) &=&x\overline{t}x%
\overline{t}=x^{2}y\overline{t}^{2}
=x^{2}yx^{2\left( k+1\right) }y^{\left( l+1\right) }=x^{2k}y^{l}=f\left( 
\overline{t}^{2}\right) \\
f\left( \overline{t}x\right) &=&f\left( xy\overline{t}\right) =xyx\overline{t%
}=x\overline{t}x=f\left( \overline{t}\right) f\left( x\right)\\
\left( f\otimes f\right) \Delta \left( \overline{t}\right) &=&\left(
f\otimes f\right) \left( \theta \left( t\right) \overline{t}\otimes 
\overline{t}\right) =\theta \left( t\right) \left( f\left( \overline{t}%
\right) \otimes f\left( \overline{t}\right) \right) =\theta \left( t\right)
\left( x\overline{t}\otimes x\overline{t}\right) \\
&=&\left( x\otimes x\right) \theta \left( t\right) \left( \overline{t}%
\otimes \overline{t}\right) =\Delta \left( x\right) \Delta \left( \overline{t%
}\right) =\Delta \left( x\overline{t}\right) =\Delta \left( f\left( 
\overline{t}\right) \right)
\end{eqnarray*}
and such an $f$ is a Hopf algebra isomorphism between $H_{k,l}$ and $%
H_{k+1,l+1}$. Thus there are at most two nonisomorphic Hopf algebras of this
type:

\begin{enumerate}
\item  $H_{a:1}=H_{0,0}$ with the trivial cocycle and $\mathbf{G}\left(
H_{a:1}^{\ast }\right) =\left\langle \chi \right\rangle \times \left\langle
\varphi \right\rangle \cong C_{4}\times C_{2}$, where $\chi (x)=i$, $\chi
(y)=\chi (t)=1$, $\varphi \left( x\right) =\varphi \left( y\right) =1$, $%
\varphi \left( t\right) =-1$. There is a degree $2$ irreducible
representation defined by
\newline
\begin{tabular}{lll}
$\pi (x)=\left( 
\begin{array}{rr}
i & 0 \\ 
0 & -i
\end{array}
\right) $ & $\pi (y)=\left( 
\begin{array}{rr}
-1 & 0 \\ 
0 & -1
\end{array}
\right) $ & $\pi (t)=\left( 
\begin{array}{rr}
0 & 1 \\ 
1 & 0
\end{array}
\right) $%
\end{tabular}
\newline
with the property $\pi ^{2}=\pi \bullet \pi =1+\chi ^{2}+\varphi +\chi
^{2}\varphi $.
\smallskip

\item  $H_{a:y}=H_{0,1}$ with the cocycle defined by $\sigma \left(
t,t\right) =y$ and $\mathbf{G}\left( H_{a:y}^{\ast }\right) =\left\langle
\chi \right\rangle \times \left\langle \varphi \right\rangle \cong
C_{4}\times C_{2}$, where $\chi (x)=i$, $\chi (y)=\chi (t)=1$, $\varphi
\left( x\right) =\varphi \left( y\right) =1$, $\varphi \left( t\right) =-1$.
There is a degree $2$ irreducible representation defined by
$\newline
\begin{tabular}{lll}
$\pi (x)=\left( 
\begin{array}{rr}
i & 0 \\ 
0 & -i
\end{array}
\right) $ & $\pi (y)=\left( 
\begin{array}{rr}
-1 & 0 \\ 
0 & -1
\end{array}
\right) $ & $\pi (t)=\left( 
\begin{array}{rr}
0 & -1 \\ 
1 & 0
\end{array}
\right) $%
\end{tabular}
\newline$
with the property $\pi ^{2}=\pi \bullet \pi =1+\chi ^{2}+\varphi +\chi
^{2}\varphi $.
\end{enumerate}

\subsubsection*{Case b)}

The action is given by 
\begin{eqnarray*}
t &\rightharpoonup &x=x^{-1} \\
t &\rightharpoonup &y=y
\end{eqnarray*}
Then $t\rightharpoonup e_{i,j}=e_{-i,j}$. Write $v=\sigma \left( t,t\right)
=\sum c_{i,j}e_{i,j}$ . By $\left( \ref{coc1}\right) $ and
$\left( \ref{coc4}\right) $  
\begin{eqnarray}
c_{i,j}&=&c_{-i,j}  \label{cond5}\\
c_{i+p,j+q}&=&c_{i,j}c_{p,q}  \label{cond6}
\end{eqnarray}
Conditions $\left( \ref{cond5}\right) $ and $\left( \ref{cond6}\right) $
imply that  $c_{1,0}=\left( -1\right) ^k$ and 
$c_{0,1}=\left( -1\right) ^l$ for $k,l=0,1$ and 
\begin{eqnarray*}
\sigma \left( t,t\right) =\sum _{p,q}\left( -1\right)^{kp+lq}e_{p,q}=
\sum \left( -1\right)^{kp}e_{p,q}\sum \left( -1\right)^{lq}e_{p,q}
=x^{2k}y^{l} \qquad k,l=0,1
\end{eqnarray*}
For $k,l=0,1$ let $H_{k,l}$ be the Hopf algebras with the structures
described above with cocycles $\sigma _{k,l}\left( t,t\right) =x^{2k}y^{l}$.
Define 
\begin{eqnarray*}
f:H_{0,0} &\rightarrow &H_{1,0}\qquad \text{by} \\
f\left( e_{r,s}\right) &=&e_{r,r+s} \\
f\left( \overline{t}\right) &=&%
{\frac12}
\left( \left( 1+i\right) 1+\left( 1-i\right) x^{2}\right) \overline{t}%
=\sum_{p=0}^{3}\sum_{q=0}^{1}i^{p^{2}}e_{p,q}\overline{t}
\end{eqnarray*}
and extend it multiplicatively to $f\left( e_{r,s}\overline{t}\right) $.
Then $f\mid _{\mathbf{G}\left( H_{0,0}\right) }$ is a group isomorphism \newline
$%
\mathbf{G}\left( H_{0,0}\right) \rightarrow \mathbf{G}\left( H_{1,0}\right) $
with $f\left( x\right) =x$, $f\left( y\right) =x^{2}y$ and 
\begin{eqnarray*}
f\left( \overline{t}\right) f\left( \overline{t}\right) &=&
{\frac14}
\left( \left( 1+i\right) 1+\left( 1-i\right) x^{2}\right) \overline{t}\left(
\left( 1+i\right) 1+\left( 1-i\right) x^{2}\right) \overline{t} \\
&=&%
{\frac14}
\left( \left( 1+i\right) 1+\left( 1-i\right) x^{2}\right) ^{2}\overline{t}%
^{2}=
{\frac14}
\left( 2i\cdot 1+4x^{2}-2i\cdot 1\right) \sigma _{1,0}\left( t,t\right) \\
&=&x^{2}x^{2}=1=f\left( \overline{t}^{2}\right) \\
f\left( \overline{t}x\right) &=&f\left( x^{-1}\overline{t}\right) =%
{\frac12}
x^{-1}\left( \left( 1+i\right) 1+\left( 1-i\right) x^{2}\right) \overline{t}
\\
&=&
{\frac12}
\left( \left( 1+i\right) 1+\left( 1-i\right) x^{2}\right) \overline{t}%
x=f\left( \overline{t}\right) f\left( x\right) \\
\left( f\otimes f\right) \Delta \left( \overline{t}\right) &=&\left(
f\otimes f\right) \left( \theta \left( t\right) \overline{t}\otimes 
\overline{t}\right) =\left( f\otimes f\right) \left( \sum \left( -1\right)
^{q_{1}p_{2}}e_{p_{1},q_{1}}\overline{t}\otimes e_{p_{2},q_{2}}\overline{t}%
\right) \\
&=&\sum \left( -1\right) ^{q_{1}p_{2}}e_{p_{1},q_{1}+p_{1}}\sum
i^{p^{2}}e_{p,q}\overline{t}\otimes e_{p_{2},q_{2}+p_{2}}\sum
i^{p^{2}}e_{p,q}\overline{t} \\
&=&\sum \left( -1\right) ^{\left( q_{1}+p_{1}\right)
p_{2}}i^{p_{1}^{2}+p_{2}^{2}}e_{p_{1},q_{1}}\overline{t}\otimes
e_{p_{2},q_{2}}\overline{t} \\
&=&\sum \left( -1\right)
^{q_{1}p_{2}}i^{p_{1}^{2}+2p_{1}p_{2}+p_{2}^{2}}e_{p_{1},q_{1}}\overline{t}%
\otimes e_{p_{2},q_{2}}\overline{t} \\
&=&\sum \left( -1\right) ^{q_{1}p_{2}}i^{\left( p_{1}+p_{2}\right)
^{2}}e_{p_{1},q_{1}}\overline{t}\otimes e_{p_{2},q_{2}}\overline{t} \\
&=&\left( \sum i^{\left( p_{1}+p_{2}\right) ^{2}}e_{p_{1},q_{1}}\otimes
e_{p_{2},q_{2}}\right) \left( \sum \left( -1\right)
^{s_{1}r_{2}}e_{r_{1},s_{1}}\overline{t}\otimes e_{r_{2},s_{2}}\overline{t}%
\right) \\
&=&\Delta \left( \sum i^{p^{2}}e_{p,q}\right) \Delta \left( \overline{t}%
\right) =\Delta \left( \sum i^{p^{2}}e_{p,q}\overline{t}\right) =\Delta
\left( f\left( \overline{t}\right) \right)
\end{eqnarray*}
and such an $f$ is a Hopf algebra isomorphism between $H_{0,0}$ and $H_{1,0}$%
. There are at most three nonisomorphic Hopf algebras of this kind:

\begin{enumerate}
\item  $H_{b:1}=H_{0,0}$ with the trivial cocycle and $\mathbf{G}\left(
H_{b:1}^{\ast }\right) =\left\langle \chi \right\rangle \times
\left\langle \varphi \right\rangle \times \left\langle \psi \right\rangle
\cong C_{2}\times C_{2}\times C_{2}$, where $\chi (x)=-1$, $\chi (y)=\chi
(t)=1$, $\varphi \left( x\right) =\varphi \left( y\right) =1$, $\varphi
\left( t\right) =-1$, $\psi \left( y\right) =-1$, $\psi \left( x\right)
=\psi \left( t\right) =1$. There is a degree $2$ irreducible representation
defined by
$\newline
\begin{tabular}{lll}
$\pi (x)=\left( 
\begin{array}{rr}
i & 0 \\ 
0 & -i
\end{array}
\right) $ & $\pi (y)=\left( 
\begin{array}{rr}
1 & 0 \\ 
0 & 1
\end{array}
\right) $ & $\pi (t)=\left( 
\begin{array}{rr}
0 & \left( -1\right) ^{p} \\ 
1 & 0
\end{array}
\right) $%
\end{tabular}
\newline $
with the property $\pi ^{2}=\pi \bullet \pi =1+\chi +\varphi +\chi \varphi $.
\smallskip

\item  $H_{b:y}=H_{0,1}$ with the cocycle defined by $\sigma \left(
t,t\right) =y$ and $\mathbf{G}\left( H_{b:y}^{\ast }\right) =\left\langle
\chi \right\rangle \times \left\langle \varphi \right\rangle \cong
C_{4}\times C_{2}$, where $\chi (x)=1$, $\chi (y)=-1$, $\chi (t)=i$, $%
\varphi \left( x\right) =-1$, $\varphi \left( y\right) =\varphi \left(
t\right) =1$. There is a degree $2$ irreducible representation defined by
$\newline
\begin{tabular}{lll}
$\pi (x)=\left( 
\begin{array}{rr}
i & 0 \\ 
0 & -i
\end{array}
\right) $ & $\pi (y)=\left( 
\begin{array}{rr}
1 & 0 \\ 
0 & 1
\end{array}
\right) $ & $\pi (t)=\left( 
\begin{array}{rr}
0 & \left( -1\right) ^{p} \\ 
1 & 0
\end{array}
\right) $%
\end{tabular}
\newline$
with the property $\pi ^{2}=\pi \bullet \pi =1+\chi ^{2}+\varphi +\chi
^{2}\varphi $.
\smallskip

\item  $H_{b:x^{2}y}=H_{1,1}$ with the cocycle defined by $\sigma \left(
t,t\right) =x^{2}y$ and $\mathbf{G}\left( H_{b:x^{2}y}^{\ast }\right)
=\left\langle \chi \right\rangle \times \left\langle \varphi \right\rangle
\cong C_{4}\times C_{2}$, where $\chi (x)=1$, $\chi (y)=-1$, $\chi (t)=i$, $%
\varphi \left( x\right) =-1$, $\varphi \left( y\right) =\varphi \left(
t\right) =1$. There is a degree $2$ irreducible representation defined by
$\newline
\begin{tabular}{lll}
$\pi (x)=\left( 
\begin{array}{rr}
i & 0 \\ 
0 & -i
\end{array}
\right) $ & $\pi (y)=\left( 
\begin{array}{rr}
1 & 0 \\ 
0 & 1
\end{array}
\right) $ & $\pi (t)=\left( 
\begin{array}{rr}
0 & \left( -1\right) ^{p} \\ 
1 & 0
\end{array}
\right) $%
\end{tabular}
\newline $
with the property $\pi ^{2}=\pi \bullet \pi =1+\chi ^{2}+\varphi +\chi
^{2}\varphi $.
\end{enumerate}

\subsubsection*{Case c)}

The action is given by 
\begin{eqnarray*}
t &\rightharpoonup &x=x \\
t &\rightharpoonup &y=x^{2}y
\end{eqnarray*}
Then $t\rightharpoonup e_{i,j}=e_{i,j+i}$. Write $v=\sigma \left( t,t\right)
=\sum c_{i,j}e_{i,j}$ . By $\left( \ref{coc1}\right) $  and
$\left( \ref{coc4}\right) $
\begin{eqnarray}
c_{i,j}&=&c_{i,i+j}  \label{cond7}\\
c_{i+p,j+q}&=&\left( -1\right) ^{ip}c_{i,j}c_{p,q}  \label{cond8}
\end{eqnarray}
Conditions $\left( \ref{cond7}\right) $ and $\left( \ref{cond8}\right) $
imply that $c_{1,0}^{4}=c_{0,1}=1$, $c_{2,0}=-c_{1,0}^{2}$ Thus $%
c_{1,0}=i^{k}$ for $k=0,1,2,3$ and 
\begin{equation*}
\sigma _{k}\left( t,t\right) =\sum _{p,q} \left( -1\right)^
{\frac{p\left( p-1\right) }{2}}i ^{kp}e_{p,q}
=x^{1-k}\left( \frac{1+i}{2}1+ \frac{1-i}{2}x^{2}\right)
\end{equation*}

For $k=0,1,2,3$ let $H_{k}$ be the Hopf algebras with the structures
described above with cocycles $\sigma _{k}$. Define 
\begin{eqnarray*}
f :H_{k+2}&\rightarrow &H_{k}\qquad \text{by} \\
f\left( e_{p,q}\right) &=&e_{p,q} \\
f\left( \overline{t}\right) &=&\sum_{p,q}\left( -1\right) ^{p+q}e_{p,q}%
\overline{t}=y\overline{t}
\end{eqnarray*}
and extend it multiplicatively to $f\left( e_{p,q}\overline{t}\right) $.
Then 
\begin{eqnarray*}
f\left( \overline{t}\right) f\left( \overline{t}\right) &=&
y\overline{t}y\overline{t}=x^{2} \overline{t}^{2}=
x^{2}x^{1-\left( k-2\right) }\left( \frac{1+i}{2}+ \frac{1-i}{2}x^{2}\right)\\
&=&x^{1-k}\left( \frac{1+i}{2}+ \frac{1-i}{2}x^{2}\right) 
=f\left( \overline{t}^{2}\right)\\
f\left( \overline{t}y\right) &=&f\left( x^{2}y
\overline{t}\right) =x^{2}yy\overline{t}=y\overline{t}y=
f\left( \overline{t}\right)f\left(y\right)\\
\left( f\otimes f\right) \Delta \left( \overline{t}\right) &=&\left(
f\otimes f\right) \left( \theta \left( t\right) \overline{t}\otimes 
\overline{t}\right) =\theta \left( t\right) \left( f\left( \overline{t}%
\right) \otimes f\left( \overline{t}\right) \right) =\theta \left( t\right)
\left( y\overline{t}\otimes y\overline{t}\right) \\
&=&\left( y\otimes y\right) \theta \left( t\right) \left( 
\overline{t}\otimes \overline{t}\right) =\Delta \left( y\right) \Delta
\left( \overline{t}\right) =\Delta \left( y\overline{t}\right) =\Delta
\left( f\left( \overline{t}\right) \right)
\end{eqnarray*}
and such a $f$ is a Hopf algebra isomorphism between $H_{k+2}$ and $%
H_{k}$. Thus there are exactly $2$ nonisomorphic Hopf algebras of
this type:

\begin{enumerate}
\item  $H_{c:\sigma _{0}}=H_{0}$ with cocycle $\sigma _{0}$
defined by $\sigma \left(
t,t\right) =\frac{1+i}{2}x+\frac{1-i}{2}x^{3}$ and $\mathbf{G}%
\left( H_{c:\sigma _{0}}^{\ast }\right) 
=\left\langle \chi \right\rangle \times
\left\langle \varphi \right\rangle \cong C_{4}\times C_{2}$, where $\chi
(x)=-1$, $\chi (y)=1$, $\chi (t)=i$, $\varphi \left( y\right) =-1$, $\varphi
\left( x\right) =\varphi \left( t\right) =1$. There is a degree $2$
irreducible representation defined by
$\newline
\begin{tabular}{lll}
$\pi (x)=\left( 
\begin{array}{rr}
i & 0 \\ 
0 & i
\end{array}
\right) $ & $\pi (y)=\left( 
\begin{array}{rr}
0 & 1 \\ 
1 & 0
\end{array}
\right) $ & $\pi (t)=\left( 
\begin{array}{rr}
i & 0 \\ 
0 & -i
\end{array}
\right) $%
\end{tabular}
\newline $
with the property $\pi ^{2}=\pi \bullet \pi =\chi +\chi ^{3}+\chi \varphi
+\chi ^{3}\varphi $.
\smallskip

\item  $H_{c:\sigma _{1}}=H_{1}$ with cocycle $\sigma _{1}$ 
defined by $\sigma\left( t,t\right) =\frac{1+i}{2}1+\frac{1-i}{2}x^{2}$ and 
$\mathbf{G}\left( H_{c:\sigma _{1}}^{\ast }\right) 
=\left\langle \chi \right\rangle \times
\left\langle \varphi \right\rangle \times \left\langle \psi \right\rangle
\cong C_{2}\times C_{2}\times C_{2}$, where $\chi (y)=-1$, $\chi (x)=\chi
(t)=1$, $\varphi \left( x\right) =\varphi \left( y\right) =1$, $\varphi
\left( t\right) =-1$, $\psi \left( x\right) =-1$, $\psi \left( y\right)
=\psi \left( t\right) =1$. There is a degree $2$ irreducible representation
defined by
$\newline
\begin{tabular}{lll}
$\pi (x)=\left( 
\begin{array}{rr}
i & 0 \\ 
0 & i
\end{array}
\right) $ & $\pi (y)=\left( 
\begin{array}{rr}
0 & 1 \\ 
1 & 0
\end{array}
\right) $ & $\pi (t)=\left( 
\begin{array}{rr}
\omega  & 0 \\ 
0 & -\omega 
\end{array}
\right) $%
\end{tabular}
\newline $
where $\omega $ is a primitive $8^{th}$-root of unity, with the property $%
\pi ^{2}=\pi \bullet \pi =\psi +\chi \psi +\varphi \psi +\chi \varphi \psi $.
\end{enumerate}

\subsection{Case of $\mathbf{G}\left( H\right) =C_{2}\times C_{2}\times
C_{2} $.}\label{3.2}
\ \newline
We will show that there are at most $4$ possible Hopf algebras of this kind.
Let $H$ be a nontrivial semisimple Hopf algebra of dimension $16$ with a
subHopfalgebra $K=k\left( C_{2}\times C_{2}\times C_{2}\right) ^{\ast } 
\cong k\left( C_{2}\times C_{2}\times C_{2}\right) $. Then
$\mathbf{G}\left( H\right) =G \cong C_{2}\times C_{2}\times C_{2}$.

Let $\mathbf{G}\left( H\right) =\left\langle x\right\rangle \times \left\langle
y\right\rangle \times \left\langle z\right\rangle $, where $\left| x\right|
=\left| y\right| = \left| z\right| =2$. Then the dual basis of $K\cong K^{\ast
}$ is given by 
\begin{equation*}
e_{p,q,r}=1/8\left( 1+\left( -1\right) ^{p}x\right) \left( 1+\left(
-1\right) ^{q}y\right) \left( 1+\left( -1\right) ^{r}z\right) ,\text{\qquad }%
p,q,r=0,1
\end{equation*}
Then 
\begin{eqnarray*}
\Delta _{H}\left( e_{p,q,r}\right) &=&\Delta _{K}\left( e_{p,q,r}\right)
=\sum_{\substack{ p_{1}+p_{2}\equiv p\mod 2  \\ %
q_{1}+q_{2}\equiv q\mod 2  \\ r_{1}+r_{2}\equiv r
\mod 2 }}e_{p_{1},q_{1},r_{1}}\otimes e_{p_{2},q_{2},r_{2}} \\
\Delta _{H}\left( \overline{t}\right) &=&\theta \left( t\right) \overline{t}%
\otimes \overline{t}
\end{eqnarray*}
where $\overline{t}=1\#t.$ Dualizing $\left( \ref{ext1}\right) $ we get
another extension 
\begin{equation*}
F^{\ast }\overset{\pi ^{\ast }}{\hookrightarrow }H^{\ast }\overset{i^{\ast }%
}{\twoheadrightarrow }K^{\ast }
\end{equation*}
and as in \cite[2.4]{Ma1}, \cite[2.11]{Ma2} or \cite[2.1]{Ma5}, since $k$ is
algebraically closed, there exist units $\overline{x}$, $\overline{y}$ and $%
\overline{z}\in H^{\ast }$, such that $\overline{x}^{2}=\overline{y}^{2}=%
\overline{z}^{2}=1_{H^{\ast }}$, $\left\langle e_{p,q,r},\overline{x}^{i}%
\overline{y}^{j}\overline{z}^{k}\right\rangle =\delta _{ip}\delta
_{jq}\delta _{kr}$ and $\alpha =\overline{z}^{-1}\overline{y}^{-1}\overline{z%
}\overline{y},\beta =\overline{z}^{-1}\overline{x}^{-1}\overline{z}\overline{%
x},\gamma =\overline{y}^{-1}\overline{x}^{-1}\overline{y}\overline{x}\in
F^{\ast }=k\left\{ e_{0},e_{1}\right\} $, where $\left\{ e_{r}\right\} $ is
a dual basis of $\left\{ t^{r}\right\} $, $r=0,1$. $\varepsilon \left(
\alpha \right) =\varepsilon \left( \beta \right) =\varepsilon \left( \gamma
\right) =1$ and therefore $\alpha =e_{0}+\xi _{3}e_{1},\beta =e_{0}+\xi
_{2}e_{1},\gamma =e_{0}+\xi _{1}e_{1}$. The right action $\rho ^{\ast
}:F^{\ast }\otimes K^{\ast }\rightarrow F^{\ast }$ is trivial, thus $F^{\ast
}$ lies in the center of $H^{\ast }$. 
\begin{equation*}
\overline{x}=\overline{y}^{2}\overline{x}=\overline{y}\overline{x}\overline{y%
}\gamma =\overline{x}\overline{y}\gamma \overline{y}\gamma =\overline{x}%
\overline{y}^{2}\gamma ^{2}=\overline{x}\gamma ^{2}
\end{equation*}
Thus $\gamma ^{2}=1$ and similarly $\alpha ^{2}=\beta ^{2}=1$. Therefore $%
\xi _{1},\xi _{2},\xi _{3}=\pm 1$ and, since $H^{\ast }$ is non-commutative,
they cannot be all equal to 1. 
\begin{eqnarray*}
\left\langle \Delta _{H}\left( \overline{t}\right) ,\overline{x}^{i}%
\overline{y}^{j}\overline{z}^{k}e_{l}\otimes \overline{x}^{p}\overline{y}^{q}%
\overline{z}^{r}e_{s}\right\rangle &=&\left\langle \overline{t},\overline{x}%
^{i}\overline{y}^{j}\overline{z}^{k}e_{l}\overline{x}^{p}\overline{y}^{q}%
\overline{z}^{r}e_{s}\right\rangle \\
&=&\delta _{ls}\left\langle \overline{t},\overline{x}^{i+p}\overline{y}^{j+q}%
\overline{z}^{k+r}\alpha ^{kq}\beta ^{kp}\gamma ^{jp}e_{l}\right\rangle =\xi
_{1}^{jp}\xi _{2}^{kp}\xi _{3}^{kq}\delta _{l1}\delta _{s1}
\end{eqnarray*}
On the other hand 
\begin{eqnarray*}
\left\langle \Delta _{H}\left( \overline{t}\right) ,\overline{x}^{i}%
\overline{y}^{j}\overline{z}^{k}e_{l}\otimes \overline{x}^{p}\overline{y}^{q}%
\overline{z}^{r}e_{s}\right\rangle &=&\left\langle \theta \left( t\right) 
\overline{t}\otimes \overline{t},\overline{x}^{i}\overline{y}^{j}\overline{z}%
^{k}e_{l}\otimes \overline{x}^{p}\overline{y}^{q}\overline{z}%
^{r}e_{s}\right\rangle \\
&=&\left\langle \theta \left( t\right) ,\overline{x}^{i}\overline{y}^{j}%
\overline{z}^{k}\otimes \overline{x}^{p}\overline{y}^{q}\overline{z}%
^{r}\right\rangle \delta _{k1}\delta _{r1}
\end{eqnarray*}
Therefore 
\begin{equation*}
\theta \left( t\right) =\sum_{ijkpqr}\xi _{1}^{jp}\xi _{2}^{kp}\xi
_{3}^{kq}e_{i,j,k}\otimes e_{p,q,r}
\end{equation*}

Action by $t$  is a Hopf algebra map and therefore $t\rightharpoonup G=G$
and $f_{t}:G\rightarrow G$ defined by $f_{t}\left( g\right)
=t\rightharpoonup g$ is a group automorphism of order $2$. Then, without
loss of generality there is only one possibility for such an automorphism:

\begin{eqnarray*}
t &\rightharpoonup &x=y \\
t &\rightharpoonup &y=x \\
t &\rightharpoonup &z=z
\end{eqnarray*}
Then $t\rightharpoonup e_{i,j,k}=e_{j,i,k}$.

Write $%
v=\sigma \left( t,t\right) =\sum c_{i,j,k}e_{i,j,k}$ then $%
c_{0,0,0}=\varepsilon \left( v\right) =1$ and $c_{i,j,k}\neq 0,$ since $v$
is a unit. By formula $\left( \ref{coc1}\right) $ 
\begin{equation}
c_{i,j,k}=c_{j,i,k}  \label{cond1}
\end{equation}
For $H$ to be a bialgebra we need $\Delta _{H}\left( \overline{t}^{2}\right)
=\Delta _{H}\left( \overline{t}\right) \Delta _{H}\left( \overline{t}\right) 
$ 
\begin{equation*}
\Delta _{H}\left( \overline{t}^{2}\right) =\Delta _{H}\left( v\right)
=\Delta _{K}\left( \sum c_{i,j,k}e_{i,j,k}\right) =\sum
c_{i+p,j+q,k+r}e_{i,j,k}\otimes e_{p,q,r}
\end{equation*}
On the other hand, 
\begin{eqnarray*}
&&\Delta _{H}\left( \overline{t}\right) \Delta _{H}\left( \overline{t}%
\right) =\left( \theta \left( t\right) \overline{t}\otimes \overline{t}%
\right) \left( \theta \left( t\right) \overline{t}\otimes \overline{t}\right)
\\
&=&\sum_{ijpq}\xi _{1}^{jp}\xi _{2}^{kp}\xi _{3}^{kq}e_{ijk}\otimes
e_{pqr} \left( \sum_{ijpq}\xi _{1}^{jp}\xi _{2}^{kp}\xi
_{3}^{kq}\left( t\rightharpoonup e_{ijk}\right) \otimes \left(
t\rightharpoonup e_{pqr}\right) \right) \sigma \left( t,t\right) \otimes
\sigma \left( t,t\right) \\
&=&\sum \xi _{1}^{jp+iq}\xi _{2}^{kp+kq}\xi
_{3}^{kq+kp}c_{i,j,k}c_{p,q,r}e_{i,j,k}\otimes e_{p,q,r}
\end{eqnarray*}
Therefore 
\begin{equation}
c_{i+p,j+q,k+r}=\xi _{1}^{jp+iq}\xi _{2}^{kp+kq}\xi
_{3}^{kq+kp}c_{i,j,k}c_{p,q,r}  \label{cond2}
\end{equation}
Conditions $\left( \ref{cond1}\right) $ and $\left( \ref{cond2}\right) $
imply that $%
c_{1,0,0}^{2}=c_{0,1,0}^{2}=c_{0,0,1}^{2}=c_{1,1,0}^{2}=c_{1,1,1}^{2}=1$ and 
$c_{0,1,0}=c_{1,0,0}$ 
\begin{eqnarray*}
c_{1,1,0} &=&\xi _{1}c_{0,1,0}c_{1,0,0}=\xi _{1}c_{1,0,0}^{2}=\xi _{1} \\
c_{1,0,1} &=&c_{1,0,0}c_{0,0,1} \\
c_{1,0,1} &=&c_{0,0,1}c_{1,0,0}\xi _{2}\xi _{3}
\end{eqnarray*}
Thus $\xi _{2}\xi _{3}=1$, that is $\xi _{2}=\xi _{3}$ and $%
c_{1,0,0}=c_{0,1,0}=\omega =\pm 1$ and $c_{0,0,1}=\tau =\pm 1$ and 
\begin{equation*}
\theta \left( t\right) =\sum_{ijkpqr}\xi _{1}^{jp}\xi
_{2}^{kp+kq}e_{i,j,k}\otimes e_{p,q,r}
\end{equation*}
\begin{eqnarray*}
\sigma \left( t,t\right) &=&e_{0,0,0}+\tau e_{0,0,1}+\xi _{1}e_{1,1,0}+\xi
_{1}\tau e_{1,1,1}+\omega \left( e_{1,0,0}+e_{0,1,0}+\tau e_{1,0,1}+\tau
e_{0,1,1}\right) \\
&=&\sum \omega ^{p+q}e_{pqr}\sum \xi _{1}^{pq}e_{pqr}\sum \tau
^{r}e_{pqr}=\left( xy\right) ^{\delta _{\omega ,-1}}\left( 1+x+y-xy\right)
^{\delta _{\xi _{1},-1}}z^{\delta _{\tau ,-1}}
\end{eqnarray*}
For $\xi _{1},\xi _{2},\tau ,\omega =\pm 1$ let $H_{d:\xi _{1},\xi _{2},\tau
,\omega }$ be the Hopf algebras with the structures described above with
cocycles $\sigma _{\xi _{1},\tau ,\omega }$. Then $\sigma _{\xi _{1},\tau
,-1}\left( t,t\right) =xy\sigma _{\xi _{1},\tau ,1}\left( t,t\right) $.
Define 
\begin{eqnarray*}
f:H_{d:\xi _{1},\xi _{2},\tau ,-1} &\rightarrow &H_{d:\xi _{1},\xi _{2},\tau
,1}\qquad \text{by} \\
f\left( e_{p,q,r}\right) &=&e_{p,q,r} \\
f\left( \overline{t}\right) &=&x\overline{t}
\end{eqnarray*}
and extend it multiplicatively to $f\left( e_{p,q,r}\overline{t}\right) $.
Then 
\begin{eqnarray*}
f\left( \overline{t}\right) f\left( \overline{t}\right) &=&x\overline{t}x%
\overline{t}=xy\overline{t}^{2}=xy\sigma _{\xi _{1},\tau ,1}\left(
t,t\right)\\
&=&\sigma _{\xi _{1},\tau ,-1}\left( t,t\right) =f\left( \sigma
_{\xi _{1},\tau ,-1}\left( t,t\right) \right) =f\left( \overline{t}%
^{2}\right) \\
f\left( \overline{t}x\right) &=&f\left( y\overline{t}\right) =yx\overline{t}%
=x\overline{t}x=f\left( \overline{t}\right) f\left( x\right) \\
f\left( \overline{t}y\right) &=&f\left( x\overline{t}\right) =x^{2}\overline{%
t}=x\overline{t}y=f\left( \overline{t}\right) f\left( y\right)\\
\left( f\circ f\right) \Delta \left( \overline{t}\right) &=&\left( f\circ
f\right) \left( \theta \left( t\right) \overline{t}\otimes \overline{t}%
\right) =\theta \left( t\right) \left( f\left( \overline{t}\right) \otimes
f\left( \overline{t}\right) \right) =\theta \left( t\right) \left( x%
\overline{t}\otimes x\overline{t}\right) \\
&=&\left( x\otimes x\right) \theta \left( t\right) \left( \overline{t}%
\otimes \overline{t}\right) =\Delta \left( x\right) \Delta \left( \overline{t%
}\right) =\Delta \left( x\overline{t}\right) =\Delta \left( f\left( 
\overline{t}\right) \right)
\end{eqnarray*}
and such a $f$ is a Hopf algebra isomorphism between $H_{d:\xi _{1},\xi
_{2},\tau ,-1}$ and $H_{d:\xi _{1},\xi _{2},\tau ,1}$. Define 
\begin{eqnarray*}
f^{\prime }:H_{d:-1,-1,\tau ,1} &\rightarrow &H_{d:-1,1,\tau ,1}\qquad \text{%
by} \\
f^{\prime }\left( e_{p,q,r}\right) &=&e_{p+r,q+r,r} \\
f^{\prime }\left( \overline{t}\right) &=&%
{\frac12}
\left( 1+z+iy-iyz\right) \overline{t}=\sum i^{r^{2}}\left( -1\right)
^{qr}e_{p,q,r}\overline{t}
\end{eqnarray*}
and extend it multiplicatively to $f^{\prime }\left( e_{p,q,r}\overline{t}%
\right) $. Then $f^{\prime }\mid _{\mathbf{G}\left( H_{d:-1,-1,\tau
,1}\right) }$ is a group isomorphism $\mathbf{G}\left( H_{d:-1,-1,\tau
,1}\right) \rightarrow \mathbf{G}\left( H_{d:-1,1,\tau ,1}\right) $ with $%
f^{\prime }\left( x\right) =xz$, $f^{\prime }\left( y\right) =yz$, $%
f^{\prime }\left( z\right) =z$ and 
\begin{eqnarray*}
f^{\prime }\left( \overline{t}\right) f^{\prime }\left( \overline{t}\right)
&=&
{\frac14}
\left( \left( 1+z\right) +iy\left( 1-z\right) \right) \overline{t}\left(
\left( 1+z\right) +iy\left( 1-z\right) \right) \overline{t} \\
&=&
{\frac14}
\left( \left( 1+z\right) +iy\left( 1-z\right) \right) \left( \left(
1+z\right) +ix\left( 1-z\right) \right) \overline{t}^{2} \\
&=&\frac{1}{8}\left( 2+2z-xy\left( 2-2z\right) \right) \left(
1+x+y-xy\right) z^{\delta _{\tau ,-1}} \\
&=&
{\frac14}
\left( \left( 1-xy\right) +z\left( 1+xy\right) \right) \left( \left(
1-xy\right) +x\left( 1+xy\right) \right) z^{\delta _{\tau ,-1}} \\
&=&
{\frac14}
\left( 2-2xy+xz\left( 2+2xy\right) \right) z^{\delta _{\tau ,-1}}=
{\frac12}
\left( 1+xz+yz-xy\right) z^{\delta _{\tau ,-1}} \\
&=&f\left( 
{\frac12}
\left( 1+x+y-xy\right) z^{\delta _{\tau ,-1}}\right) =f^{\prime }\left( 
\overline{t}^{2}\right) \\
f^{\prime }\left( \overline{t}x\right) &=&f^{\prime }\left( y\overline{t}%
\right) =
{\frac12}
yz\left( 1+z+iy-iyz\right) \overline{t}\\
&=&
{\frac12}
\left( 1+z+iy-iyz\right) \overline{t}xz=f^{\prime }\left( \overline{t}%
\right) f^{\prime }\left( x\right) \\
f^{\prime }\left( \overline{t}y\right) &=&f^{\prime }\left( x\overline{t}%
\right) =
{\frac12}
xz\left( 1+z+iy-iyz\right) \overline{t}\\
&=&
{\frac12}
\left( 1+z+iy-iyz\right) \overline{t}yz=f^{\prime }\left( \overline{t}%
\right) f^{\prime }\left( y\right)
\end{eqnarray*}
\begin{eqnarray*}
\left( f^{\prime }\otimes f^{\prime }\right) \Delta \left( \overline{t}%
\right) &=&\left( f^{\prime }\otimes f^{\prime }\right) \left( \theta \left(
t\right) \overline{t}\otimes \overline{t}\right) \\
&=&\left( f^{\prime }\otimes
f^{\prime }\right) \left( \sum \left( -1\right) ^{bp}\left( -1\right)
^{cp+cq}e_{a,b,c}\overline{t}\otimes e_{p,q,r}\overline{t}\right) \\
&=&\left( \sum \left( -1\right) ^{bp}\left( -1\right) ^{c\left( p+q\right)
}e_{a+c,b+c,c}\otimes e_{p+r,q+r,r}\right) \\
&&\times \left( \sum i^{n^{2}}\left( -1\right) ^{mn}e_{l,m,n}\overline{t}%
\otimes \sum i^{n^{2}}\left( -1\right) ^{mn}e_{l,m,n}\overline{t}\right) \\
&=&\left( \sum \left( -1\right) ^{\left( b+c\right) \left( p+r\right)
}\left( -1\right) ^{c\left( p+q\right) }e_{a,b,c}\otimes e_{p,q,r}\right) \\
&&\times \left( \sum i^{n^{2}}\left( -1\right) ^{mn}e_{l,m,n}\overline{t}%
\otimes \sum i^{n^{2}}\left( -1\right) ^{mn}e_{l,m,n}\overline{t}\right) \\
&=&\sum \left( -1\right) ^{bp+cp+br+cr}\left( -1\right)
^{cp+cq}i^{c^{2}}\left( -1\right) ^{bc}i^{r^{2}}\left( -1\right)
^{qr}e_{a,b,c}\overline{t}\otimes e_{p,q,r}\overline{t} \\
&=&\sum \left( -1\right) ^{bp+br+cq+bc+qr}\left( -1\right)
^{cr}i^{c^{2}}i^{r^{2}}e_{a,b,c}\overline{t}\otimes e_{p,q,r}\overline{t} \\
&=&\sum \left( -1\right) ^{bp+br+cq+bc+qr}\left( -1\right)
^{cr}i^{c^{2}}i^{r^{2}}e_{a,b,c}\overline{t}\otimes e_{p,q,r}\overline{t} \\
&=&\sum i^{\left( c+r\right) ^{2}}\left( -1\right) ^{\left( b+q\right)
\left( c+r\right) }\left( -1\right) ^{bp}e_{a,b,c}\overline{t}\otimes
e_{p,q,r}\overline{t} \\
&=&\sum i^{n^{2}}\left( -1\right) ^{mn}\sum_{\substack{ l_{1}+l_{2}=l  \\ %
m_{1}+m_{2}=m  \\ n_{1}+n_{2}=n}}e_{l_{1},m_{1},n_{1}}\otimes
e_{l_{2},m_{2},n_{2}} \\
&&\times \sum \left( -1\right) ^{bp}e_{a,b,c}\overline{t}\otimes e_{p,q,r}%
\overline{t} \\
&=&\sum i^{n^{2}}\left( -1\right) ^{mn}\Delta \left( e_{l,m,n}\right) \Delta
\left( \overline{t}\right) \\
&=&\Delta \left( \sum i^{n^{2}}\left( -1\right)
^{mn}e_{l,m,n}\overline{t}\right) =\Delta f^{\prime }\left( \overline{t}%
\right)
\end{eqnarray*}
and such an $f^{\prime }$ is a Hopf algebra isomorphism between $%
H_{d:-1,-1,\tau ,1}$ and $H_{d:-1,1,\tau ,1}$. Thus we may assume that $%
\omega =1$ and $\xi _{2}=-1$. Therefore there are at most four nonisomorphic
Hopf algebras $H_{d:\xi _{1},\tau }$ of this kind, $H_{d:1,1}$, $H_{d:1,-1}$%
, $H_{d:-1,1}$ and $H_{d:-1,-1}$:

\begin{enumerate}
\item  $H_{d:1,1}$ with the trivial cocycle and $\mathbf{G}\left(
H_{d:1,1}^{\ast }\right) =\left\langle \chi \right\rangle \times
\left\langle \varphi \right\rangle \times \left\langle \psi \right\rangle
\cong C_{2}\times C_{2}\times C_{2}$, where $\chi \left( x\right) =\chi
\left( y\right) =\chi \left( z\right) =1$, $\chi \left( t\right) =-1$, $%
\varphi (x)=\varphi (y)=-1$, $\varphi (z)=\varphi (t)=1$, $\psi \left(
z\right) =-1$, $\psi \left( x\right) =\psi \left( y\right) =\psi \left(
t\right) =1$. There is a degree $2$ irreducible representation defined by
\newline
\noindent 
\begin{tabular}{llll}
$\pi (x)=\left( 
\begin{array}{rr}
1 & 0 \\ 
0 & -1
\end{array}
\right) $ & $\pi (y)=\left( 
\begin{array}{rr}
-1 & 0 \\ 
0 & 1
\end{array}
\right) $ & $\pi (z)=\left( 
\begin{array}{rr}
1 & 0 \\ 
0 & 1
\end{array}
\right) $ & $\pi (t)=\left( 
\begin{array}{rr}
0 & 1 \\ 
1 & 0
\end{array}
\right) $%
\end{tabular}
\newline
with the property $\pi ^{2}=\pi \bullet \pi =1+\chi +\varphi +\chi \varphi $%
.
\smallskip

\item  $H_{d:1,-1}$ with the cocycle defined by $\sigma \left(
t,t\right) =z$ and $\mathbf{G}\left( H_{d:1,-1}^{\ast }\right) =\left\langle
\chi \right\rangle \times \left\langle \varphi \right\rangle \cong
C_{4}\times C_{2}$, where $\chi (x)=\chi (y)=1$, $\chi (z)=-1$, $\chi (t)=i$%
, $\varphi \left( x\right) =\varphi \left( y\right) =-1$, $\varphi \left(
t\right) =\varphi \left( z\right) =1$. There is a degree $2$ irreducible
representation defined by
\newline
\noindent 
\begin{tabular}{llll}
$\pi (x)=\left( 
\begin{array}{rr}
1 & 0 \\ 
0 & -1
\end{array}
\right) $ & $\pi (y)=\left( 
\begin{array}{rr}
-1 & 0 \\ 
0 & 1
\end{array}
\right) $ & $\pi (z)=\left( 
\begin{array}{rr}
1 & 0 \\ 
0 & 1
\end{array}
\right) $ & $\pi (t)=\left( 
\begin{array}{rr}
0 & 1 \\ 
1 & 0
\end{array}
\right) $%
\end{tabular}
\newline
with the property $\pi ^{2}=\pi \bullet \pi =1+\chi ^{2}+\varphi +\chi
^{2}\varphi $.
\smallskip

\item  $H_{d:-1,1}$ with the cocycle defined by $\sigma \left(
t,t\right) =1+x+y-xy$ and $\mathbf{G}\left( H_{d:-1,1}^{\ast }\right)
=\left\langle \chi \right\rangle \times \left\langle \varphi \right\rangle
\times \left\langle \psi \right\rangle \cong C_{2}\times C_{2}\times C_{2}$,
where $\chi \left( x\right) =\chi \left( y\right) =\chi \left( z\right) =1$, 
$\chi \left( t\right) =-1$, $\varphi (x)=\varphi (y)=-1$, $\varphi (z)=1$, $%
\varphi (t)=i$, $\psi \left( z\right) =-1$, $\psi \left( x\right) =\psi
\left( y\right) =\psi \left( t\right) =1$. There is a degree $2$ irreducible
representation defined by
\newline
\noindent 
\begin{tabular}{llll}
$\pi (x)=\left( 
\begin{array}{rr}
1 & 0 \\ 
0 & -1
\end{array}
\right) $ & $\pi (y)=\left( 
\begin{array}{rr}
-1 & 0 \\ 
0 & 1
\end{array}
\right) $ & $\pi ( z) =\left( 
\begin{array}{rr}
1 & 0 \\ 
0 & 1
\end{array}
\right) $ & $\pi ( t) =\left( 
\begin{array}{rr}
0 & 1 \\ 
1 & 0
\end{array}
\right) $%
\end{tabular}
\newline
with the property $\pi ^{2}=\pi \bullet \pi =1+\chi +\varphi +\chi \varphi $%
.
\smallskip

\item  $H_{d:-1,-1}$ with the cocycle defined by $\sigma \left(
t,t\right) =\left( 1+x+y-xy\right) z$ and \newline$\mathbf{G}\left(
H_{d:-1,-1}^{\ast }\right)=\left\langle \chi \right\rangle \times
\left\langle \varphi \right\rangle \cong C_{4}\times C_{2}$, where $\chi
\left( x\right) =\chi \left( y\right) =1$, $\chi \left( z\right) =-1$, $\chi
\left( t\right) =i$, $\varphi (x)=\varphi (y)=-1$, $\varphi (z)=1$, $\varphi
(t)=i$. There is a degree $2$ irreducible representation defined by
\newline
\noindent 
\begin{tabular}{llll}
$\pi (x)=\left( 
\begin{array}{rr}
1 & 0 \\ 
0 & -1
\end{array}
\right) $ & $\pi (y)=\left( 
\begin{array}{rr}
-1 & 0 \\ 
0 & 1
\end{array}
\right) $ & $\pi (z) =\left( 
\begin{array}{rr}
1 & 0 \\ 
0 & 1
\end{array}
\right) $ & $\pi (t) =\left( 
\begin{array}{rr}
0 & 1 \\ 
1 & 0
\end{array}
\right) $%
\end{tabular}
\newline
with the property $\pi ^{2}=\pi \bullet \pi =1+\chi ^{2}+\varphi +\chi
^{2}\varphi $.
\end{enumerate}

\subsection{Case of $G=D_{8}$.}\label{3.3}
\ \newline Let $G=D_{8}=\left\langle x,y\mid 
x^{4}=y^{2}=1,yx=x^{-1}y\right\rangle $. Let $\left\{ e_{pq}\right\}
_{p=0,1,2,3;q=0,1}$ be the basis of $K$, dual to the basis $\left\{
x^{p}y^{q}\right\} _{p=0,1,2,3;q=0,1}$ of $K^{\ast }=kD_{8}$. Then 
\[
\Delta _{H}\left( e_{pq}\right) =\Delta _{K}\left( e_{pq}\right)
=\sum _{\substack{p_{1}+p_{2}+2q_{1}p_{2}\equiv p\mod 4 \\
q_{1}+q_{2}\equiv q\mod 2}}
e_{p_{1}q_{1}}\otimes e_{p_{2}q_{2}} 
\]
and it is easy to check that elements 
\begin{eqnarray*}
X &=&\sum_{pq}\left( -1\right) ^{p}e_{pq} \\
Y &=&\sum_{pq}\left( -1\right) ^{q}e_{pq}
\end{eqnarray*}
are grouplike of order $2$. For $\overline{t}=1\#t$

\[
\Delta _{H}\left( \overline{t}\right) =\theta \left( t\right) \overline{t}%
\otimes \overline{t} 
\]
Dualizing $\left( \ref{ext1}\right) $ we get another extension 
\[
F^{\ast }\stackrel{\pi ^{\ast }}{\hookrightarrow }H^{\ast }\stackrel{i^{\ast
}}{\twoheadrightarrow }K^{\ast } 
\]
and as in \cite[2.4]{Ma1}, \cite[2.11]{Ma2} or \cite[2.1]{Ma5}, since $k$ is
algebraically closed, there exist units $\overline{x}$ and $\overline{y}\in
H^{\ast }$, such that $\overline{x}^{4}=\overline{y}^{2}=1_{H^{\ast }}$, $%
\left\langle e_{pq},\overline{x}^{i}\overline{y}^{j}\right\rangle =\delta
_{ip}\delta _{jq}$ and $\alpha =\overline{y}\overline{x}^{2}\overline{y}%
\overline{x}^{2}\in F^{\ast }=k\left\{ e_{0},e_{1}\right\} $, where $\left\{
e_{r}\right\} $ is a dual basis of $\left\{ t^{r}\right\} $, $r=0,1$. The
right action $\rho ^{\ast }:F^{\ast }\otimes K^{\ast }\rightarrow F^{\ast }$
is trivial, thus $F^{\ast }$ lies in the center of $H^{\ast }$. 
\[
\overline{x}^{2}=\overline{y}^{2}\overline{x}^{2}=\overline{y}\overline{x}%
^{2}\overline{y}\alpha =\overline{x}^{2}\overline{y}\alpha \overline{y}%
\alpha =\overline{x}^{2}\overline{y}^{2}\alpha ^{2}=\overline{x}^{2}\alpha
^{2} 
\]
Thus $\alpha ^{2}=1$.

Consider $\beta =\overline{y}\overline{x}\overline{y}\overline{x}\in F^{\ast
}=k\left\{ e_{0},e_{1}\right\} $. $\varepsilon \left( \beta \right)
=\varepsilon \left( \overline{y}^{-1}\overline{x}^{-1}\overline{y}\overline{x%
}\right) =1$ and therefore $\beta =e_{0}+\xi e_{1}$. Moreover, $\overline{y}%
\overline{x}\overline{y}\overline{x}^{-1}=\beta \overline{x}^{2}$ and 
\[
\overline{x}=\overline{y}^{2}\overline{x}=\overline{y}\beta \overline{x}^{2}%
\overline{x}\overline{y}=\overline{y}\beta \overline{x}^{2}\overline{y}\beta 
\overline{x}^{2}\overline{x}=\overline{y}\overline{x}^{2}\overline{y}%
\overline{x}^{3}\beta ^{2}=\overline{y}^{2}\alpha \overline{x}^{2}\overline{x%
}^{3}\beta ^{2}=\overline{x}\alpha \beta ^{2} 
\]
Thus $\beta ^{2}=\alpha ^{-1}=\alpha $, implying $\beta ^{4}=1$ and $\xi
=\pm 1$ or $\pm i$. 
\begin{eqnarray*}
\left\langle \Delta _{H}\left( \overline{t}\right) ,\overline{x}^{i}%
\overline{y}^{j}e_{k}\otimes \overline{x}^{p}\overline{y}^{q}e_{r}\right%
\rangle &=&\left\langle \overline{t},\overline{x}^{i}\overline{y}^{j}e_{k}%
\overline{x}^{p}\overline{y}^{q}e_{r}\right\rangle =\delta _{kr}\left\langle 
\overline{t},\overline{x}^{i+p}\beta ^{jp}\overline{x}^{2jp}\overline{y}%
^{j+q}e_{k}\right\rangle \\
&=&\delta _{kr}\left\langle \overline{t},\overline{x}^{i+p+2jp}\overline{y}%
^{j+q}\beta ^{jp}e_{k}\right\rangle =\xi ^{jp}\delta _{k1}\delta _{r1}
\end{eqnarray*}
On the other hand 
\[
\left\langle \Delta _{H}\left( \overline{t}\right) ,\overline{x}^{i}%
\overline{y}^{j}e_{k}\otimes \overline{x}^{p}\overline{y}^{q}e_{r}\right%
\rangle =\left\langle \theta \left( t\right) \overline{t}\otimes \overline{t}%
,\overline{x}^{i}\overline{y}^{j}e_{k}\otimes \overline{x}^{p}\overline{y}%
^{q}e_{r}\right\rangle =\left\langle \theta \left( t\right) ,\overline{x}^{i}%
\overline{y}^{j}\otimes \overline{x}^{p}\overline{y}^{q}\right\rangle \delta
_{k1}\delta _{r1} 
\]
Therefore 
\[
\theta \left( t\right) =\sum_{ijpq}\xi ^{jp}e_{ij}\otimes e_{pq} 
\]

It is easy to check that if $\xi =\pm i$ then $1,X,Y$ and $XY$ are the only
grouplikes of $H$.

Write $v=\sigma \left( t,t\right) =\sum c_{i,j}e_{i,j}$ then $%
c_{0,0}=\varepsilon \left( v\right) =1$ and $c_{i,j}\neq 0,$ since $v$ is a
unit and 
\[
\Delta _{H}\left( \overline{t}^{2}\right) =\Delta _{H}\left( v\right)
=\Delta _{K}\left( \sum c_{i,j}e_{i,j}\right) =\sum
c_{p+r+2rq,q+s}e_{p,q}\otimes e_{r,s} 
\]
On the other hand, if we write 
\begin{equation*}
t\rightharpoonup e_{p,q} =e_{\alpha _{1}\left( p,q\right) ,\alpha _{2}\left(
p,q\right) }
\end{equation*}
\begin{eqnarray*}
\Delta \left( \overline{t}\right)\Delta \left( \overline{t}\right) &=&
\sum_{pqrs}\xi ^{qr}e_{p,q}\overline{t}\otimes e_{r,s}\overline{t}
\sum_{pqrs}\xi ^{qr}e_{p,q}\overline{t}\otimes e_{r,s}\overline{t}\\
&=& \sum_{pqrs}\xi ^{qr}e_{p,q}\otimes e_{r,s}
\sum_{pqrs}\xi ^{qr}
e_{\alpha _{1}\left( p,q\right) ,\alpha _{2}\left( p,q\right) }
\overline{t}^{2}\otimes 
e_{\alpha _{1}\left( r,s\right) ,\alpha _{2}\left( r,s\right) }
\overline{t}^{2}\\
&=& \sum_{pqrs}\xi ^{qr}e_{p,q}\otimes e_{r,s}
\sum_{pqrs}\xi ^{\alpha _{2}\left( p,q\right) \alpha _{1}\left( r,s\right) }
e_{p,q}\overline{t}^{2}\otimes e_{r,s}\overline{t}^{2}\\ 
&=&\sum_{pqrs}\xi ^{qr+\alpha _{2}\left( p,q\right) 
\alpha _{1}\left( r,s\right) }c_{pq}c_{rs}e_{p,q}\otimes e_{r,s}
\end{eqnarray*}
Thus for $H$ to be a bialgebra we should have
\begin{equation}
c_{p+r+2rq,q+s}=\xi ^{qr+\alpha _{2}\left( p,q\right) 
\alpha _{1}\left( r,s\right) }c_{pq}c_{rs} \label{Coc4}
\end{equation}

Action by $t$ is a Hopf algebra map and therefore it induces a group
automorphism $f_{t}:G\rightarrow G$ defined by $\left\langle
e_{p,q},f_{t}\left( g\right) \right\rangle =\left\langle t\rightharpoonup
e_{p,q},g\right\rangle $, which has order $2$.

$f_{t}\left( x\right) =x$ or $x^{-1}$ since order of $x$\ is $4$. If $%
f_{t}\left( x\right) =x$ then in order for $f_{t}$ to be of order $2$ we
should have $f_{t}\left( y\right) =x^{2}y$. If $f_{t}\left( x\right) =x^{-1}$
then renaming generators we are down to two choices for $f_{t}\left(
y\right) $, namely $f_{t}\left( y\right) =y$ or $xy$. Thus there are three
possibilities for the action of $t$; we consider them below:

\subsubsection*{Case A)}

The action is given by $t\rightharpoonup e_{p,q}=e_{p+2q,q}$, corresponding
to 
\begin{eqnarray*}
f_{t}\left( x\right) &=&x \\
f_{t}\left( y\right) &=&x^{2}y
\end{eqnarray*}
Then $X$ and $Y$ are central grouplikes of $H$. Write $v=\sigma \left(
t,t\right) =\sum c_{p,q}e_{p,q}$. By $\left( \ref{coc1}\right) $ and
$\left( \ref{Coc4}\right) $ 
\begin{eqnarray}
c_{p,q}&=&c_{p+2q,q}  \label{Cond3} \\
c_{p+r+2rq,q+s}&=&\xi ^{qr+q\left( r+2s\right) }c_{p,q}c_{r,s}
\xi ^{2q\left( r+s\right) }c_{p,q}c_{r,s}  \label{Cond4}
\end{eqnarray}
Conditions $\left( \ref{Cond3}\right) $ and $\left( \ref{Cond4}\right) $
imply that 
\begin{eqnarray*}
\xi ^{2}c_{0,1}c_{1,0} &=&c_{3,1}=c_{1,1}=c_{1,0}c_{0,1} \\
\xi ^{2}c_{0,1}c_{0,1} &=&c_{0,0}=1 \\
c_{1,0}c_{1,0} &=&c_{2,0} \\
c_{2,0}c_{0,1} &=&c_{2,1}=c_{0,1}
\end{eqnarray*}
Thus $\xi ^{2}=1$, $c_{2,0}=1$ and $c_{1,0}^{2}=c_{0,1}^{2}=1.$ Therefore 
$c_{1,0}=\left( -1\right) ^k$ and $c_{0,1}=\left( -1\right) ^l$ for 
$k,l=0,1$ and
\begin{equation}
\sigma \left( t,t\right)
=\sum\left( -1\right) ^{kp}\left( -1\right) ^{lq}e_{p,q}
=\sum\left( -1\right) ^{kp} e_{p,q}\sum\left( -1\right)  
^{ls}e_{r,s}=X^{k}Y^{l},\qquad k,l=0,1
\end{equation}

If $\xi =1$ then $\overline{t}$ is a grouplike of $H$, if $\xi =-1$ then $%
\sum i^{p}e_{p,q}\overline{t}$ is a grouplike of $H$. In both cases $\mathbf{%
G}\left( H\right) $ is abelian of order $8$ and $H$ was described in Section 
\ref{3.1} or Section \ref{3.2}.

\subsubsection*{Case B)}

The action is given by $t\rightharpoonup e_{p,q}=e_{-p,q}$, corresponding to 
\begin{eqnarray*}
f_{t}\left( x\right) &=&x^{-1} \\
f_{t}\left( y\right) &=&y
\end{eqnarray*}
Then $X$ and $Y$ are central grouplikes of $H$. Write $v=\sigma \left(
t,t\right) =\sum c_{p,q}e_{p,q}$ . By $\left( \ref{coc1}\right) $  and
$\left( \ref{Coc4}\right) $ 
\begin{eqnarray}
c_{p,q}&=&c_{-p,q}  \label{Cond5} \\
c_{p+r+2rq,q+s}&=&\xi ^{qr-qr}c_{p,q}c_{r,s}=c_{p,q}c_{r,s}  \label{Cond6}
\end{eqnarray}
Conditions $\left( \ref{Cond5}\right) $ and $\left( \ref{Cond6}\right) $
imply that $c_{1,0}=\left( -1\right) ^k$ and $c_{0,1}=\left( -1\right) ^l$ for 
$k,l=0,1$ and
\begin{equation}
\sigma \left( t,t\right)
=\sum\left( -1\right) ^{kp}\left( -1\right) ^{lq}e_{p,q}
=\sum\left( -1\right) ^{kp} e_{p,q}\sum\left( -1\right)  
^{ls}e_{r,s}=X^{k}Y^{l},\qquad k,l=0,1
\end{equation} 
If $\xi =1$ then $\overline{t}$ is a grouplike of $H$, if $\xi =-1$ then $%
\sum i^{p}e_{p,q}\overline{t}$ is a grouplike of $H$. In both cases $\mathbf{%
G}\left( H\right) $ is abelian of order $8$ and $H$ was described in Section 
\ref{3.1} or Section \ref{3.2}. So now we will consider only $\xi =\pm i$.

For $k,l=0,1$ let $H_{\xi ,X^{k}Y^{l}}$ be the Hopf algebras with the
structures described above with cocycles $\sigma _{k,l}\left( t,t\right)
=X^{k}Y^{l}$. Define 
\begin{eqnarray*}
f:H_{\xi ,X^{k}Y^{l}} &\rightarrow &H_{-\xi ,X^{k}Y^{l}}\qquad \text{by} \\
f\left( e_{r,s}\right) &=&e_{r,s} \\
f\left( \overline{t}\right) &=&\sum i^{p}e_{p,q}\overline{t}
\end{eqnarray*}
and extend it multiplicatively to $f\left( e_{r,s}\overline{t}\right) $.
Then 

\begin{eqnarray*}
f\left( \overline{t}\right) f\left( \overline{t}\right) &=&\sum i^{p}e_{p,q}%
\overline{t}\sum i^{p}e_{p,q}\overline{t}=\sum i^{p}e_{p,q}\sum i^{-p}e_{p,q}%
\overline{t}^{2}=\overline{t}^{2}=f\left( \overline{t}^{2}\right) \\
f\left( \overline{t}e_{r,s}\right) &=&f\left( e_{-r,s}\overline{t}\right)
=e_{-r,s}\sum i^{p}e_{p,q}\overline{t}=\sum i^{p}e_{p,q}\overline{t}%
e_{r,s}=f\left( \overline{t}\right) f\left( e_{r,s}\right) \\
\Delta \left( f\left( \overline{t}\right) \right) &=&\Delta \left( \sum
i^{p}e_{p,q}\overline{t}\right) =\Delta \left( \sum i^{p}e_{p,q}\right)
\Delta \left( \overline{t}\right) = \\
&=&\left( \sum i^{p_{1}+p_{2}+2q_{1}p_{2}}e_{p_{1},q_{1}}\otimes
e_{p_{2},q_{2}}\right) \left( \sum \xi ^{s_{1}r_{2}}e_{r_{1},s_{1}}\overline{%
t}\otimes e_{r_{2},s_{2}}\overline{t}\right) \\
&=&\left( \sum i^{p_{1}+p_{2}+2q_{1}p_{2}}\xi ^{q_{1}p_{2}}e_{p_{1},q_{1}}%
\overline{t}\otimes e_{p_{2},q_{2}}\overline{t}\right) \\
&=&\sum i^{p_{1}+p_{2}}\left( -\xi \right) ^{q_{1}p_{2}}e_{p_{1},q_{1}}%
\overline{t}\otimes e_{p_{2},q_{2}}\overline{t} \\
&=&\sum \left( -\xi \right) ^{q_{1}p_{2}}e_{p_{1},q_{1}}f\left( \overline{t}%
\right) \otimes e_{p_{2},q_{2}}f\left( \overline{t}\right) \\
&=&\left( f\otimes f\right) \left( \sum \left( -\xi \right)
^{q_{1}p_{2}}e_{p_{1},q_{1}}\overline{t}\otimes e_{p_{2},q_{2}}\overline{t}%
\right) =\left( f\otimes f\right) \Delta \left( \overline{t}\right)
\end{eqnarray*}
and such an $f$ is a
Hopf algebra isomorphism between $H_{\xi ,X^{k}Y^{l}}$ and $H_{-\xi
,X^{k}Y^{l}}$. Thus we may assume that $\xi =i$ and write $H_{i
,X^{k}Y^{l}}=H_{X^{k}Y^{l}}$. Define  \begin{eqnarray*}
f':H_{X^{k}Y} &\rightarrow &H_{X^{k}}\qquad \text{by} \\
f'\left( e_{r,s}\right) &=&e_{r+2s,s} \\
f'\left( \overline{t}\right) &=&\sum i^{q}e_{p,q}\overline{t}=
\left(\frac{1+i}{2}1+\frac{1-i}{2}Y\right)\overline{t}
\end{eqnarray*}
and extend it multiplicatively to $f'\left( e_{r,s}\overline{t}\right) $. Note that 
restriction $f' | _{\left( kD_{8}\right) ^{\ast }}$ corresponds to the group 
automorphism $f_{t}$ described in Case A) and $f'\left( X\right)=X, f'\left(
Y\right)=Y$. Then 
\begin{eqnarray*}
f'\left( \overline{t}\right) f'\left( \overline{t}\right) &=&\sum i^{q}e_{p,q}%
\overline{t}\sum i^{q}e_{p,q}\overline{t}=\sum i^{q}e_{p,q}\sum i^{q}e_{-p,q}%
\overline{t}^{2} \\
&=&\sum (-1)^{q}e_{p,q}X^{k}=YX^{k}=f'\left( YX^{k}\right) =
f'\left( \overline{t}^{2}\right) \\
f' \left( \overline{t}e_{r,s}\right) &=&f' \left( e_{-r,s}\overline{t}\right)
=e_{-r+2s,s}\sum i^{q}e_{p,q}\overline{t}=\sum i^{q}e_{p,q}\overline{t}%
e_{r-2s,s}=f' \left( \overline{t}\right) f' \left( e_{r,s}\right) \\
\end{eqnarray*}
There are at most two nonisomorphic Hopf algebras of this kind:

\begin{enumerate}
\item  $H_{B:1}$ with trivial cocycle and $\mathbf{G}\left( H_{B:1}^{\ast
}\right) =\left\langle \chi ,\varphi \right\rangle \cong D_{8}$, where $\varphi
(e_{r,s})=\delta _{r,2}\delta _{s,0}$, $\varphi \left( \overline{t}\right) =1$,
$\chi (e_{r,s})=\delta _{r,2}\delta _{s,1}$, $\chi (\overline{t})=-1$. 
There are two degree $2$ irreducible representations defined by \newline
\begin{tabular}{lll}
$\pi _{1}\left( e_{1,0}\right) =\left( 
\begin{array}{rr}
1 & 0 \\ 
0 & 0
\end{array}
\right) $ & $\pi _{1}\left( e_{3,0}\right) =\left( 
\begin{array}{rr}
0 & 0 \\ 
0 & 1
\end{array}
\right) $ & $\pi _{1}\left( \overline{t}\right) =\left( 
\begin{array}{rr}
0 & 1 \\ 
1 & 0
\end{array}
\right) $ \\ 
$\pi _{2}\left( e_{1,1}\right) =\left( 
\begin{array}{rr}
1 & 0 \\ 
0 & 0
\end{array}
\right) $ & $\pi _{2}\left( e_{3,1}\right) =\left( 
\begin{array}{rr}
0 & 0 \\ 
0 & 1
\end{array}
\right) $ & $\pi _{2}\left( \overline{t}\right) =\left( 
\begin{array}{rr}
0 & 1 \\ 
1 & 0
\end{array}
\right) $%
\end{tabular}
\newline
with the property $\pi _{k} ^{2}=\pi _{k} \bullet \pi _{k} =
1+\chi ^{2}+\varphi +\chi ^{2}\varphi $.
\medskip

\item  $H_{B:X}$ with the cocycle defined by $\sigma _{X}\left( t,t\right)
=X$ and $\mathbf{G}\left( H_{B:X}^{\ast }\right) =\left\langle \chi ,\varphi
\right\rangle \cong D_{8}$, where $\chi (e_{r,s})=\delta _{r,2}\delta _{s,1}$%
, $\chi (\overline{t})=-1$, $\varphi (e_{r,s})=\delta _{r,2}\delta _{s,0}$, $%
\varphi \left( \overline{t}\right) =1$.  There are two degree $2$
irreducible representations defined by \newline
\begin{tabular}{lll}
$\pi _{1}\left( e_{1,0}\right) =\left( 
\begin{array}{rr}
1 & 0 \\ 
0 & 0
\end{array}
\right) $ & $\pi _{1}\left( e_{3,0}\right) =\left( 
\begin{array}{rr}
0 & 0 \\ 
0 & 1
\end{array}
\right) $ & $\pi _{1}\left( \overline{t}\right) =\left( 
\begin{array}{rr}
0 & 1 \\ 
-1 & 0
\end{array}
\right) $ \\ 
$\pi _{2}\left( e_{1,1}\right) =\left( 
\begin{array}{rr}
1 & 0 \\ 
0 & 0
\end{array}
\right) $ & $\pi _{2}\left( e_{3,1}\right) =\left( 
\begin{array}{rr}
0 & 0 \\ 
0 & 1
\end{array}
\right) $ & $\pi _{2}\left( \overline{t}\right) =\left( 
\begin{array}{rr}
0 & 1 \\ 
-1 & 0
\end{array}
\right) $%
\end{tabular}
\newline
with the property $\pi _{k}^{2}=\pi _{k} \bullet \pi _{k} =
1+\chi ^{2}+\varphi +\chi ^{2}\varphi $.
\end{enumerate}

\subsubsection*{Case C)}

The action is given by $t\rightharpoonup e_{p,q}=e_{-p+q,q}$, corresponding
to 
\begin{eqnarray*}
f_{t}\left( x\right) &=&x^{-1} \\
f_{t}\left( y\right) &=&xy
\end{eqnarray*}
Then $Y$ is a central grouplike of $H$. Write $v=\sigma \left( t,t\right)
=\sum c_{p,q}e_{p,q}$ . By $\left( \ref{coc1}\right) $ and
$\left( \ref{Coc4}\right) $  
\begin{eqnarray}
c_{p,q}&=&c_{-p+q,q}  \label{Cond7}\\
c_{p+r+2rq,q+s}&=&\xi ^{qr+q\left( -r+s\right) }c_{p,q}c_{r,s}
=\xi ^{qs}c_{p,q}c_{r,s}  \label{Cond8}
\end{eqnarray}
Conditions $\left( \ref{Cond7}\right) $ and $\left( \ref{Cond8}\right) $
imply that 
\begin{eqnarray*}
c_{0,1} &=&c_{1,1} \\
c_{2,1} &=&c_{3,1} \\
c_{1,0} &=&c_{3,0} \\
c_{1,0}c_{0,1} &=&c_{1,1}=c_{0,1} \\
c_{0,1}c_{1,0} &=&c_{3,1} \\
\xi c_{0,1}c_{0,1} &=&c_{0,0}=1
\end{eqnarray*}
Thus $c_{1,0}=1$ and $c_{0,1}=\omega ^{k}$, where $\omega $ is a primitive $%
8 $-th root of $1$ and $\xi =\omega ^{-2k}$. Therefore 
\[
\sigma _{k}\left( t,t\right) =\sum \omega ^{kq}e_{p,q}=\frac{1+Y}{2}+\frac{%
\omega ^{k}\left( 1-Y\right) }{2} \ , \qquad k=0,\ldots ,7
\]
For $k=0,\ldots ,7$ let $H_{k}$ be the Hopf algebra with the structure 
described above with cocycle $\sigma _{k}\left( t,t\right) $. Define 
\begin{eqnarray*}
f :H_{k+2}&\rightarrow &H_{k}\qquad \text{by} \\
f\left( e_{p,q}\right) &=&e_{p,q} \\
f\left( \overline{t}\right) &=&\sum_{p,q}i^{p}e_{p,q}\overline{t}
\end{eqnarray*}
and extend it multiplicatively to $f\left( e_{p,q}\overline{t}\right) $.
Then 
\begin{eqnarray*}
f\left( \overline{t}\right) f\left( \overline{t}\right) &=&\sum i^{p}e_{p,q}%
\overline{t}\sum i^{p}e_{p,q}\overline{t}=\sum i^{p}e_{p,q}\sum
i^{p}e_{-p+q,q}\overline{t}^{2} \\
&=&\sum i^{p}i^{-p+q}e_{p,q}\sigma _{k}\left( t,t\right) =\sum
i^{q}e_{p,q}\sum \omega ^{kq}e_{p,q} \\
&=&\sum \omega ^{kq+2q}e_{p,q}=\sigma _{k+2}\left( t,t\right) =f\left(
\sigma _{k+2}\left( t,t\right) \right) =f\left( \overline{t}^{2}\right) \\
f\left( \overline{t}e_{p,q}\right) &=&f\left( e_{-p+q,q}\overline{t}\right)
=e_{-p+q,q}\sum i^{p}e_{p,q}\overline{t}=\sum i^{p}e_{p,q}\overline{t}%
e_{p,q}=f\left( \overline{t}\right) f\left( e_{p,q}\right) \\
\Delta \left( f\left( \overline{t}\right) \right) &=&\Delta \left( \sum
i^{p}e_{p,q}\overline{t}\right) =\Delta \left( \sum i^{p}e_{p,q}\right)
\Delta \left( \overline{t}\right) = \\
&=&\left( \sum i^{p_{1}+p_{2}+2q_{1}p_{2}}e_{p_{1},q_{1}}\otimes
e_{p_{2},q_{2}}\right) \left( \sum \omega ^{-2ks_{1}r_{2}}e_{r_{1},s_{1}}%
\overline{t}\otimes e_{r_{2},s_{2}}\overline{t}\right) \\
&=&\left( \sum i^{p_{1}+p_{2}+2q_{1}p_{2}}\omega
^{-2kq_{1}p_{2}}e_{p_{1},q_{1}}\overline{t}\otimes e_{p_{2},q_{2}}\overline{t%
}\right) \\
&=&\sum i^{p_{1}+p_{2}}\omega ^{-2\left( k+2\right)
q_{1}p_{2}}e_{p_{1},q_{1}}\overline{t}\otimes e_{p_{2},q_{2}}\overline{t} \\
&=&\sum \omega ^{-2\left( k+2\right) q_{1}p_{2}}e_{p_{1},q_{1}}f\left( 
\overline{t}\right) \otimes e_{p_{2},q_{2}}f\left( \overline{t}\right) \\
&=&\left( f\otimes f\right) \left( \sum \omega ^{-2\left( k+2\right)
q_{1}p_{2}}e_{p_{1},q_{1}}\overline{t}\otimes e_{p_{2},q_{2}}\overline{t}%
\right) =\left( f\otimes f\right) \Delta \left( \overline{t}\right)
\end{eqnarray*}
and such a $f$ is a Hopf algebra isomorphism between $H_{k}$ and $H_{k+2}$.
Thus there are exactly two nonisomorphic Hopf algebras of this type:

\begin{enumerate}
\item  $H_{C:1}=H_{0}$ with a trivial cocycle and $\xi =1$. Then $\mathbf{G}\left( 
H_{C:1}\right) =\left\langle X\overline{t},X\right\rangle
\cong D_{8}$ and $\mathbf{G}\left( H_{C:1}^{\ast }\right) =\left\langle \chi
\right\rangle \times \left\langle \varphi \right\rangle \cong C_{2}\times
C_{2}$, where $
\chi (e_{p,q})=\delta _{p,2}\delta _{q,0}$, $\chi (\overline{t})=1$, $
\varphi (e_{p,q})=\delta _{p,0}\delta _{q,0}$, $\varphi
\left( \overline{t}
\right) =-1$. There are three degree $2$ irreducible
representations defined by \newline
\begin{tabular}{lll}
$\pi _{1}\left( e_{0,1}\right) =\left( 
\begin{array}{rr}
1 & 0 \\ 
0 & 0
\end{array}
\right) $ & $\pi _{1}\left( e_{1,1}\right) =\left( 
\begin{array}{rr}
0 & 0 \\ 
0 & 1
\end{array}
\right) $ & $\pi _{1}\left( \overline{t}\right) =\left( 
\begin{array}{rr}
0 & 1 \\ 
1 & 0
\end{array}
\right) $ \\ 
$\pi _{2}\left( e_{1,0}\right) =\left( 
\begin{array}{rr}
1 & 0 \\ 
0 & 0
\end{array}
\right) $ & $\pi _{2}\left( e_{3,0}\right) =\left( 
\begin{array}{rr}
0 & 0 \\ 
0 & 1
\end{array}
\right) $ & $\pi _{2}\left( \overline{t}\right) =\left( 
\begin{array}{rr}
0 & 1 \\ 
1 & 0
\end{array}
\right) $ \\ 
$\pi _{3}\left( e_{2,1}\right) =\left( 
\begin{array}{rr}
1 & 0 \\ 
0 & 0
\end{array}
\right) $ & $\pi _{3}\left( e_{3,1}\right) =\left( 
\begin{array}{rr}
0 & 0 \\ 
0 & 1
\end{array}
\right) $ & $\pi _{3}\left( \overline{t}\right) =\left( 
\begin{array}{rr}
0 & 1 \\ 
1 & 0
\end{array}
\right) $%
\end{tabular}
\newline
with the property $\pi _{2}^{2}=\pi _{2}\bullet \pi _{2}=1+\chi +\varphi
+\chi \varphi $, $\pi _{1}^{2}=\pi _{3}^{2}=1+\varphi +\pi _{2}$.
\smallskip

\item  $H_{C:\sigma _{1}}$ with cocycle $\sigma _{1}$ defined by 
$\sigma _{1}\left(
t,t\right) =\sum \omega ^{q}e_{p,q}$ and $\xi =\omega ^{-2}$, where $\omega $
is a primitive $8$-th root of $1$. Then 
$\mathbf{G}\left( H_{C:\sigma _{1}}
\right) =\left\langle X\right\rangle \times \left\langle Y\right\rangle
\cong C_{2}\times C_{2}$ and 
$\mathbf{G}\left( H_{C:\sigma _{1}}^{\ast }\right)
=\left\langle \chi \right\rangle \times \left\langle \varphi \right\rangle
\cong C_{2}\times C_{2}$, where $\chi (e_{p,q})=\delta _{p,2}\delta _{q,0}$, 
$\chi (\overline{t})=1$, $\varphi (e_{p,q})=\delta _{p,0}\delta _{q,0}$, $%
\varphi \left( \overline{t}\right) =-1$. There are three degree $2$
irreducible representations defined by \newline
\begin{tabular}{lll}
$\pi _{1}\left( e_{0,1}\right) =\left( 
\begin{array}{rr}
1 & 0 \\ 
0 & 0
\end{array}
\right) $ & $\pi _{1}\left( e_{1,1}\right) =\left( 
\begin{array}{rr}
0 & 0 \\ 
0 & 1
\end{array}
\right) $ & $\pi _{1}\left( \overline{t}\right) =\left( 
\begin{array}{rr}
0 & \sqrt{\omega } \\ 
\sqrt{\omega } & 0
\end{array}
\right) $ \\ 
$\pi _{2}\left( e_{1,0}\right) =\left( 
\begin{array}{rr}
1 & 0 \\ 
0 & 0
\end{array}
\right) $ & $\pi _{2}\left( e_{3,0}\right) =\left( 
\begin{array}{rr}
0 & 0 \\ 
0 & 1
\end{array}
\right) $ & $\pi _{2}\left( \overline{t}\right) =\left( 
\begin{array}{rr}
0 & 1 \\ 
1 & 0
\end{array}
\right) $ \\ 
$\pi _{3}\left( e_{2,1}\right) =\left( 
\begin{array}{rr}
1 & 0 \\ 
0 & 0
\end{array}
\right) $ & $\pi _{3}\left( e_{3,1}\right) =\left( 
\begin{array}{rr}
0 & 0 \\ 
0 & 1
\end{array}
\right) $ & $\pi _{3}\left( \overline{t}\right) =\left( 
\begin{array}{rr}
0 & \sqrt{\omega } \\ 
\sqrt{\omega } & 0
\end{array}
\right) $%
\end{tabular}
\newline
with the property $\pi _{2}^{2}=\pi _{2}\bullet \pi _{2}=1+\chi +\varphi
+\chi \varphi $, $\pi _{1}^{2}=\pi _{3}^{2}=1+\varphi +\pi _{2}$.
\end{enumerate}

\subsection{Case of $G=Q_{8}$.}\label{3.4}
\ \newline
Let $G=Q_{8}=\left\langle x,y\mid 
x^{4}=1,y^{2}=x^{2},yx=x^{-1}y\right\rangle $. Let $\left\{ e_{pq}\right\}
_{p=0,1,2,3;q=0,1}$ be the basis of $K$, dual to the basis $\left\{
x^{p}y^{q}\right\} _{p=0,1,2,3;q=0,1}$ of $K^{\ast }=kQ_{8}$. Then 
\[
\Delta _{H}\left( e_{pq}\right) =\Delta _{K}\left( e_{pq}\right)
=\sum _{\substack {p_{1}+p_{2}+2q_{1}\left( p_{2}+q_{2}\right)\equiv p\mod4 \\   
q_{1}+q_{2}\equiv q\mod2 }}%
e_{p_{1}q_{1}}\otimes e_{p_{2}q_{2}} 
\]
and it is easy to check that elements 
\begin{eqnarray*}
X &=&\sum_{pq}\left( -1\right) ^{p}e_{pq} \\
Y &=&\sum_{pq}\left( -1\right) ^{q}e_{pq}
\end{eqnarray*}
are grouplike of order $2$. For $\overline{t}=1\#t$

\[
\Delta _{H}\left( \overline{t}\right) =\theta \left( t\right) \overline{t}%
\otimes \overline{t} 
\]
Dualizing $\left( \ref{ext1}\right) $ we get another extension 
\[
F^{\ast }\stackrel{\pi ^{\ast }}{\hookrightarrow }H^{\ast }\stackrel{i^{\ast
}}{\twoheadrightarrow }K^{\ast } 
\]
and as in \cite[2.4]{Ma1}, \cite[2.11]{Ma2} or \cite[2.1]{Ma5}, since $k$ is
algebraically closed, there exist units $\overline{x}$ and $\overline{y}\in
H^{\ast }$, such that such that $\overline{x}^{4}=1_{H^{\ast }}$, $\overline{%
y}^{2}=\overline{x}^{2}$, $\left\langle e_{pq},\overline{x}^{i}\overline{y}%
^{j}\right\rangle =\delta _{ip}\delta _{jq}$ and $\alpha =\overline{x}%
\overline{y}\overline{x}\overline{y}^{-1}\in F^{\ast }=k\left\{
e_{0},e_{1}\right\} $, where $\left\{ e_{r}\right\} $ is a dual basis of $%
\left\{ t^{r}\right\} $, $r=0,1$. $\varepsilon \left( \alpha \right)
=\varepsilon \left( \overline{x}\overline{y}\overline{x}\overline{y}%
^{-1}\right) =1$ and therefore $\alpha =e_{0}+\xi e_{1}$. The right action $%
\rho ^{\ast }:F^{\ast }\otimes K^{\ast }\rightarrow F^{\ast }$ is trivial,
thus $F^{\ast }$ lies in the center of $H^{\ast }$. Moreover, $\overline{x}%
^{2}=\overline{y}^{2}$ also lies in the center of $H^{\ast }$. Then 
\begin{eqnarray*}
\overline{x}\overline{y}\overline{x}^{-1}\overline{y}^{-1} &=&\overline{x}%
\overline{y}\overline{x}^{3}\overline{y}^{-1}=\overline{x}\overline{y}%
\overline{x}\overline{y}^{-1}\overline{x}^{2}=\alpha \overline{x}^{2} \\
\overline{x}^{3} &=&\overline{x}\overline{x}^{2}=\overline{x}\overline{y}%
^{2}=\alpha \overline{x}^{2}\overline{y}\overline{x}\overline{y}=\alpha 
\overline{x}^{2}\overline{y}\alpha \overline{x}^{2}\overline{y}\overline{x}%
=\alpha ^{2}\overline{x}^{4}\overline{y}^{2}\overline{x}=\alpha ^{2}%
\overline{x}^{3}
\end{eqnarray*}
Thus $\alpha ^{2}=1$ and $\xi =\pm 1$. 
\begin{eqnarray*}
\left\langle \Delta _{H}\left( \overline{t}\right) ,\overline{x}^{i}%
\overline{y}^{j}e_{k}\otimes \overline{x}^{p}\overline{y}^{q}e_{r}\right%
\rangle &=&\left\langle \overline{t},\overline{x}^{i}\overline{y}^{j}e_{k}%
\overline{x}^{p}\overline{y}^{q}e_{r}\right\rangle =\delta _{kr}\left\langle 
\overline{t},\overline{x}^{i+p}\left( \alpha \overline{x}^{2}\right) ^{-jp}%
\overline{y}^{j}\overline{y}^{q}e_{k}\right\rangle \\
&=&\delta _{kr}\left\langle \overline{t},\overline{x}^{i+p+2jp+2jq}\overline{%
y}^{j+q-2jq}\alpha ^{jp}e_{k}\right\rangle =\xi ^{jp}\delta _{k1}\delta _{r1}
\end{eqnarray*}
On the other hand 
\[
\left\langle \Delta _{H}\left( \overline{t}\right) ,\overline{x}^{i}%
\overline{y}^{j}e_{k}\otimes \overline{x}^{p}\overline{y}^{q}e_{r}\right%
\rangle =\left\langle \theta \left( t\right) \overline{t}\otimes \overline{t}%
,\overline{x}^{i}\overline{y}^{j}e_{k}\otimes \overline{x}^{p}\overline{y}%
^{q}e_{r}\right\rangle =\left\langle \theta \left( t\right) ,\overline{x}^{i}%
\overline{y}^{j}\otimes \overline{x}^{p}\overline{y}^{q}\right\rangle \delta
_{k1}\delta _{r1} 
\]
Therefore 
\[
\theta \left( t\right) =\sum_{ijpq}\xi ^{jp}e_{ij}\otimes e_{pq} 
\]

If $\xi =1$ then $\overline{t}$ is a grouplike of $H$, if $\xi =-1$ then $%
\sum i^{p+q^{2}}e_{p,q}\overline{t}$ is a grouplike of $H$. Thus $\mathbf{G}%
\left( H\right) $ has always order $8$.

Write $v=\sigma \left( t,t\right) =\sum c_{i,j}e_{i,j}$ then $%
c_{0,0}=\varepsilon \left( v\right) =1$ and $c_{i,j}\neq 0,$ since $v$ is a
unit and 
\[
\Delta _{H}\left( \overline{t}^{2}\right) =\Delta _{H}\left( v\right)
=\Delta _{K}\left( \sum c_{i,j}e_{i,j}\right) =\sum
c_{p+r+2rq+2sq,q+s}e_{p,q}\otimes e_{r,s} 
\]

Action by $t$  is a Hopf algebra map and therefore it induces a group
automorphism $f_{t}:G\rightarrow G$ defined by $\left\langle
e_{p,q},f_{t}\left( g\right) \right\rangle =\left\langle t\rightharpoonup
e_{p,q},g\right\rangle $, which has order $2$. Renaming generators we are
down to two choices for $f_{t}$; we consider them below:

\subsubsection*{Case D)}

The action is given by $t\rightharpoonup e_{i,j}=e_{i+2j,j}$, corresponding
to 
\begin{eqnarray*}
f_{t}\left( x\right) &=&x \\
f_{t}\left( y\right) &=&x^{2}y
\end{eqnarray*}
Then $X$ and $Y$ are central grouplikes of $H$. Thus $\mathbf{G}\left(
H\right) $ is abelian of order $8$ and $H$ was described in Section \ref{3.1} 
or Section \ref{3.2}.

\subsubsection*{Case E)}

The action is given by $t\rightharpoonup e_{i,j}=e_{-i+j,j}$, corresponding
to 
\begin{eqnarray*}
f_{t}\left( x\right)  &=&x^{-1} \\
f_{t}\left( y\right)  &=&xy
\end{eqnarray*}
Then $Y$ is a central grouplike of $H$. Write $v=\sigma \left(
t,t\right) =\sum c_{i,j}e_{i,j}$ . By $\left( \ref{coc1}\right) $ 
\begin{equation}
c_{i,j}=c_{-i+j,j}  \label{Cond1}
\end{equation}
On the other hand, for $H$ to be a bialgebra 
\begin{eqnarray*}
\Delta _{H}\left( \overline{t}^{2}\right)  &=&\Delta _{H}\left( \overline{t}%
\right) \Delta _{H}\left( \overline{t}\right) =\left( \theta \left( t\right) 
\overline{t}\otimes \overline{t}\right) \left( \theta \left( t\right) 
\overline{t}\otimes \overline{t}\right)  \\
&=&\left( \sum_{pqrs}\xi ^{rq}e_{pq}\otimes e_{rs}\right) \left(
\sum_{pqrs}\xi ^{rq}\left( t\rightharpoonup e_{pq}\right) \otimes \left(
t\rightharpoonup e_{rs}\right) \right) \sigma \left( t,t\right) \otimes
\sigma \left( t,t\right)  \\
&=&\sum \xi ^{rq}\xi ^{\left( -r+s\right) q}c_{p,q}c_{r,s}e_{p,q}\otimes
e_{r,s}=\sum \xi ^{qs}c_{p,q}c_{r,s}e_{p,q}\otimes e_{r,s}
\end{eqnarray*}
Therefore 
\begin{equation}
c_{p+r+2\left( r+s\right) q,q+s}=\xi ^{qs}c_{p,q}c_{r,s}  \label{Cond2}
\end{equation}
Conditions $\left( \ref{Cond1}\right) $ and $\left( \ref{Cond2}\right) $
imply that 
\begin{eqnarray*}
c_{0,1} &=&c_{1,1} \\
c_{2,1} &=&c_{3,1} \\
c_{1,0} &=&c_{3,0} \\
c_{1,0}c_{0,1} &=&c_{1,1}=c_{0,1} \\
c_{0,1}c_{1,0} &=&c_{3,1} \\
\xi c_{0,1}c_{0,1} &=&c_{2,0}=c_{1,0}c_{1,0}
\end{eqnarray*}
Thus $c_{1,0}=1$ and $c_{0,1}=i^{k}$, where $\xi =i^{2k}$  and $k=0,1,2,3$. 
Therefore 
\[
\sigma _{k}\left( t,t\right) =\sum i^{kq}e_{p,q}=\frac{1+Y}{2}+\frac{%
i^{k}\left( 1-Y\right) }{2}
\]
Let $H_{k}$ be the Hopf algebra with the structure described above with
cocycle $\sigma _{k}\left( t,t\right) $. Define 
\begin{eqnarray*}
f :H_{k}&\rightarrow &H_{k+1}\qquad \text{by} \\
f\left( e_{p,q}\right)  &=&e_{p,q} \\
f\left( \overline{t}\right)  &=&\sum_{p,q}i^{p+q^{2}}e_{p,q}\overline{t}
\end{eqnarray*}
and extend it multiplicatively to $f\left( e_{p,q}\overline{t}\right) $.
Then 
\begin{eqnarray*}
f\left( \overline{t}\right) f\left( \overline{t}\right)  &=&\sum
i^{p+q^{2}}e_{p,q}\overline{t}\sum i^{p+q^{2}}e_{p,q}\overline{t}=\sum
i^{p+q^{2}}e_{p,q}\sum i^{p+q^{2}}e_{-p+q,q}\overline{t}^{2} \\
&=&\sum i^{p+q^{2}}i^{-p+q+q^{2}}e_{p,q}\sigma _{k+1}\left( t,t\right) =\sum
i^{3q}e_{p,q}\sum i^{\left( k+1\right) q}e_{p,q} \\
&=&\sum i^{\left( k+1\right) q+3q}e_{p,q}=\sum i^{kq}e_{p,q}=\sigma
_{k}\left( t,t\right) =f\left( \sigma _{k}\left( t,t\right) \right) =f\left( 
\overline{t}^{2}\right)  \\
f\left( \overline{t}e_{p,q}\right)  &=&f\left( e_{-p+q,q}\overline{t}\right)
=e_{-p+q,q}\sum i^{p+q^{2}}e_{p,q}\overline{t}=\sum i^{p+q^{2}}e_{p,q}%
\overline{t}e_{p,q}=f\left( \overline{t}\right) f\left( e_{p,q}\right)
\end{eqnarray*}
\begin{eqnarray*} 
\Delta f\left( \overline{t}\right)  &=&\Delta \left( \sum
i^{p+q^{2}}e_{p,q}\overline{t}\right) =\Delta \left( \sum
i^{p+q^{2}}e_{p,q}\right) \Delta \left( \overline{t}\right) = \\
&=&\left( \sum i^{p_{1}+p_{2}+2q_{1}\left( p_{2}+q_{2}\right) +\left(
q_{1}+q_{2}\right) ^{2}}e_{p_{1},q_{1}}\otimes e_{p_{2},q_{2}}\right) \left(
\sum i^{2\left( k+1\right) s_{1}r_{2}}e_{r_{1},s_{1}}\overline{t}\otimes
e_{r_{2},s_{2}}\overline{t}\right)  \\
&=&\left( \sum i^{p_{1}+p_{2}+2q_{1}\left( p_{2}+q_{2}\right) +\left(
q_{1}+q_{2}\right) ^{2}+2\left( k+1\right) q_{1}p_{2}}e_{p_{1},q_{1}}%
\overline{t}\otimes e_{p_{2},q_{2}}\overline{t}\right)  \\
&=&\sum i^{p_{1}+p_{2}+q_{1}^{2}+q_{2}^{2}}i^{2kq_{1}p_{2}}e_{p_{1},q_{1}}%
\overline{t}\otimes e_{p_{2},q_{2}}\overline{t} \\
&=&\sum i^{2kq_{1}p_{2}}e_{p_{1},q_{1}}f\left( \overline{t}\right) \otimes
e_{p_{2},q_{2}}f\left( \overline{t}\right)  \\
&=&\left( f\otimes f\right) \left( \sum i^{2kq_{1}p_{2}}e_{p_{1},q_{1}}%
\overline{t}\otimes e_{p_{2},q_{2}}\overline{t}\right) =\left( f\otimes
f\right) \Delta \left( \overline{t}\right) 
\end{eqnarray*}
and such an $f$ is a Hopf algebra isomorphism between $H_{k}$ and $H_{k+1}$.
Thus there is exactly one Hopf algebra of this type: $H_{E}=H_{0}$
with a trivial cocycle and $\xi =1$. Then $\mathbf{G}\left( H_{E}
\right) =\left\langle X\overline{t},X\right\rangle \cong D_{8}$ and $%
\mathbf{G}\left( H_{E}^{\ast }\right) =\left\langle \chi \right\rangle
\times \left\langle \varphi \right\rangle \cong C_{2}\times C_{2}$, where $%
\chi (e_{p,q})=\delta _{p,2}\delta _{q,0}$, $\chi (\overline{t})=1$, $%
\varphi (e_{p,q})=\delta _{p,0}\delta _{q,0}$, $\varphi \left( \overline{t}%
\right) =-1$. There are three degree $2$ irreducible representations defined
by \newline
\begin{tabular}{lll}
$\pi _{1}\left( e_{0,1}\right) =\left( 
\begin{array}{rr}
1 & 0 \\ 
0 & 0
\end{array}
\right) $ & $\pi _{1}\left( e_{1,1}\right) =\left( 
\begin{array}{rr}
0 & 0 \\ 
0 & 1
\end{array}
\right) $ & $\pi _{1}\left( \overline{t}\right) =\left( 
\begin{array}{rr}
0 & 1 \\ 
1 & 0
\end{array}
\right) $ \\ 
$\pi _{2}\left( e_{1,0}\right) =\left( 
\begin{array}{rr}
1 & 0 \\ 
0 & 0
\end{array}
\right) $ & $\pi _{2}\left( e_{3,0}\right) =\left( 
\begin{array}{rr}
0 & 0 \\ 
0 & 1
\end{array}
\right) $ & $\pi _{2}\left( \overline{t}\right) =\left( 
\begin{array}{rr}
0 & 1 \\ 
1 & 0
\end{array}
\right) $ \\ 
$\pi _{3}\left( e_{2,1}\right) =\left( 
\begin{array}{rr}
1 & 0 \\ 
0 & 0
\end{array}
\right) $ & $\pi _{3}\left( e_{3,1}\right) =\left( 
\begin{array}{rr}
0 & 0 \\ 
0 & 1
\end{array}
\right) $ & $\pi _{3}\left( \overline{t}\right) =\left( 
\begin{array}{rr}
0 & 1 \\ 
1 & 0
\end{array}
\right) $%
\end{tabular}
\newline
with the property $\pi _{2}^{2}=\pi _{2}\bullet \pi _{2}=1+\chi +\varphi
+\chi \varphi $, $\pi _{1}^{2}=\pi _{3}^{2}=\chi +\chi \varphi +\pi _{2}$.
\subsection{Summary}
\begin{proposition}
\label{th3}Let $H$ be a nontrivial semisimple Hopf algebra of dimension $16$. 
Then $\mathbf{G}\left( H\right) $ is abelian of order $8$ if and only if  
$\mathbf{G}\left( H^{\ast }\right) $ is abelian of order $8$.
\end{proposition}
\proof
All nontrivial Hopf algebras with abelian groups of grouplikes were described
in Sections \ref{3.1} and \ref{3.2} and their duals have abelian groups of 
grouplikes.
\endproof
\begin{proposition}\label{4x2}
There are exactly $7$ nonisomorphic nontrivial semisimple Hopf algebras
of dimension $16$ with $\mathbf{G}\left( H\right) 
\cong C_{4}\times C_{2}$.
\end{proposition}
\proof
All nontrivial Hopf algebras with $\mathbf{G}\left( H\right) 
\cong C_{4}\times C_{2}$ were described in Section \ref{3.1}. There are at
most $7$ nonisomorphic Hopf algebras with $\mathbf{G}\left( H\right) \cong
C_{4}\times C_{2}$, namely $H_{a:1}$, $H_{a:y}$, $H_{b:1}$, $H_{b:y}$, $%
H_{b:x^{2}y}$, $H_{c:\sigma _{0}}$, $H_{c:\sigma _{1}}$. 

Assume 
$f$ is a Hopf algebra isomorphism between Hopf algebras $H_{1}$ and $H_{2}$
with $\mathbf{G}\left( H_{1}\right) \cong \mathbf{G}\left( H_{2}\right)
\cong C_{4}\times C_{2}$. Then we get a group isomorphism 
\begin{equation*}
f\mid _{\mathbf{G}\left( H_{1}\right) }:\mathbf{G}\left( H_{1}\right)
\rightarrow \mathbf{G}\left( H_{2}\right)
\end{equation*}
Write $\mathbf{G}\left( H_{1}\right) \cong \mathbf{G}\left( H_{2}\right)
=\left\langle x\right\rangle \times \left\langle y\right\rangle $, where $%
\left| x\right| =4$ and $\left| y\right| =2$. Then the dual basis of $k%
\mathbf{G}\left( H_{1}\right) \cong k\mathbf{G}\left( H_{2}\right) $ is
given by 
\begin{equation*}
e_{pq}=1/8\left( 1+i^{p}x+i^{2p}x^{2}+i^{3p}x^{3}\right) \left( 1+\left(
-1\right) ^{q}y\right) ,\text{\qquad }p=0,1,2,3;q=0,1
\end{equation*}
Write 
\begin{eqnarray*}
f\left( e_{p,q}\right) &=&e_{\alpha _{1}\left( p,q\right) ,\alpha _{2}\left(
p,q\right) } \\
f^{-1}\left( e_{p,q}\right) &=&e_{\beta _{1}\left( p,q\right) ,\beta
_{2}\left( p,q\right) }
\end{eqnarray*}
where $\alpha _{1}\left( p,q\right) ,\beta _{1}\left( p,q\right) \in \left\{
0,1,2,3\right\} $ and $\alpha _{2}\left( p,q\right) ,\beta _{2}\left(
p,q\right) \in \left\{ 0,1\right\} $.

Write $\left\{ e_{pq}\overline{t^{r}}\right\} _{p=0,1,2,3;q=0,1;r=0,1}$ and $%
\left\{ e_{pq}\overline{T^{r}}\right\} _{p=0,1,2,3;q=0,1;r=0,1}$ for the
bases of $H_{1}$ and $H_{2}$, respectively. Write 
\begin{equation*}
f\left( \overline{t}\right) =\sum_{p,q,r}\lambda _{p,q,r}e_{pq}\overline{%
T^{r}}
\end{equation*}

Then 
\begin{eqnarray*}
\Delta f\left( \overline{t}\right) &=&\Delta \left( \sum_{p,q,r}\lambda
_{p,q,r}e_{pq}\overline{T^{r}}\right) =\sum_{p,q}\lambda _{p,q,0}\Delta
\left( e_{pq}\right) +\sum_{p,q}\lambda _{p,q,1}\Delta \left( e_{pq}\right)
\Delta \left( \overline{T}\right) \\
&=&\sum \lambda _{p_{1}+p_{2},q_{1}+q_{2},0}e_{p_{1}q_{1}}\otimes
e_{p_{2}q_{2}} \\
&&+\left( \sum \lambda _{p_{1}+p_{2},q_{1}+q_{2},1}e_{p_{1}q_{1}}\otimes
e_{p_{2}q_{2}}\right) \left( \sum \left( -1\right) ^{bc}e_{ab}\overline{T}%
\otimes e_{cd}\overline{T}\right) \\
&=&\sum \lambda _{p_{1}+p_{2},q_{1}+q_{2},0}e_{p_{1}q_{1}}\otimes
e_{p_{2}q_{2}}+\sum \left( -1\right) ^{p_{2}q_{1}}\lambda
_{p_{1}+p_{2},q_{1}+q_{2},1}e_{p_{1}q_{1}}\overline{T}\otimes e_{p_{2}q_{2}}%
\overline{T}
\end{eqnarray*}
\begin{eqnarray*}
\left( f\otimes f\right) \Delta \left( \overline{t}\right) &=&\left(
f\otimes f\right) \left( \sum \left( -1\right) ^{p_{2}q_{1}}e_{p_{1}q_{1}}%
\overline{T}\otimes e_{p_{2}q_{2}}\overline{T}\right) \\
&=&\sum \left( -1\right) ^{p_{2}q_{1}}f\left( e_{p_{1}q_{1}}\right)
\sum_{p,q,r}\lambda _{p,q,r}e_{pq}\overline{T^{r}}\otimes f\left(
e_{p_{2}q_{2}}\right) \sum_{p,q,r}\lambda _{p,q,r}e_{pq}\overline{T^{r}} \\
&=&\sum \left( -1\right) ^{\beta _{1}\left( p_{2},q_{2}\right) \beta
_{2}\left( p_{1},q_{1}\right) }\lambda _{p_{1},q_{1},r_{1}}\lambda
_{p_{2},q_{2},r_{2}}e_{p_{1}q_{1}}\overline{T^{r_{1}}}\otimes e_{p_{2}q_{2}}%
\overline{T^{r_{2}}}
\end{eqnarray*}
Since $f$ is a coalgebra map, 
\begin{equation*}
\Delta f\left( \overline{t}\right) =\left( f\otimes f\right) \Delta \left( 
\overline{t}\right)
\end{equation*}
and therefore $\lambda _{p_{1},q_{1},0}\lambda _{p_{2},q_{2},1}=0$ for all $%
p_{1},p_{2}\in \left\{ 0,1,2,3\right\} $, $q_{1},q_{2}\in \left\{
0,1\right\} $. Thus either $\lambda _{p,q,0}=0$ for all $p\in \left\{
0,1,2,3\right\} $, $q\in \left\{ 0,1\right\} $ or $\lambda _{p,q,1}=0$ for
all $p\in \left\{ 0,1,2,3\right\} $, $q\in \left\{ 0,1\right\} $. In the
latter case $f\left( \overline{t}\right) =\sum \lambda _{p,q,0}e_{pq}\in k%
\mathbf{G}\left( H_{2}\right) $, which contradicts the bijectivity of $f$.
Therefore $\lambda _{p,q,0}=0$ for all $p\in \left\{ 0,1,2,3\right\} $, $%
q\in \left\{ 0,1\right\} $. Write $\lambda _{p,q}=\lambda _{p,q,1}$. Then 
\begin{equation*}
f\left( \overline{t}\right) =\sum_{p,q}\lambda _{p,q}e_{pq}\overline{T}
\end{equation*}
and so, applying $\varepsilon $, also 
\begin{equation*}
\lambda _{0,0}=\varepsilon \left( \overline{t}\right) =1
\end{equation*}
Moreover, since 
$\sum \left( -1\right) ^{p_{2}q_{1}}\lambda
_{p_{1}+p_{2},q_{1}+q_{2}}e_{p_{1}q_{1}}\overline{T}\otimes e_{p_{2}q_{2}}%
\overline{T}$\newline 
$=\sum \left( -1\right) ^{\beta _{1}\left( p_{2},q_{2}\right)
\beta _{2}\left( p_{1},q_{1}\right) }\lambda _{p_{1},q_{1}}\lambda
_{p_{2},q_{2}}e_{p_{1}q_{1}}\overline{T}\otimes e_{p_{2}q_{2}}\overline{T}$
we get 
\begin{equation}
\lambda _{p_{1}+p_{2},q_{1}+q_{2}}=\left( -1\right) ^{p_{2}q_{1}}\left(
-1\right) ^{\beta _{1}\left( p_{2},q_{2}\right) \beta _{2}\left(
p_{1},q_{1}\right) }\lambda _{p_{1},q_{1}}\lambda _{p_{2},q_{2}}
\label{f_coalg}
\end{equation}
for any $p_{1},p_{2}\in \left\{ 0,1,2,3\right\} $, $q_{1},q_{2}\in \left\{
0,1\right\} $.

Let $u\in k\mathbf{G}\left( H_{1}\right) $. Then 
\begin{eqnarray*}
f\left( t\rightharpoonup _{1}u\right) f\left( \overline{t}\right) &=&f\left(
\left( t\rightharpoonup _{1}u\right) \overline{t}\right) =f\left( \overline{t%
}u\right) =f\left( \overline{t}\right) f\left( u\right) =\sum \lambda
_{p,q}e_{p,q}\overline{T}f\left( u\right) \\
&=&\left( t\rightharpoonup _{2}f\left( u\right) \right) \sum \lambda
_{p,q}e_{p,q}\overline{T}=\left( t\rightharpoonup _{2}f\left( u\right)
\right) f\left( \overline{t}\right)
\end{eqnarray*}
Thus, since $\overline{t}$ is a unit ($\overline{t}^{2}=\sigma \left(
t,t\right) $ is a unit), 
\begin{equation}
f\left( t\rightharpoonup _{1}u\right) =t\rightharpoonup _{2}f\left( u\right)
\label{f_t}
\end{equation}

Let's show that Hopf algebras from types $H_{a}$, $H_{b}$ and $H_{c}$ cannot
be isomorphic to each other. $K_{0}\left( H_{c}\right) \ncong K_{0}\left(
H_{a}\right) $ or $K_{0}\left( H_{b}\right) $, thus $H_{c}\ncong H_{a}$ or $%
H_{b}$. If $f:H_{a}\rightarrow H_{b}$ then by formula (\ref{f_t}) 
\begin{equation*}
f\left( x\right) f\left( y\right) =f\left( xy\right) =f\left(
t\rightharpoonup _{1}x\right) =t\rightharpoonup _{2}f\left( x\right)
=f\left( x\right) ^{-1}
\end{equation*}
and therefore $f\left( y\right) =f\left( x^{2}\right) $, which is impossible
if $f$ is an isomorphism.

$H_{b:1}\ncong H_{b:y}$ or $H_{b:x^{2}y}$ and $H_{c:\sigma _{0}}\ncong 
H_{c:\sigma _{1}}$ since
their duals have nonisomorphic groups of grouplikes. Thus there are at least 
$5$ nonisomorphic Hopf algebras with $\mathbf{G}\left( H\right) \cong
C_{4}\times C_{2}$, namely $H_{a:1}$, $H_{b:1}$, $H_{b:y}$, $H_{c:\sigma _{0}}$
and $H_{c:\sigma _{1}}$.

Now we prove that $H_{a:1}\ncong H_{a:y}$ and $H_{b:y}\ncong H_{b:x^{2}y}$.
If $f$ is a Hopf algebra isomorphism as before, then since $f\mid _{\mathbf{%
G}\left( H_{1}\right) }$ is a group isomorphism $f\left( x\right) \in
\left\{ x,x^{-1},xy,x^{-1}y\right\} $, $f\left( y\right) \in \left\{
y,x^{2}y\right\} $ and $f\left( x^{2}\right) =x^{2}$.

If $f\left( x\right) =x^{2k+1}y^{l}$ and $f\left( y\right) =y$, where $k,l=0,1$
then
\begin{equation*}
f^{-1}\left( e_{p,q}\right) =f\left( e_{p,q}\right) =
e_{\left( 2k+1\right)p+2lq,q}
\end{equation*}
and by formula (\ref{f_coalg}) 
\begin{equation*}
\lambda _{p_{1}+p_{2},q_{1}+q_{2}}=\left( -1\right) ^{p_{2}q_{1}}\left(
-1\right) ^{\left( \left( 2k+1\right) p_{2}+2lq_{2}\right) q_{1}}
\lambda _{p_{1},q_{1}}\lambda _{p_{2},q_{2}}=\lambda
_{p_{1},q_{1}}\lambda _{p_{2},q_{2}}
\end{equation*}

If $f\left( x\right) =x^{2k+1}y^{l}$ and $f\left( y\right) =x^{2}y$, where 
$k,l=0,1$ then
\begin{eqnarray*}
f\left( e_{p,q}\right) &=&e_{\left( 2k+2l+1\right) p+2lq,p+q} \\
f^{-1}\left( e_{p,q}\right) &=&e_{\left( 2k+1\right)p+2lq,p+q}
\end{eqnarray*}
and by formula (\ref{f_coalg}) 
\begin{equation*}
\lambda _{p_{1}+p_{2},q_{1}+q_{2}}=\left( -1\right) ^{p_{2}q_{1}}\left(
-1\right) ^{\left(\left( 2k+1\right) p_{2}+2lq_{2}\right) 
\left( q_{1}+p_{1}\right) }\lambda
_{p_{1},q_{1}}\lambda _{p_{2},q_{2}}=\left( -1\right) ^{p_{1}p_{2}}\lambda
_{p_{1},q_{1}}\lambda _{p_{2},q_{2}}
\end{equation*}

Now assume $f$ is a Hopf algebra isomorphism 
\begin{equation*}
f:H_{a:1}\rightarrow H_{a:y}
\end{equation*}
If $f\left( y\right) =x^{2}y$ then by formula (\ref
{f_t}) 
\begin{equation*}
f\left( x\right) y=t\rightharpoonup _{2}f\left( x\right) =f\left(
t\rightharpoonup _{1}x\right) =f\left( xy\right) =f\left( x\right) f\left(
y\right) =f\left( x\right) x^{2}y
\end{equation*}
that is $x^{2}=1$, which contradicts the fact that $\left| x\right| =4$.
Thus $f\left( y\right) =y$ and therefore 
\begin{equation*}
\lambda _{p_{1}+p_{2},q_{1}+q_{2}}=\lambda _{p_{1},q_{1}}\lambda
_{p_{2},q_{2}}
\end{equation*}
and thus 
\begin{equation*}
\lambda _{1,0}^{4}=\lambda _{2,0}^{2}=\lambda _{0,1}^{2}=\lambda _{0,0}=1
\end{equation*}
Then 
\begin{eqnarray*}
f\left( \overline{t}\right) f\left( \overline{t}\right) &=&\sum \lambda
_{p,q}e_{pq}\overline{T}\sum \lambda _{r,s}e_{rs}\overline{T} \\
&=&\sum \lambda _{p,q}e_{pq}\sum \lambda _{r,s}e_{r+2s,s}\overline{T}%
^{2}=\sum \lambda _{p,q}\lambda _{p+2q,q}e_{pq}\sigma _{H_{a:y}}\left(
t,t\right) \\
&=&\sum \lambda _{2\left( p+q\right) ,0}e_{pq}\sigma _{H_{a:y}}\left(
t,t\right) =\sum \lambda _{2,0}^{p+q}e_{pq}\sigma _{H_{a:y}}\left( t,t\right)
\end{eqnarray*}
If $\lambda _{2,0}=1$, $f\left( \overline{t}\right) f\left( \overline{t}%
\right) =\sigma _{H_{a:y}}\left( t,t\right) =y\neq f\left( \overline{t}%
^{2}\right) =1$.\newline
If $\lambda _{2,0}=-1$, $f\left( \overline{t}\right) f\left( \overline{t}%
\right) =x^{2}y\sigma _{H_{a:y}}\left( t,t\right) =x^{2}\neq f\left( 
\overline{t}^{2}\right) =1$.

Therefore, there is no Hopf algebra isomorphism between $H_{a:1}$ and $%
H_{a:y}$.

Now assume $f$ is a Hopf algebra isomorphism 
\begin{equation*}
f:H_{b:y}\rightarrow H_{b:x^{2}y}
\end{equation*}
Then 
\begin{eqnarray*}
f\left( \overline{t}\right) f\left( \overline{t}\right) &=&\sum \lambda
_{p,q}e_{pq}\overline{T}\sum \lambda _{r,s}e_{rs}\overline{T} \\
&=&\sum \lambda _{p,q}e_{pq}\sum \lambda _{r,s}e_{-r,s}\overline{T}^{2}=\sum
\lambda _{p,q}\lambda _{-p,q}e_{pq}\sigma _{H_{b:x^{2}y}}\left( t,t\right)
\end{eqnarray*}

$f\left( y\right) =y$ is not possible, since then we have 
\begin{equation*}
\lambda _{p,q}\lambda _{-p,q}=\lambda _{0,0}=1
\end{equation*}
and thus 
\begin{equation*}
f\left( \overline{t}\right) f\left( \overline{t}\right) =\left( \sum
e_{pq}\right) \sigma _{H_{b:x^{2}y}}\left( t,t\right) =\sigma
_{H_{b:x^{2}y}}\left( t,t\right) =x^{2}y\neq y=f\left( y\right) =f\left( 
\overline{t}^{2}\right) .
\end{equation*}

$f\left( y\right) =x^{2}y$ is not possible, since then we get 
\begin{equation*}
\lambda _{p_{1}+p_{2},q_{1}+q_{2}}=\left( -1\right) ^{p_{1}p_{2}}\lambda
_{p_{1},q_{1}}\lambda _{p_{2},q_{2}}
\end{equation*}
so 
\begin{equation*}
\lambda _{p,q}\lambda _{-p,q}=\left( -1\right) ^{p^{2}}\lambda _{0,0}=\left(
-1\right) ^{p^{2}}
\end{equation*}
and 
\begin{equation*}
f\left( \overline{t}\right) f\left( \overline{t}\right) =\left( \sum \left(
-1\right) ^{p^{2}}e_{pq}\right) \sigma _{H_{b:x^{2}y}}\left( t,t\right)
=x^{2}x^{2}y =y=f\left( x^{2}y\right) \neq
f\left( y\right) =f\left( \overline{t}^{2}\right) .
\end{equation*}

Therefore, there is no Hopf algebra isomorphism between $H_{b:y}$ and $%
H_{b:x^{2}y}$.

Thus there are exactly $7$ nonisomorphic Hopf algebras with $\mathbf{G}%
\left( H\right) \cong C_{4}\times C_{2}$, namely $H_{a:1}$, $H_{a:y}$, $%
H_{b:1}$, $H_{b:y}$, $H_{b:x^{2}y}$, $H_{c:\sigma _{0}}$ and 
$H_{c:\sigma _{1}}$.
\endproof
\begin{proposition}\label{2x2x2}
There are at least $2$ and at most $4$ nonisomorphic nontrivial semisimple 
Hopf algebras of dimension $16$ with $\mathbf{G}\left( H\right) 
\cong C_{2}\times C_{2}\times C_{2}$.
\end{proposition}
\proof
All nontrivial Hopf algebras with $\mathbf{G}\left( H\right) 
\cong C_{2}\times C_{2}\times C_{2}$ were described in Section \ref{3.2}. 
There are at
most $4$ nonisomorphic Hopf algebras with $\mathbf{G}\left( H\right) \cong
C_{2}\times C_{2}\times C_{2}$, namely $H_{d:1,1}$, $H_{d:1,-1}$, $H_{d:-1,1}$ 
and $H_{d:-1,-1}$. At least two of them are not isomorphic, since
$\mathbf{G}\left(H_{d:1,1}^{\ast }\right) \ncong 
\mathbf{G}\left(H_{d:1,-1}^{\ast }\right)$.
\endproof
\begin{proposition}\label{D_8}
There are exactly $2$ nonisomorphic nontrivial semisimple 
Hopf algebras of dimension $16$ with a commutative subHopfalgebra of dimension
$8$ and nonabelian $\mathbf{G}\left( H^{\ast }\right)$. In this case
$\mathbf{G}\left( H^{\ast }\right) \cong D_{8}$, $\mathbf{G}\left( H\right) 
\cong C_{2}\times C_{2}$ and $H^{\ast }$ also has a commutative 
subHopfalgebra of dimension $8$.
\end{proposition}
\proof
All nontrivial Hopf algebras with a commutative subHopfalgebra of dimension
$8$ and nonabelian $\mathbf{G}\left( H^{\ast }\right)$ were described in 
Section \ref{3.3}, Case B). There are at most $2$ of them, namely
$H_{B:1}$ and $H_{B:X}$.

Let's compute all the possible $8$-dimensional Hopf quotients of $H_{B}$. 
There is a one-to-one correspondence between hereditary subrings of $K_{0}(H)$
and Hopf quotients of $H$ (see \cite[Theorem 6]{NR} or 
\cite[Proposition 3.11]{N1}). Thus $H_{B}$ has $3$ quotients
of dimension $8$ corresponding to the hereditary subrings 
$R_{1}=\{ a1 + b\varphi +c\chi ^{2} + d\chi ^{2}\varphi+e\pi_{1} \in
K_{0}(H): a,b,c,d,e \in \mathbb{Z} \} $, 
$R_{2}=\{ a1 + b\varphi +c\chi ^{2} + d\chi ^{2}\varphi+e\pi_{2} \in
K_{0}(H): a,b,c,d,e \in \mathbb{Z} \}$ and 
$R_{3}=\{ \sum _{i=0}^{4} \sum _{j=0}^{2} a_{i,j}\chi ^{i}\varphi ^{j}\in 
K_{0}(H): a_{i,j}\in \mathbb{Z} \}$. They are obtained by factoring modulo 
normal 
ideals $\left( Y-1\right) H$, $\left( X-1\right) H$ and $\left( XY-1\right) H$,
where $X$, $Y$ and 
$XY$ are central grouplikes of $H$. It is easy to see that $H/\left(
Y-1\right) H$ is 
cocommutative (in fact, $H_{B:1}/\left( Y-1\right) H_{B:1}
\cong kD_{8}$ and $H_{B:X}/\left( Y-1\right) H_{B:X} \cong kQ_{8}$),
$H/\left( X-1\right) H$ is 
commutative (therefore $H/\left(X-1\right) H
\cong \left( kD_{8}\right)^{\ast }$ since $\mathbf{G}\left( H^{\ast }\right)
\cong D_{8}$) and $H/\left( XY-1\right) H$ is neither commutative nor
cocommutative (therefore $H/\left( X-1\right) H \cong H_{8}$). Therefore we 
see that $H_{B:1} \ncong H_{B:X}$ since they have different sets of
quotients. Both $ H_{B:1}$ and $H_{B:X}$ have cocommutative Hopf quotients of 
dimension 
$8$, $kD_{8}$ and $kQ_{8}$ respectively. Thus their duals  were described in 
Section \ref{3.4}. In particular,  $ H_{B:1}\cong H_{C:1}^{\ast }$,
$ H_{B:X}\cong H_{E}^{\ast }$ and $\mathbf{G}\left( H_{C:1}^{\ast }\right)\cong
\mathbf{G}\left( H_{E}^{\ast }\right)\cong C_{2}\times C_{2}$.
\endproof
\section{Nonabelian groups of order 16.} \label{sec4}

There are nine nonabelian groups of order $16$ (see \cite[118]{B}). The
first four of them are of exponent $8$, the last five of exponent $4$ (we
denote the quaternion group of order $8$ by $Q_{8}$ and the quasiquaternion
group of order $16$ by $Q_{16}$):

\begin{enumerate}
\item  $G_{1}=\left\langle a,b:a^{8}=b^{2}=1,ba=a^{5}b\right\rangle $\newline
$\mathbf{G}\left( \left( kG_{1}\right) ^{\ast }\right) =\left\langle \chi
\right\rangle \times \left\langle \varphi \right\rangle \cong C_{4}\times
C_{2}$, where $\chi (a)=i$, $\chi (b)=1$, $\varphi \left( a\right) =1$, $%
\varphi \left( b\right) =-1$. Degree $2$ irreducible representations of $%
G_{1}$ are defined by
\newline
\begin{tabular}{ll}
$\pi _{1}(a)=\left( 
\begin{array}{rr}
\omega  & 0 \\ 
0 & -\omega 
\end{array}
\right) $ & $\pi _{1}(b)=\left( 
\begin{array}{rr}
0 & 1 \\ 
1 & 0
\end{array}
\right) $ \\ 
$\pi _{2}(a)=\left( 
\begin{array}{rr}
\omega ^{3} & 0 \\ 
0 & -\omega ^{3}
\end{array}
\right) $ & $\pi _{2}(b)=\left( 
\begin{array}{rr}
0 & 1 \\ 
1 & 0
\end{array}
\right) $%
\end{tabular}
\newline
where $\omega $ is a primitive $8^{th}$-root of unity and $\pi _{1}^{2}=\chi
+\chi ^{3}+\chi \varphi +\chi ^{3}\varphi =\pi _{2}^{2}.$
\smallskip

\item  $G_{2}=\left\langle a,b:a^{8}=b^{2}=1,ba=a^{3}b\right\rangle $\newline
$\mathbf{G}\left( \left( kG_{2}\right) ^{\ast }\right) =\left\langle \chi
\right\rangle \times \left\langle \varphi \right\rangle \cong C_{2}\times
C_{2}$, where $\chi (a)=-1$, $\chi (b)=1$, $\varphi \left( a\right) =1$, $%
\varphi \left( b\right) =-1$. Degree $2$ irreducible representations of $%
G_{2}$ are defined by \newline
\begin{tabular}{lll}
$\pi _{1}(a)=\left( 
\begin{array}{rr}
\omega  & 0 \\ 
0 & \omega ^{3}
\end{array}
\right) $ &
$\pi _{2}(a)=\left( 
\begin{array}{rr}
i & 0 \\ 
0 & -i
\end{array}
\right) $ & 
$\pi _{3}(a)=\left( 
\begin{array}{rr}
\omega ^{5} & 0 \\ 
0 & \omega ^{7}
\end{array}
\right) $ \\
$\pi _{1}(b)=\left( 
\begin{array}{rr}
0 & 1 \\ 
1 & 0
\end{array}
\right) $ & $\pi _{2}(b)=\left( 
\begin{array}{rr}
0 & 1 \\ 
1 & 0
\end{array}
\right) $ &  $\pi _{3}(b)=\left( 
\begin{array}{rr}
0 & 1 \\ 
1 & 0
\end{array}
\right) $%
\end{tabular}
\newline
where $\omega $ is a primitive $8^{th}$-root of unity, and representations
satisfy the properties \newline
\begin{tabular}{lll}
$\pi _{1}^{2}=\chi +\chi \varphi +\pi _{2}=\pi _{3}^{2}$ & $\chi \bullet \pi
_{1}=\pi _{3}$ & $\chi \bullet \pi _{3}=\pi _{1}$ \\ 
$\pi _{2}^{2}=1+\chi +\varphi +\chi \varphi $ & $\varphi \bullet \pi
_{1}=\pi _{1}$ & $\varphi \bullet \pi _{3}=\pi _{3}$%
\end{tabular}
\smallskip

\item  $G_{3}=\left\langle a,b:a^{8}=b^{2}=1,ba=a^{-1}b\right\rangle =D_{16}$%
, the dihedral group\newline
$\mathbf{G}\left( \left( kG_{3}\right) ^{\ast }\right) =\left\langle \chi
\right\rangle \times \left\langle \varphi \right\rangle \cong C_{2}\times
C_{2}$, where $\chi (a)=-1$, $\chi (b)=1$, $\varphi \left( a\right) =1$, $%
\varphi \left( b\right) =-1$. Degree $2$ irreducible representations of $%
G_{3}$ are defined by \newline
\begin{tabular}{lll}
$\pi _{1}(a)=\left( 
\begin{array}{rr}
\omega  & 0 \\ 
0 & \omega ^{7}
\end{array}
\right) $ &
$\pi _{2}(a)=\left( 
\begin{array}{rr}
i & 0 \\ 
0 & -i
\end{array}
\right) $ &
$\pi _{3}(a)=\left( 
\begin{array}{rr}
\omega ^{3} & 0 \\ 
0 & \omega ^{5}
\end{array}
\right) $\\
 $\pi _{1}(b)=\left( 
\begin{array}{rr}
0 & 1 \\ 
1 & 0
\end{array}
\right) $  & $\pi _{2}(b)=\left( 
\begin{array}{rr}
0 & 1 \\ 
1 & 0
\end{array}
\right) $ & $\pi _{3}(b)=\left( 
\begin{array}{rr}
0 & 1 \\ 
1 & 0
\end{array}
\right) $%
\end{tabular}
\newline
where $\omega $ is a primitive $8^{th}$-root of unity, and representations
satisfy the properties \newline
\begin{tabular}{lll}
$\pi _{1}^{2}=1+\varphi +\pi _{2}=\pi _{3}^{2}$ & $\chi \bullet \pi _{1}=\pi
_{3}$ & $\chi \bullet \pi _{3}=\pi _{1}$ \\ 
$\pi _{2}^{2}=1+\chi +\varphi +\chi \varphi $ & $\varphi \bullet \pi
_{1}=\pi _{1}$ & $\varphi \bullet \pi _{3}=\pi _{3}$%
\end{tabular}
\smallskip

\item  $G_{4}=\left\langle a,b:a^{8}=1,b^{2}=a^{4},ba=a^{-1}b\right\rangle
=Q_{16}$, the quasiquaternion group\newline
$\mathbf{G}\left( \left( kG_{4}\right) ^{\ast }\right) =\left\langle \chi
\right\rangle \times \left\langle \varphi \right\rangle \cong C_{2}\times
C_{2}$, where $\chi (a)=-1$, $\chi (b)=1$, $\varphi \left( a\right) =1$, $%
\varphi \left( b\right) =-1$. Degree $2$ irreducible representations of $%
G_{4}$ are defined by \newline
\begin{tabular}{lll}
$\pi _{1}(a)=\left( 
\begin{array}{rr}
\omega  & 0 \\ 
0 & \omega ^{7}
\end{array}
\right) $ & 
$\pi _{2}(a)=\left( 
\begin{array}{rr}
i & 0 \\ 
0 & -i
\end{array}
\right) $ & 
$\pi _{3}(a)=\left( 
\begin{array}{rr}
\omega ^{3} & 0 \\ 
0 & \omega ^{5}
\end{array}
\right) $ \\
$\pi _{1}(b)=\left( 
\begin{array}{rr}
0 & -1 \\ 
1 & 0
\end{array}
\right) $ & $\pi _{2}(b)=\left( 
\begin{array}{rr}
0 & 1 \\ 
1 & 0
\end{array}
\right) $   & $\pi _{3}(b)=\left( 
\begin{array}{rr}
0 & -1 \\ 
1 & 0
\end{array}
\right) $%
\end{tabular}
\newline
where $\omega $ is a primitive $8^{th}$-root of unity, and representations
satisfy the properties \newline
\begin{tabular}{lll}
$\pi _{1}^{2}=1+\varphi +\pi _{2}=\pi _{3}^{2}$ & $\chi \bullet \pi _{1}=\pi
_{3}$ & $\chi \bullet \pi _{3}=\pi _{1}$ \\ 
$\pi _{2}^{2}=1+\chi +\varphi +\chi \varphi $ & $\varphi \bullet \pi
_{1}=\pi _{1}$ & $\varphi \bullet \pi _{3}=\pi _{3}$%
\end{tabular}
\smallskip

\item  $G_{5}=\left\langle a,b:a^{4}=b^{4}=1,ba=a^{-1}b\right\rangle $%
\newline
$\mathbf{G}\left( \left( kG_{5}\right) ^{\ast }\right) =\left\langle \chi
\right\rangle \times \left\langle \varphi \right\rangle \cong C_{4}\times
C_{2}$, where $\chi (a)=1$, $\chi (b)=i$, $\varphi \left( a\right) =-1$, $%
\varphi \left( b\right) =1$. Degree $2$ irreducible representations of $G_{5}
$ are defined by \newline
\begin{tabular}{ll}
$\pi _{1}(a)=\left( 
\begin{array}{rr}
0 & 1 \\ 
1 & 0
\end{array}
\right) $ & $\pi _{1}(b)=\left( 
\begin{array}{rr}
i & 0 \\ 
0 & -i
\end{array}
\right) $ \\ 
$\pi _{2}(a)=\left( 
\begin{array}{rr}
0 & i \\ 
i & 0
\end{array}
\right) $ & $\pi _{2}(b)=\left( 
\begin{array}{rr}
i & 0 \\ 
0 & -i
\end{array}
\right) $%
\end{tabular}
\newline
with the property $\pi _{1}^{2}=1+\chi ^{2}+\varphi +\chi ^{2}\varphi =\pi
_{2}^{2}.$
\smallskip

\item  $G_{6}=\left\langle a,b,c:a^{4}=b^{2}=c^{2}=1,bab=ac\right\rangle $%
\newline
$\mathbf{G}\left( \left( kG_{6}\right) ^{\ast }\right) =\left\langle \chi
\right\rangle \times \left\langle \varphi \right\rangle \cong C_{4}\times
C_{2}$, where $\chi (a)=i$, $\chi (b)=\chi (c)=1$, $\varphi \left( a\right)
=\varphi \left( c\right) =1$, $\varphi \left( b\right) =-1$. Degree $2$
irreducible representations of $G_{6}$ are defined by \newline
\begin{tabular}{lll}
$\pi _{1}(a)=\left( 
\begin{array}{rr}
0 & -1 \\ 
1 & 0
\end{array}
\right) $ & $\pi _{1}(b)=\left( 
\begin{array}{rr}
0 & 1 \\ 
1 & 0
\end{array}
\right) $ & $\pi _{1}(c)=\left( 
\begin{array}{rr}
-1 & 0 \\ 
0 & -1
\end{array}
\right) $ \\ 
$\pi _{2}(a)=\left( 
\begin{array}{rr}
0 & -i \\ 
i & 0
\end{array}
\right) $ & $\pi _{2}(b)=\left( 
\begin{array}{rr}
0 & 1 \\ 
1 & 0
\end{array}
\right) $ & $\pi _{2}(c)=\left( 
\begin{array}{rr}
-1 & 0 \\ 
0 & -1
\end{array}
\right) $%
\end{tabular}
\newline
with the property $\pi _{1}^{2}=1+\chi ^{2}+\varphi +\chi ^{2}\varphi =\pi
_{2}^{2}.$
\smallskip

\item  $G_{7}=\left\langle a,b,c:a^{4}=b^{2}=c^{2}=1,cbc=a^{2}b\right\rangle 
$\newline
$\mathbf{G}\left( \left( kG_{7}\right) ^{\ast }\right) =\left\langle \chi
\right\rangle \times \left\langle \varphi \right\rangle \times \left\langle
\psi \right\rangle \cong C_{2}\times C_{2}\times C_{2}$, where $\chi (a)=-1$%
, $\chi (b)=\chi (c)=1$, $\varphi \left( a\right) =\varphi \left( b\right)
=-1$, $\varphi \left( c\right) =1$, $\psi \left( a\right) =\psi \left(
b\right) =1$, $\psi \left( c\right) =-1$. Degree $2$ irreducible
representations of $G_{7}$ are defined by \newline
\begin{tabular}{lll}
$\pi _{1}(a)=\left( 
\begin{array}{rr}
i & 0 \\ 
0 & -i
\end{array}
\right) $ & $\pi _{1}(b)=\left( 
\begin{array}{rr}
-1 & 0 \\ 
0 & 1
\end{array}
\right) $ & $\pi _{1}(c)=\left( 
\begin{array}{rr}
0 & 1 \\ 
1 & 0
\end{array}
\right) $ \\ 
$\pi _{2}(a)=\left( 
\begin{array}{rr}
i & 0 \\ 
0 & -i
\end{array}
\right) $ & $\pi _{2}(b)=\left( 
\begin{array}{rr}
1 & 0 \\ 
0 & -1
\end{array}
\right) $ & $\pi _{2}(c)=\left( 
\begin{array}{rr}
0 & 1 \\ 
1 & 0
\end{array}
\right) $%
\end{tabular}
\newline
with the property $\pi _{1}^{2}=\chi +\chi \varphi +\chi \psi +\chi \varphi
\psi =\pi _{2}^{2}.$
\smallskip

\item  $G_{8}=\left\langle a,b,c:a^{4}=b^{2}=c^{2}=1,ba=a^{-1}b\right\rangle
=D_{8}\times C_{2}$\newline
$\mathbf{G}\left( \left( kG_{8}\right) ^{\ast }\right) =\left\langle \chi
\right\rangle \times \left\langle \varphi \right\rangle \times \left\langle
\psi \right\rangle \cong C_{2}\times C_{2}\times C_{2}$, where $\chi
(a)=\chi (b)=1$, $\chi (c)=-1$, $\varphi \left( a\right) =-1$, $\varphi
\left( b\right) =\varphi \left( c\right) =1$, $\psi \left( a\right) =\psi
\left( c\right) =1$, $\psi \left( b\right) =-1$. Degree $2$ irreducible
representations of $G_{8}$ are defined by \newline
\begin{tabular}{lll}
$\pi _{1}(a)=\left( 
\begin{array}{rr}
i & 0 \\ 
0 & -i
\end{array}
\right) $ & $\pi _{1}(b)=\left( 
\begin{array}{rr}
0 & 1 \\ 
1 & 0
\end{array}
\right) $ & $\pi _{1}(c)=\left( 
\begin{array}{rr}
1 & 0 \\ 
0 & 1
\end{array}
\right) $ \\ 
$\pi _{2}(a)=\left( 
\begin{array}{rr}
i & 0 \\ 
0 & -i
\end{array}
\right) $ & $\pi _{2}(b)=\left( 
\begin{array}{rr}
0 & 1 \\ 
1 & 0
\end{array}
\right) $ & $\pi _{2}(c)=\left( 
\begin{array}{rr}
-1 & 0 \\ 
0 & -1
\end{array}
\right) $%
\end{tabular}
\newline
with the property $\pi _{1}^{2}=1+\varphi +\psi +\varphi \psi =\pi _{2}^{2}$.
\smallskip

\item  $G_{9}=\left\langle
a,b,c:a^{4}=c^{2}=1,b^{2}=a^{2},ba=a^{-1}b\right\rangle =Q_{8}\times C_{2}$%
\newline
$\mathbf{G}\left( \left( kG_{9}\right) ^{\ast }\right) =\left\langle \chi
\right\rangle \times \left\langle \varphi \right\rangle \times \left\langle
\psi \right\rangle \cong C_{2}\times C_{2}\times C_{2}$, where $\chi
(a)=\chi (b)=1$, $\chi (c)=-1$, $\varphi \left( a\right) =-1$, $\varphi
\left( b\right) =\varphi \left( c\right) =1$, $\psi \left( a\right) =\psi
\left( c\right) =1$, $\psi \left( b\right) =-1$. Degree $2$ irreducible
representations of $G_{9}$ are defined by \newline
\begin{tabular}{lll}
$\pi _{1}(a)=\left( 
\begin{array}{rr}
i & 0 \\ 
0 & -i
\end{array}
\right) $ & $\pi _{1}(b)=\left( 
\begin{array}{rr}
0 & 1 \\ 
-1 & 0
\end{array}
\right) $ & $\pi _{1}(c)=\left( 
\begin{array}{rr}
1 & 0 \\ 
0 & 1
\end{array}
\right) $ \\ 
$\pi _{2}(a)=\left( 
\begin{array}{rr}
i & 0 \\ 
0 & -i
\end{array}
\right) $ & $\pi _{2}(b)=\left( 
\begin{array}{rr}
0 & 1 \\ 
-1 & 0
\end{array}
\right) $ & $\pi _{2}(c)=\left( 
\begin{array}{rr}
-1 & 0 \\ 
0 & -1
\end{array}
\right) $%
\end{tabular}
\newline
with the property $\pi _{1}^{2}=1+\varphi +\psi +\varphi \psi =\pi _{2}^{2}$.
\end{enumerate}

\section{Computations in $K_{0}(H)$ in the case of $\left| \mathbf{G}\left(
H^{\ast }\right) \right| =8$.} \label{sec5}

In this case we have $8$ one-dimensional irreducible representations $\chi
_{1}=1_{K_{0}(H)},\ldots ,$ $\chi _{8}\in \mathbf{G}\left( H^{\ast }\right) $
and two 2-dimensional ones $\pi _{1}$ and $\pi _{2}$. Then, since $\chi
_{i}\bullet \chi _{i}^{-1}=1_{K_{0}(H)}$, $\quad \chi _{i}^{\ast }=\chi
_{i}^{-1}$.

$\deg \left( \chi _{i}\bullet \pi _{k}\right) =2$ thus there are two
possibilities:

\begin{enumerate}
\item[i)]  $\chi _{i}\bullet \pi _{k}=\chi _{j}+\chi _{l}$

\item[ii)]  $\chi _{i}\bullet \pi _{k}=\pi _{l}$
\end{enumerate}

Case i) cannot happen, since otherwise $\pi _{k}=\chi _{i}^{-1}\bullet \chi
_{j}+\chi _{i}^{-1}\bullet \chi _{l}$ is not irreducible. Thus $\chi
_{i}\bullet \pi _{k}=\pi _{l}$. Then it is impossible to have $\chi
_{i}\bullet \pi _{k}=\pi _{k}$ and $\chi _{i}\bullet \pi _{l}=\pi _{k}$ for $%
k\neq l$, since otherwise $\pi _{l}=\chi _{i}^{8}\bullet \pi _{l}=\chi
_{i}^{7}\bullet \pi _{k}=\pi _{k}$. Thus $\chi _{i}$ either fixes both $\pi
_{1}$ and $\pi _{2}$ or interchanges them. It is easy to check that either
all $\chi _{i}$ fix $\pi _{k}$, $k=1,2$, or half of $\chi _{i}$ fixes $\pi
_{k}$ and half of $\chi _{i}$ interchanges them.

Suppose $\chi _{i}\bullet \pi _{k}=\pi _{k}$, for $i=1,\ldots ,8$ and $k=1,2$%
. Then 
\begin{equation*}
1=m\left( \pi _{k},\chi _{i}\bullet \pi _{k}\right) =m\left( \chi _{i}^{\ast
},\pi _{k}\bullet \pi _{k}^{\ast }\right) =m\left( \chi _{i}^{-1},\pi
_{k}\bullet \pi _{k}^{\ast }\right)
\end{equation*}
Thus 
\begin{equation*}
\pi _{k}\bullet \pi _{k}^{\ast }=\sum_{i=1}^{8}\chi
_{i}^{-1}=\sum_{i=1}^{8}\chi _{i}
\end{equation*}
but $\deg \left( \pi _{k}\bullet \pi _{k}^{\ast }\right) =4$ and $\deg
\left( \sum_{i=1}^{8}\chi _{i}\right) =8$. Thus
all $\chi _{i}$ cannot fix $\pi _{k}$.

Now let half of $\chi _{i}$ fix $\pi _{k}$ and half of $\chi _{i}$
interchange them, say 
\begin{eqnarray*}
\chi _{i}\bullet \pi _{k} &=&\pi _{k}\text{ for }i\text{ odd} \\
\chi _{i}\bullet \pi _{k} &=&\pi _{l}\text{ for }i\text{ even}   
\end{eqnarray*}
if $k\neq l$. It is clear that $\chi _{i}$ and $\chi _{i}^{\ast
}=\chi _{i}^{-1}$ fix or interchange $\pi _{k}$ simultaneously. Then for
$k\neq l$%
\begin{eqnarray*}
1 &=&m\left( \pi _{k},\chi _{i}\bullet \pi _{k}\right) =m\left( \chi
_{i}^{\ast },\pi _{k}\bullet \pi _{k}^{\ast }\right) =m\left( \chi
_{i}^{-1},\pi _{k}\bullet \pi _{k}^{\ast }\right) \text{ for }i\text{ odd} \\
1 &=&m\left( \pi _{l},\chi _{i}\bullet \pi _{k}\right) =m\left( \chi
_{i}^{\ast },\pi _{k}\bullet \pi _{l}^{\ast }\right) =m\left( \chi
_{i}^{-1},\pi _{k}\bullet \pi _{l}^{\ast }\right) \text{ for }i\text{ even}
\end{eqnarray*}
and therefore
\begin{eqnarray*}
\pi _{k}\bullet \pi _{k}^{\ast } &=&\chi _{1}+\chi _{3}+\chi _{5}+\chi _{7}
\\
\pi _{k}\bullet \pi _{l}^{\ast } &=&\chi _{2}+\chi _{4}+\chi _{6}+\chi _{8}
\end{eqnarray*}
There are two possibilities for the involution:
either $\pi _{1}^{\ast }=\pi _{1}$ and $\pi _{2}^{\ast }=\pi _{2}$,
or $\pi _{1}^{\ast }=\pi _{2}$ and
$\pi _{2}^{\ast }=\pi _{1}$.
It is easy to check that when $\mathbf{G}\left( H^{\ast }\right) $ is
isomorphic to $C_{2}\times C_{2}\times C_{2}$ or $Q_{8}$ it does not matter
which generators we choose to fix $\pi _{k}$ and which to interchange
them. In the case of $\mathbf{G}\left( H^{\ast }\right) \cong D_{8}$ or $%
C_{4}\times C_{2}$ it matters and should give us two more nonisomorphic
structures for $K_{0}(H)$ for each of them, but due to the results of the
Section \ref{sec3} we can see that in the case of $\mathbf{G}\left( H^{\ast
}\right) \cong C_{4}\times C_{2}$, $\pi _{k}$ can be fixed only by elements
of order $1$ or $2$ (since $\left( \pi _{k}\right) ^{2}$ is either the sum
of all elements of order $1$ or $2,$ or the sum of all elements of order $4$%
). 

Now assume that $\mathbf{G}\left( H^{\ast }\right) =\left\langle \chi ,
\varphi :\chi ^{4}=1,\varphi ^{2}=1,\varphi \chi =\chi ^{-1}\varphi 
\right\rangle \cong D_{8}$ or $\mathbf{G}\left( H^{\ast }\right) =\newline
\left\langle \chi ,
\varphi :\chi ^{4}=1,\varphi ^{2}=\chi ^{2},\varphi \chi =\chi ^{-1}\varphi 
\right\rangle \cong Q_{8}$. Then $\chi ^{2}$ is the only nontrivial central
element of $\mathbf{G}\left( H^{\ast }\right)$. Since by \cite[Theorem1]{Ma3} 
$H^{\ast }$ has a nontrivial central grouplike, this grouplike 
should be equal to $\chi ^{2}$. Since $\chi ^{2}$ is a central grouplike of 
order $2$, which fixes all $\pi _{k}$, by Proposition \ref{com_quot} $H$ has
a commutative subHopfalgebra of dimension $8$. Therefore,
it should have the same $K_{0}$-ring as one of the Hopf algebras described in
Section \ref{sec3}. Thus the only possible $K_{0}$-ring structure corresponds 
to $\mathbf{G}\left( H^{\ast }\right)\cong D_{8}$ with 
$\left( \pi _{k}\right) ^{2}=1+\chi^{2}+\varphi +\chi ^{2}\varphi $ and
$\pi _{i}^{\ast }=\pi _{i}$. 

Now let's list all the possible ring structures of $K_{0}(H)$. \medskip

\n $\mathbf{G}\left( H^{\ast }\right) =\left\langle \chi \right\rangle \times
\left\langle \varphi \right\rangle \times \left\langle \psi \right\rangle
\cong C_{2}\times C_{2}\times C_{2}$ where $\pi _{1}^{\ast }=\pi _{1}$ and $%
\pi _{2}^{\ast }=\pi _{2}$ and \label{5.1}\newline
\begin{tabular}{ll}
$\chi \bullet \pi _{1}=\pi _{2}=\pi _{1}\bullet \chi $ & $\psi \bullet \pi
_{1}=\pi _{1}=\pi _{1}\bullet \psi $ \\ 
$\chi \bullet \pi _{2}=\pi _{1}=\pi _{2}\bullet \chi $ & $\psi \bullet \pi
_{2}=\pi _{2}=\pi _{2}\bullet \psi $ \\ 
$\varphi \bullet \pi _{1}=\pi _{1}=\pi _{1}\bullet \varphi $ & $\pi
_{1}^{2}=1+\varphi +\psi +\varphi \psi =\pi _{2}^{2}$ \\ 
$\varphi \bullet \pi _{2}=\pi _{2}=\pi _{2}\bullet \varphi $ & $\pi
_{1}\bullet \pi _{2}=\chi +\chi \varphi +\chi \psi +\chi \varphi \psi =\pi
_{2}\bullet \pi _{1}$%
\end{tabular}\newline
Examples: $H_{b:1}$, $H_{d:1,1}$, $H_{d:-1,1}$, 
$k\left( D_{8}\times C_{2}\right) $, $k\left( Q_{8}\times
C_{2}\right) $, $H_{8}\otimes kC_{2}$. \medskip

\n $\mathbf{G}\left( H^{\ast }\right) =\left\langle \chi \right\rangle \times
\left\langle \varphi \right\rangle \times \left\langle \psi \right\rangle
\cong C_{2}\times C_{2}\times C_{2}$ where $\pi _{1}^{\ast }=\pi _{2}$ and $%
\pi _{2}^{\ast }=\pi _{1}$ and \label{5.2}\newline
\begin{tabular}{ll}
$\chi \bullet \pi _{1}=\pi _{2}=\pi _{1}\bullet \chi $ & $\psi \bullet \pi
_{1}=\pi _{1}=\pi _{1}\bullet \psi $ \\ 
$\chi \bullet \pi _{2}=\pi _{1}=\pi _{2}\bullet \chi $ & $\psi \bullet \pi
_{2}=\pi _{2}=\pi _{2}\bullet \psi $ \\ 
$\varphi \bullet \pi _{1}=\pi _{1}=\pi _{1}\bullet \varphi $ & $\pi
_{1}^{2}=\chi +\chi \varphi +\chi \psi +\chi \varphi \psi =\pi _{2}^{2}$ \\ 
$\varphi \bullet \pi _{2}=\pi _{2}=\pi _{2}\bullet \varphi $ & $\pi
_{1}\bullet \pi _{2}=1+\varphi +\psi +\varphi \psi =\pi _{2}\bullet \pi _{1}$%
\end{tabular}\newline
Examples: $H_{c:\sigma _{1}}$ and $kG_{7}$, where $%
G_{7}=\left\langle a,b,c:a^{4}=b^{2}=c^{2}=1,cbc=a^{2}b\right\rangle $.
\medskip

\n $\mathbf{G}\left( H^{\ast }\right) =\left\langle \chi \right\rangle \times
\left\langle \varphi \right\rangle \cong C_{4}\times C_{2}$ where $\pi
_{1}^{\ast }=\pi _{1}$ and $\pi _{2}^{\ast }=\pi _{2}$ and \label{5.3}\newline
\begin{tabular}{lll}
$\chi \bullet \pi _{1}=\pi _{2}=\pi _{1}\bullet \chi $ & $\varphi \bullet
\pi _{1}=\pi _{1}=\pi _{1}\bullet \varphi $ & $\pi _{1}^{2}=1+\chi
^{2}+\varphi +\chi ^{2}\varphi =\pi _{2}^{2}$ \\ 
$\chi \bullet \pi _{2}=\pi _{1}=\pi _{2}\bullet \chi $ & $\varphi \bullet
\pi _{2}=\pi _{2}=\pi _{2}\bullet \varphi $ & $\pi _{1}\bullet \pi _{2}=\chi
+\chi ^{3}+\chi \varphi +\chi ^{3}\varphi =\pi _{2}\bullet \pi _{1}$%
\end{tabular}\newline
Examples: $H_{a:1}$, $H_{a:y}$, $H_{b:y}$, $H_{b:x^{2}y}$, $H_{d:1,-1}$, $%
H_{d:-1,-1}$, $kG_{5}$ and $kG_{6}$, where \newline $%
G_{5}=\left\langle a,b:a^{4}=b^{4}=1,b^{-1}ab=a^{-1}\right\rangle $ and 
$G_{6}=\left\langle a,b,c:a^{4}=b^{2}=c^{2}=1,bab=ac\right\rangle $.
\medskip

\n $\mathbf{G}\left( H^{\ast }\right) =\left\langle \chi \right\rangle \times
\left\langle \varphi \right\rangle \cong C_{4}\times C_{2}$ where $\pi
_{1}^{\ast }=\pi _{2}$ and $\pi _{2}^{\ast }=\pi _{1}$ and \label{5.4}\newline
\begin{tabular}{lll}
$\chi \bullet \pi _{1}=\pi _{2}=\pi _{1}\bullet \chi $ & $\varphi \bullet
\pi _{1}=\pi _{1}=\pi _{1}\bullet \varphi $ & $\pi _{1}^{2}=1+\chi
^{2}+\varphi +\chi ^{2}\varphi =\pi _{2}^{2}$ \\ 
$\chi \bullet \pi _{2}=\pi _{1}=\pi _{2}\bullet \chi $ & $\varphi \bullet
\pi _{2}=\pi _{2}=\pi _{2}\bullet \varphi $ & $\pi _{1}\bullet \pi _{2}=\chi
+\chi ^{3}+\chi \varphi +\chi ^{3}\varphi =\pi _{2}\bullet \pi _{1}$%
\end{tabular}\newline
Examples: $H_{c:\sigma _{0}}$, $kG_{1}$, where $G_{1}=\left\langle
a,b:a^{8}=b^{2}=1,bab=a^{5}\right\rangle $.
\medskip

\n $\mathbf{G}\left( H^{\ast }\right) =\left\langle \chi ,\varphi :\chi
^{4}=1,\varphi ^{2}=1,\varphi \chi =\chi ^{-1}\varphi \right\rangle \cong
D_{8}$, where $\pi _{1}^{\ast }=\pi _{1}$ and $\pi _{2}^{\ast }=\pi _{2}$
and \label{5.5}\newline
\begin{tabular}{lll}
$\chi \bullet \pi _{1}=\pi _{2}=\pi _{1}\bullet \chi $ & $\varphi \bullet
\pi _{1}=\pi _{1}=\pi _{1}\bullet \varphi $ & $\pi _{1}^{2}=1+\chi
^{2}+\varphi +\chi ^{2}\varphi =\pi _{2}^{2}$ \\ 
$\chi \bullet \pi _{2}=\pi _{1}=\pi _{2}\bullet \chi $ & $\varphi \bullet
\pi _{2}=\pi _{2}=\pi _{2}\bullet \varphi $ & $\pi _{1}\bullet \pi _{2}=\chi
+\chi ^{3}+\chi \varphi +\chi ^{3}\varphi =\pi _{2}\bullet \pi _{1}$%
\end{tabular}\newline
Examples: $H_{B:1}$, $H_{B:X}$ and $kQ_{8}\#^{\alpha }kC_{2}$.

\begin{remark}
Noncommutative $K_{0}(H)$ should have the structure \ref{sec5}.\ref{5.5}.
\end{remark}

\section{Computations in $K_{0}(H)$ in the case of $\left| \mathbf{G}\left(
H^{\ast }\right) \right| =4$.} \label{sec6}

In this case by Theorem \ref{th1} $\mathbf{G}\left( H^{\ast }\right) \cong
C_{2} \times C_{2}$ and
we have $4$ one-dimensional irreducible representations $\chi
_{1}=1_{K_{0}(H)},\ldots ,$ $\chi _{4}\in \mathbf{G}\left( H^{\ast }\right) $
and three two-dimensional ones $\pi _{1}$, $\pi _{2}$ and $\pi _{3}$. Then,
since $\chi _{i}\bullet \chi _{i}=1_{K_{0}(H)}$, $\quad \chi _{i}^{\ast
}=\chi _{i}$. The involution is an antihomomorphism of $K_{0}(H)$ of
order $2$, thus it either fixes all $\pi _{k}$ or interchanges two of them
and fixes the third one. Assume that we always have $\pi _{2}^{\ast }=\pi
_{2}$.

$\chi _{i}\bullet \pi _{k}\neq \chi _{j}+\chi _{l}$ as in the case of $%
\left| \mathbf{G}\left( H^{\ast }\right) \right| =8$. Thus multiplication by 
$\chi _{i}$ permutes $\pi _{k}$. Since $o\left( \chi _{i}\right) =1$ or $2$ 
then each $\chi _{i}$ either fixes all $\pi _{k}$ or interchanges two of
them and fixes the third one. There are two possible cases:

\begin{enumerate}
\item[i)]  $\chi _{i}\bullet \pi _{k}=\pi _{k}$ for $i=1,\ldots ,4$ and $%
k=1,2,3$. Then
\begin{eqnarray}
m\left( \chi _{i},\pi _{k}\bullet \pi _{k}^{\ast }\right)&=&
m\left( \pi _{k},\chi _{i}\bullet \pi _{k}\right) =1\ \notag \\
\pi _{k}\bullet \pi _{k}^{\ast }&=&\sum_{i=1}^{4}\chi _{i}\ \text{ for }
k=1,2,3 \label{fix}
\end{eqnarray}
By \cite[Theorem 1]{Ma3}, one of the $\chi _{i}$ is central of order $2$ and 
therefore by Proposition
\ref{com_quot}, $H$ has a commutative subHopfalgebra of order $8$. Therefore,
it should have the same $K_{0}$-ring as one of the Hopf algebras described in
Section \ref{sec3}. But none of these $K_{0}$-rings satisfies (\ref{fix}).
Therefore this case is not possible.
\item[ii)]  $\chi _{i}\bullet \pi _{k}\neq \pi _{k}$ for some $i\in \left\{
1,\ldots ,4\right\} $ and $k\in \left\{ 1,2,3\right\} $. Then, say, 
\[
\chi _{1}\bullet \pi _{k}=\chi _{3}\bullet \pi _{k}=\pi _{k}\text{ for
}k=1,2,3 
\]
but $\chi _{2}\bullet \pi _{k}\neq \pi _{k}$, $\chi _{4}\bullet \pi _{k}\neq
\pi _{k}$ for some $k\in \left\{ 1,2,3\right\} $.

Assume that $\pi _{1}^{\ast }=\pi _{3}$, $\quad \pi _{3}^{\ast }=\pi _{1}$
and $\chi _{2}\bullet \pi _{2}\neq \pi _{2}$. Then 
\begin{eqnarray*}
1&=&m\left( \pi _{2},\chi _{i}\bullet \pi _{2}\right) =m\left( \chi _{i},\pi
_{2}\bullet \pi _{2}\right) \text{ for }i=1,3 \\
0&=&m\left( \pi _{2},\chi _{i}\bullet \pi _{2}\right) =m\left( \chi _{i},\pi
_{2}\bullet \pi _{2}\right) \text{ for }i=2,4 
\end{eqnarray*}
Therefore 
\[
\pi _{2}\bullet \pi _{2}^{\ast }=\chi _{1}+\chi _{3}+\pi _{r}=\pi
_{2}\bullet \pi _{2}  
\]
Since $\left( \pi _{k}\bullet \pi _{k}^{\ast }\right) ^{\ast }=\pi
_{k}\bullet \pi _{k}^{\ast }$ , we get $\left( \pi _{r}\right) ^{\ast }=\pi
_{r}$ and and thus $\pi _{r}=\pi _{2}$, that is 
\[
\pi _{2}\bullet \pi _{2}^{\ast }=\chi _{1}+\chi _{3}+\pi _{2}  
\]

Therefore $R=\left\{ a\chi _{1} + b\chi _{3} +c\pi_{2} \in K_{0}(H):
a,b,c \in \mathbb{Z} \right\} $ is a hereditary subring of $K_{0}(H)$ (see
\cite[Definition 3.10]{N1}). There is a one-to-one correspondence between
hereditary subrings of $K_{0}(H)$ and Hopf quotients of $H$, that is between
hereditary subrings of $K_{0}(H)$ and subHopfalgebras of $H^{\ast }$ (see 
\cite[Theorem 6]{NR} or \cite[Proposition 3.11]{N1}). Thus $H^{\ast }$ has a
subHopfalgebra of dimension $1+1+4=6$, which contradicts Nichols-Zoeller
Theorem \cite{NZ}. 

Thus without loss of generality
$\chi _{2}\bullet \pi _{2}=\pi _{2}$. Then 
\begin{eqnarray*}
\chi _{i}\bullet \pi _{2}&=&\pi _{2}\text{ for }i=1,\ldots ,4 \\
\chi _{i}\bullet \pi _{1}&=&\pi _{1}\text{ for }i=1,3 \\
\chi _{i}\bullet \pi _{1}&=&\pi _{3}\text{ for }i=2,4 \\
\chi _{i}\bullet \pi _{3}&=&\pi _{3}\text{ for }i=1,3 \\
\chi _{i}\bullet \pi _{3}&=&\pi _{1}\text{ for }i=2,4
\end{eqnarray*}
and therefore 
\begin{eqnarray*}
1&=&m\left( \pi _{2},\chi _{i}\bullet \pi _{2}\right) =m\left( \chi _{i},\pi
_{2}\bullet \pi _{2}\right) \text{ for }i=1,\ldots ,4 \\
0&=&m\left( \pi _{2},\chi _{i}\bullet \pi _{k}\right) =m\left( \chi _{i},\pi
_{k}\bullet \pi _{2}^{\ast }\right) \text{ for }i=1,\ldots ,4\text{, }k\neq 2
\\
1&=&m\left( \pi _{k},\chi _{i}\bullet \pi _{k}\right) =m\left( \chi _{i},\pi
_{k}\bullet \pi _{k}^{\ast }\right) \text{ for }i=1,3 \\
0&=&m\left( \pi _{k},\chi _{i}\bullet \pi _{k}\right) =m\left( \chi _{i},\pi
_{k}\bullet \pi _{k}^{\ast }\right) \text{ for }i=2,4 \\
0&=&m\left( \pi _{3},\chi _{i}\bullet \pi _{1}\right) =m\left( \chi _{i},\pi
_{1}\bullet \pi _{3}^{\ast }\right) \text{ for }i=1,3 \\
1&=&m\left( \pi _{3},\chi _{i}\bullet \pi _{1}\right) =m\left( \chi _{i},\pi
_{1}\bullet \pi _{3}^{\ast }\right) \text{ for }i=2,4
\end{eqnarray*}
Thus we get 
\begin{eqnarray*}
\pi _{2}\bullet \pi _{2}^{\ast }&=&\sum_{i=1}^{4}\chi _{i} \\
\pi _{1}\bullet \pi _{1}^{\ast }&=&\chi _{1}+\chi _{3}+\pi _{2}=\pi
_{3}\bullet \pi _{3}^{\ast } \\
\pi _{1}\bullet \pi _{3}^{\ast }&=&\chi _{2}+\chi _{4}+\pi _{t} \\
\pi _{k}\bullet \pi _{2}^{\ast }&=&\alpha _{1}\pi _{1}+\alpha _{2}\pi
_{2}+\alpha _{3}\pi _{3}\text{ for }k\neq 2
\end{eqnarray*}
and then 
\begin{eqnarray*}
1&=&m\left( \pi _{2},\pi _{k}^{\ast }\bullet \pi _{k}\right) =m\left( \pi
_{k},\pi _{k}\bullet \pi _{2}^{\ast }\right) \text{ for }k\neq 2 \\
0&=&m\left( \pi _{k},\pi _{l}^{\ast }\bullet \pi _{l}\right) =m\left( \pi
_{l},\pi _{l}\bullet \pi _{k}^{\ast }\right) \text{ for }k\neq 2
\end{eqnarray*}
Therefore 
\begin{eqnarray*}
\pi _{1}\bullet \pi _{2}^{\ast }&=&\pi _{1}+\pi _{3} \\
\pi _{3}\bullet \pi _{2}^{\ast }&=&\pi _{1}+\pi _{3}
\end{eqnarray*}
and 
\begin{equation*}
1=m\left( \pi _{3},\pi _{1}\bullet \pi _{2}^{\ast }\right) =m\left( \pi
_{2},\pi _{3}^{\ast }\bullet \pi _{1}\right) =m\left( \pi _{2},\pi
_{1}\bullet \pi _{3}^{\ast }\right) 
\end{equation*}
So, finally, 
\begin{equation*}
\pi _{1}\bullet \pi _{3}^{\ast }=\chi _{2}+\chi _{4}+\pi _{2}
\end{equation*}
\end{enumerate}

Now let's list all the possible ring structures of $K_{0}(H)$. \medskip

\n $\mathbf{G}\left( H^{\ast }\right) =\left\langle \chi \right\rangle \times
\left\langle \varphi \right\rangle \cong C_{2}\times C_{2}$ where $\pi
_{1}^{\ast }=\pi _{1}$, $\pi _{2}^{\ast }=\pi _{2}$ and $\pi _{3}^{\ast
}=\pi _{3}$ and \label{6.3}\newline
\begin{tabular}{lll}
$\chi \bullet \pi _{1}=\pi _{3}=\pi _{1}\bullet \chi $ & $\varphi \bullet
\pi _{1}=\pi _{1}=\pi _{1}\bullet \varphi $ & $\pi _{1}\bullet \pi _{2}=\pi
_{1}+\pi _{3}=\pi _{2}\bullet \pi _{1}$ \\ 
$\chi \bullet \pi _{2}=\pi _{2}=\pi _{2}\bullet \chi $ & $\varphi \bullet
\pi _{2}=\pi _{2}=\pi _{2}\bullet \varphi $ & $\pi _{1}\bullet \pi _{3}=\chi
+\chi \varphi +\pi _{2}=\pi _{3}\bullet \pi _{1}$ \\ 
$\chi \bullet \pi _{3}=\pi _{1}=\pi _{3}\bullet \chi $ & $\varphi \bullet
\pi _{3}=\pi _{3}=\pi _{3}\bullet \varphi $ & $\pi _{2}\bullet \pi _{3}=\pi
_{1}+\pi _{3}=\pi _{3}\bullet \pi _{2}$ \\ 
$\pi _{1}^{2}=1+\varphi +\pi _{2}=\pi _{3}^{2}$ & $\pi _{2}^{2}=1+\chi
+\varphi +\chi \varphi $ & 
\end{tabular}\newline
Examples: $H_{C:1}$, $H_{C:\sigma _{1}},$ $kD_{16}$ and $kQ_{16}$.\medskip

\n $\mathbf{G}\left( H^{\ast }\right) =\left\langle \chi \right\rangle \times
\left\langle \varphi \right\rangle \cong C_{2}\times C_{2}$ where $\pi
_{1}^{\ast }=\pi _{3}$, $\pi _{2}^{\ast }=\pi _{2}$ and $\pi _{3}^{\ast
}=\pi _{1}$ and \label{6.4} \newline
\begin{tabular}{lll}
$\chi \bullet \pi _{1}=\pi _{3}=\pi _{1}\bullet \chi $ & $\varphi \bullet
\pi _{1}=\pi _{1}=\pi _{1}\bullet \varphi $ & $\pi _{1}\bullet \pi _{2}=\pi
_{1}+\pi _{3}=\pi _{2}\bullet \pi _{1}$ \\ 
$\chi \bullet \pi _{2}=\pi _{2}=\pi _{2}\bullet \chi $ & $\varphi \bullet
\pi _{2}=\pi _{2}=\pi _{2}\bullet \varphi $ & $\pi _{1}\bullet \pi
_{3}=1+\varphi +\pi _{2}=\pi _{3}\bullet \pi _{1}$ \\ 
$\chi \bullet \pi _{3}=\pi _{1}=\pi _{3}\bullet \chi $ & $\varphi \bullet
\pi _{3}=\pi _{3}=\pi _{3}\bullet \varphi $ & $\pi _{2}\bullet \pi _{3}=\pi
_{1}+\pi _{3}=\pi _{3}\bullet \pi _{2}$ \\ 
$\pi _{1}^{2}=\chi +\chi \varphi +\pi _{2}=\pi _{3}^{2}$ & $\pi
_{2}^{2}=1+\chi +\varphi +\chi \varphi $ & 
\end{tabular} \newline
Examples: $H_{E}$ and $kG_{2}$, where $G_{2}=\left\langle
a,b:a^{8}=b^{2}=1,bab=a^{3}\right\rangle $.

We can now prove \textbf{Theorem \ref{th4}}:

\begin{proof}
In Sections \ref{sec5} and \ref{sec6} we have described all possible
Grothendieck ring structures of non-commutative semisimple Hopf algebras of
dimension $16$ and there are exactly $7$ of them. Only one of these 
$K_{0}$-rings is not
commutative, namely $K_{\ref{sec5}.\ref{5.5}}$, which corresponds to 
nonabelian $\mathbf{G}\left( H^{\ast }\right) \cong D_{8}$. Therefore, by
\cite[Theorem 4.1]{N1} all Hopf algebras with non-commutative $K_{0}$-ring
are twistings of each other with a $2$-pseudo-cocycle. Moreover,
by \cite[4.5]{N1}, Hopf algebras with non-commutative $K_{0}$-rings are 
not twistings of group algebras.

If $\mathbf{G}\left( H^{\ast }\right)$ is abelian then there are $6$ 
possibilities for the $K_{0}$-ring structure, all of which are commutative, 
namely $K_{\ref{sec5}.\ref{5.1}}=K_{0}\left( k\left( D_{8}\times C_{2}\right)
\right) $, $K_{\ref{sec5}.\ref{5.2}}=K_{0}\left( kG_{7}\right) $, $%
K_{\ref{sec5}.\ref{5.3}}=K_{0}\left( kG_{5}\right) $, 
$K_{\ref{sec5}.\ref{5.4}}=K_{0}\left( kG_{1}\right) $, 
$K_{\ref{sec6}.\ref{6.3}}=K_{0}\left( kD_{16}\right) $ and
$K_{\ref{sec6}.\ref{6.4}}=K_{0}\left( kG_{2}\right) $. Thus
by \cite[Theorem 4.1]{N1} $H$ is
a twisting of one of these group algebras with a $2$-pseudo-cocycle.
Since $H$ is semisimple, 
$K_{0}\left( H\right) \otimes _{\mathbb{Z}}k$ is also semisimple by
\cite[Lemma 2]{Z}. Therefore, if $K_{0}\left( H\right)$ is commutative,
as algebras $K_{0}\left( H\right) \otimes _{\mathbb{Z}}k\cong k^{\left( 
10\right) }$ when $\left| \mathbf{G}\left( H^{\ast }\right) \right| =8$ and
$K_{0}\left( H\right) \otimes _{\mathbb{Z}}k\cong k^{\left( 7\right) }$
when $\left| \mathbf{G}\left( H^{\ast }\right) \right| =4$. If 
$K_{0}\left( H\right)$ is not commutative, that is $K_{0}\left( H\right)=
K_{\ref{sec5}.\ref{5.5}}$, it is easy to see that $\dim Z\left( K_{0}\left( 
H\right)\right) = 7$ and thus $K_{0}\left( H\right)\otimes _{\mathbb{Z}}k
\cong k^{\left(6\right) }\oplus M_{2}\left( k\right) $.
\end{proof}
 
\section{Twistings of group algebras with a $2$-cocycle.} \label{seccoc}

All nonabelian groups $G$, considered in this section, have an abelian
subgroup $F=\left\{ 1,c,b,cb\right\} \cong C_{2}\times C_{2}$. $kF\cong
\left( kF\right) ^{\ast }$ thus we can identify $\delta _{x}\in kF$ with the
elements of the dual basis. Now define $J\in kF\otimes kF$ as follows: 
\begin{eqnarray}
J &=&\delta _{1}\otimes \delta _{1}+\delta _{1}\otimes \delta _{c}+\delta
_{1}\otimes \delta _{b}+\delta _{1}\otimes \delta _{cb}  \notag \\
&&+\delta _{c}\otimes \delta _{1}+\delta _{c}\otimes \delta _{c}+i\delta
_{c}\otimes \delta _{b}-i\delta _{c}\otimes \delta _{cb}  \notag \\
&&+\delta _{b}\otimes \delta _{1}-i\delta _{b}\otimes \delta _{c}+\delta
_{b}\otimes \delta _{b}+i\delta _{b}\otimes \delta _{cb}  \notag \\
&&+\delta _{cb}\otimes \delta _{1}+i\delta _{cb}\otimes \delta _{c}-i\delta
_{cb}\otimes \delta _{b}+\delta _{cb}\otimes \delta _{cb}  \label{omega}
\end{eqnarray}
where 
\begin{eqnarray*}
\delta _{1} &=&\frac{1}{4}\left( 1+c+b+cb\right) \\
\delta _{c} &=&\frac{1}{4}\left( 1+c-b-cb\right) \\
\delta _{b} &=&\frac{1}{4}\left( 1-c+b-cb\right) \\
\delta _{cb} &=&\frac{1}{4}\left( 1-c-b+cb\right)
\end{eqnarray*}
We can rewrite $J$ as 
\begin{eqnarray*}
J &=&\frac{1}{8}(5\cdot 1\otimes 1+c\otimes 1+b\otimes 1+cb\otimes 1 \\
&&+1\otimes c+c\otimes c+\left( -1-2i\right) b\otimes c+\left( -1+2i\right)
cb\otimes c \\
&&+1\otimes b+\left( -1+2i\right) c\otimes b+b\otimes b+\left( -1-2i\right)
cb\otimes b \\
&&+1\otimes cb+\left( -1-2i\right) c\otimes cb+\left( -1+2i\right) b\otimes
cb+cb\otimes cb)
\end{eqnarray*}

Such a $J$ is a 2-cocycle for $kF$ and since $J\in kG\otimes kG$, it is also
a 2-cocycle for $kG$. Thus we can form $\left( kG\right) _{J}$ which is a
Hopf algebra by \cite{V}, \cite[2.8]{N1}. By \cite[6.4]{V}, $\left(
kG\right) _{J}$ is non-cocommutative if and only if $J^{-1}\left( \tau
J\right) $ does not lie in the centralizer of $\Delta \left( kG\right) $ in $%
kG\otimes kG$. Moreover, by \cite[Theorem 4.1]{N1} $K_{0}\left( \left(
kG\right) _{J}\right) \cong K_{0}\left( kG\right) $. Since $J$ is a
2-cocycle, then by \cite{EG} $\left( kG\right) _{J}$ is triangular.

We now discuss Examples $2$, $12$ and $13$ in the tables. We used GAP to
compute $\mathbf{G}\left( H\right) $ in Examples $12$ and $13$.

\noindent \textbf{Example }$\mathbf{2}$\textbf{. }$H=\left( k\left(
D_{8}\times C_{2}\right) \right) _{J}$, where $F=\left\{ 1,c,b,cb\right\}
\cong C_{2}\times C_{2}$ is a subgroup of $D_{8}\times C_{2}=\left\langle
a,b,c:a^{4}=b^{2}=c^{2}=1,ba=a^{-1}b\right\rangle $ and $J$ is given by the
formula $\left( \ref{omega}\right) $. Then

\begin{enumerate}
\item  $\mathbf{G}\left( H\right) =\left\langle a^{2},b,c\right\rangle \cong
C_{2}\times C_{2}\times C_{2}$.

\item  $\mathbf{G}\left( H^{\ast }\right) \cong \mathbf{G}\left( k\left(
D_{8}\times C_{2}\right) \right) \cong C_{2}\times C_{2}\times C_{2}$ and $%
K_{0}\left( H\right) \cong K_{0}\left( k\left( D_{8}\times C_{2}\right)
\right) \cong K_{\ref{sec5}.\ref{5.1}}$.
\end{enumerate}

\noindent \textbf{Example }$\mathbf{12}$\textbf{. }$H=\left( kD_{16}\right)
_{J}$, where $F=\left\{ 1,a^{4},b,a^{4}b\right\} \cong C_{2}\times C_{2}$ is
a subgroup of \newline
$D_{16}=\left\langle a,b:a^{8}=b^{2}=1,ba=a^{-1}b\right\rangle $ and $J$ is
given by the formula $\left( \ref{omega}\right) $. Then

\begin{enumerate}
\item  $\mathbf{G}\left( H\right) =\left\langle b,g=\frac{1}{2}\left(
-a^{2}+a^{2}b+a^{6}+a^{6}b\right) :g^{4}=b^{2}=1,bgb=g^{-1}\right\rangle
\cong D_{8}$.

\item  $\mathbf{G}\left( H^{\ast }\right) \cong \mathbf{G}\left( \left(
kD_{16}\right) ^{\ast }\right) \cong C_{2}\times C_{2}$ and $K_{0}\left(
H\right) \cong K_{0}\left( kD_{16}\right) \cong K_{\ref{sec6}.\ref{6.3}}$.
\end{enumerate}

\noindent \textbf{Example }$\mathbf{13}$\textbf{. }$H=\left( kG_{2}\right)
_{J}$, where $F=\left\{ 1,a^{4},b,a^{4}b\right\} \cong C_{2}\times C_{2}$ is
a subgroup of $G_{2}=\left\langle a,b:a^{8}=b^{2}=1,ba=a^{3}b\right\rangle $
and $J$ is given by the formula $\left( \ref{omega}\right) $. Then

\begin{enumerate}
\item  $\mathbf{G}\left( H\right) \cong D_{8}$.

\item  $\mathbf{G}\left( H^{\ast }\right) \cong \mathbf{G}\left( \left(
kG_{2}\right) ^{\ast }\right) \cong C_{2}\times C_{2}$ and $K_{0}\left(
H\right) \cong K_{0}\left( kG_{2}\right) \cong K_{\ref{sec6}.\ref{6.4}}$.
\end{enumerate}

\section{A construction using smash coproducts.} \label{seccoprod}

Let\textbf{\ }$%
H=kQ_{8}\#^{\alpha }kC_{2}$, a smash coproduct of $kQ_{8}$ and $kC_{2}$ (see 
\cite[10.6.1]{M} or \cite[Proposition 3.8]{N1}), where $Q_{8}=\left\langle
a,b:a^{4}=1,b^{2}=a^{2},ba=a^{-1}b\right\rangle $, $C_{2}=\left\{
1,g\right\} $. $H$ has the algebra structure of $kQ_{8}\otimes kC_{2}$ and
the following comultiplication, antipode and counit: 
\begin{eqnarray*}
\Delta \left( x\#\delta _{g^{k}}\right) &=&\underset{r+t=k}{\sum }\left(
x_{1}\#\delta _{g^{r}}\right) \otimes \left( \alpha _{g^{r}}\left(
x_{2}\right) \#\delta _{g^{t}}\right) \\
S\left( x\#\delta _{g^{k}}\right) &=&\alpha _{g^{k}}\left( S\left( x\right)
\right) \#\delta _{g^{-k}} \\
\varepsilon \left( x\#\delta _{g^{k}}\right) &=&\varepsilon \left( x\right)
\delta _{g^{k},1}
\end{eqnarray*}
where $\delta _{1}=\frac{1}{2}\left( 1+g\right) $, $\delta _{g}=\frac{1}{2}%
\left( 1-g\right) $, $x\in kQ_{8}$ and $\alpha :G\rightarrow Aut\left(
kQ_{8}\right) $ is defined by 
\begin{eqnarray*}
\alpha _{1}\left( x\right) &=&x \\
\alpha _{g}\left( a\right) &=&b \\
\alpha _{g}\left( b\right) &=&a,
\end{eqnarray*}
see \cite[Erratum]{N1}. It follows from the above that 
\begin{eqnarray}
\Delta \left( a\#1\right) &=&\Delta \left( a\#\delta _{1}\right) +\Delta
\left( a\#\delta _{g}\right)  \notag \\
&=&a\#\delta _{1}\otimes a\#\delta _{1}+a\#\delta _{g}\otimes b\#\delta
_{g}+a\#\delta _{1}\otimes a\#\delta _{g}+a\#\delta _{g}\otimes b\#\delta
_{1}  \notag \\
&=&a\#\delta _{1}\otimes a\#1+a\#\delta _{g}\otimes b\#1  \notag \\
&=&
{\frac12}
\left( a\#1\otimes a\#1+a\#g\otimes a\#1+a\#1\otimes b\#1-a\#g\otimes
b\#1\right)  \label{del_a}
\end{eqnarray}
and 
\begin{equation}
\Delta \left( b\#1\right) =
{\frac12}
\left( b\#1\otimes b\#1+b\#g\otimes b\#1+b\#1\otimes a\#1-b\#g\otimes
a\#1\right)  \label{del_b}
\end{equation}
Let's describe $\mathbf{G}\left( H\right) $, $\mathbf{G}\left( H^{\ast
}\right) $ and $K_{0}\left( H\right) $. By straightforward computations,
using $\left( \ref{del_a}\right) $ and $\left( \ref{del_b}\right) $, $%
\mathbf{G}\left( H\right) =\left\langle a^{2}\#1\right\rangle \times
\left\langle 1\#g\right\rangle $. %
$\mathbf{G}\left( H^{\ast }\right) \ $is generated by the multiplicative
characters $\chi $ and $\varphi $, defined by $\chi (a)=\chi (g)=-1$, $\chi
(b)=1$ and $\varphi \left( a\right) =-1$, $\varphi \left( g\right) =\varphi
\left( b\right) =1$. Then $\chi ^{-1}(a)=1$, $\chi ^{-1}(b)=\chi ^{-1}(g)=-1$
and $\varphi \chi \varphi =\chi ^{-1}$. Therefore $\mathbf{G}\left( H^{\ast
}\right) \cong D_{8}$. Degree $2$ irreducible representations of $H$ are
defined by \newline
\begin{tabular}{lll}
$\pi _{1}(a)=\left( 
\begin{array}{rr}
i & 0 \\ 
0 & -i
\end{array}
\right) $ & $\pi _{1}(b)=\left( 
\begin{array}{rr}
0 & 1 \\ 
-1 & 0
\end{array}
\right) $ & $\pi _{1}(g)=\left( 
\begin{array}{rr}
1 & 0 \\ 
0 & 1
\end{array}
\right) $ \\ 
$\pi _{2}(a)=\left( 
\begin{array}{rr}
i & 0 \\ 
0 & -i
\end{array}
\right) $ & $\pi _{2}(b)=\left( 
\begin{array}{rr}
0 & 1 \\ 
-1 & 0
\end{array}
\right) $ & $\pi _{2}(g)=\left( 
\begin{array}{rr}
-1 & 0 \\ 
0 & -1
\end{array}
\right) $%
\end{tabular}
\newline
with the property $\pi _{1}^{2}=\pi _{2}^{2}=1+\chi ^{2}+\varphi +\chi
^{2}\varphi $ and therefore $K_{0}\left( H\right) \cong K_{\ref{sec5}.\ref{5.5}}$.

\section{Main results.} \label{summary}
\begin{theorem} \label{com_sub}
Every nontrivial semisimple Hopf algebra $H$ of dimension $16$ has a 
commutative subHopfalgebra of dimension $8$.
\end{theorem}
\proof
By \cite[Theorem 1]{Ma3} $H^{\ast }$ has a central grouplike $g$ of order $2$.
Thus we get a short exact sequence of Hopf algebras 
\begin{equation}
k\left\langle g\right\rangle\overset{i}{\hookrightarrow }H^{\ast }
\overset{\pi }{\twoheadrightarrow }K
\end{equation}
If $K$ is cocommutative, $K^{\ast } \subset H$ is commutative and we are 
done. If $K$ is commutative, but not cocommutative, then 
$K^{\ast }=k\mathbf{G}\left( H\right)\subset H$ and $\mathbf{G}\left( H\right)$
is nonabelian of order $8$. Applying Proposition \ref{com_quot} and results of
Section \ref{sec5} we get that $H^{\ast}$ has a 
commutative subHopfalgebra of dimension $8$. Therefore $H^{\ast}$ was described
in Section \ref{3.3}, case B), and it has a group algebra of dimension $8$ 
as a quotient. Therefore $H$ has a commutative subHopfalgebra of dimension $8$.

Now assume that $K$ is neither commutative nor cocommutative. Then 
$K\cong K^{\ast }\cong H_{8}$ and as algebras 
$H^{\ast} \cong k\left\langle g\right\rangle\# _{\sigma }K$, a 
crossed product of Hopf algebras with an action  $
\rightharpoonup :K\otimes k\left\langle g\right\rangle\to 
k\left\langle g\right\rangle$ and a cocycle 
$\sigma :K\otimes K\to k\left\langle g\right\rangle$. Since $g$ is central 
in $H^{\ast}$, the action
$\rightharpoonup $ is trivial. Let's prove that $\sigma $ is also trivial.
By \cite[Theorem 4.8]{Ma7} $K\cong H_{8}$ doesn't have nontrivial right Galois 
objects, thus for any 
$2$-cocycle $\alpha :K\otimes K \to k$ the crossed product 
$_{\alpha }K=k\# _{\alpha }K$ is trivial, that is $ab=\sum \alpha
\left(a_{1}, b_{1}\right) a_{2}b_{2}$. Write
$k\left\langle g\right\rangle=k\frac{1+g}{2}\oplus k\frac{1-g}{2}$ and 
$\sigma \left( 
a\otimes b\right) =\alpha _{1}\left( a\otimes b\right)\frac{1+g}{2}+
\alpha _{2}\left( a\otimes b\right)\frac{1-g}{2}$. Then for $j=1,2$, 
$\alpha _{j}:K\otimes K \to k$ are $2$-cocycles and 
\begin{eqnarray*}
\left(1\# a\right)\left(1\# b\right)&=&\sum \sigma \left(a_{1}, b_{1}\right) 
\# a_{2}b_{2}
=\frac{1+g}{2}\# \sum \alpha _{1} \left(a_{1}, b_{1}\right) 
a_{2}b_{2} \\
&&+\frac{1-g}{2}\#\sum \alpha _{2} \left(a_{1}, b_{1}\right) a_{2}b_{2}
=\frac{1+g}{2}\# ab +\frac{1-g}{2}\# ab = 1\# ab
\end{eqnarray*}
Therefore as algebras $H^{\ast} \cong 
k\left\langle g\right\rangle\otimes K$. 
As coalgebras $H \cong k\left\langle g\right\rangle\otimes K^{\ast} \cong 
k\left\langle g\right\rangle\otimes H_{8}$, 
$H$ has $8$ grouplikes and we are done by the previous argument. 
\endproof

We now prove \textbf{Theorem \ref{th2}}:
\proof
\begin{enumerate}
\item Assume $\mathbf{G}\left( H\right)$ is abelian of order $8$. By Theorem 
\ref{th1}, $%
\mathbf{G}\left( H\right) \cong C_{2}\times C_{2}\times C_{2}$ or $\mathbf{G}%
\left( H\right) \cong C_{4}\times C_{2}$ and by Proposition \ref{th3} in this 
case $\mathbf{G}\left( H^{\ast }\right)$ is also abelian of order $8$. 
In Propositions \ref{4x2} and \ref{2x2x2} we have shown that there are exactly
$7$ nonisomorphic Hopf algebras with $\mathbf{G}\left( H\right) \cong
C_{4}\times C_{2}$ and at most $4$ nonisomorphic Hopf
algebras with $\mathbf{G}\left( H\right) \cong C_{2}\times C_{2}\times C_{2}$.
Now we show that there are $4$ distinct Hopf
algebras with $\mathbf{G}\left( H\right) \cong C_{2}\times C_{2}\times C_{2}$.
There are $2$ nonisomorphic examples of Hopf algebras with $%
\mathbf{G}\left( H\right) \cong \mathbf{G}\left( H^{\ast }\right) \cong
C_{2}\times C_{2}\times C_{2}$, namely $H_{8}\otimes kC_{2}$, which is not
triangular (if it were triangular, so would be $H_{8}$), and 
$\left( k\left( D_{8}\times C_{2}\right) \right) _{J}$,
which is triangular (see Section \ref{seccoc}), and $2$
more nonisomorphic Hopf algebras with $\mathbf{G}\left( H\right) \cong
C_{2}\times C_{2}\times C_{2}$ and $\mathbf{G}\left( H^{\ast }\right) \cong
C_{4}\times C_{2}$, namely $\left( H_{b:1}\right) ^{\ast }$ and $\left(
H_{c:\sigma _{1}}\right) ^{\ast }$. Comparing the structures of $H_{8}$ (see 
\cite[2.3, 2.4, 2.8]{Ma1}) and $H_{d:-1,1}$ we see that $H_{8}\otimes
kC_{2}\cong H_{d:-1,1}$, and therefore $H_{d:1,1}\cong \left( k\left(
D_{8}\times C_{2}\right) \right) _{J}$.
\item Assume that $\mathbf{G}\left( H\right)$ is nonabelian. Then, by Theorem
\ref{com_sub}, $H$ has a commutative subHopfalgebra of dimension $8$. By
Proposition \ref{D_8} $\mathbf{G}\left( H\right)\cong D_{8}$,
$\mathbf{G}\left( H^{\ast }\right) =C_{2}\times C_{2}$ and there are exactly
$2$ such Hopf algebras, $H_{C:1}\cong H_{B:1}^{\ast }$ and
$ H_{E}\cong H_{B:X}^{\ast }$. Comparing their $K_{0}$-rings with $K_{0}$-rings
of examples described in Section \ref{seccoc} we see that
$H_{C:1}\cong \left( kD_{16}\right) _{J}$ and 
$H_{E}\cong \left( kG_{2}\right) _{J}$.
\item Assume that $\mathbf{G}\left( H\right)$ is abelian of order $4$. 
By Theorem \ref{th1}, $\mathbf{G}\left( H\right) \cong C_{2}\times C_{2}$.
By Theorem \ref{com_sub}, $H$ has a commutative subHopfalgebra of dimension 
$8$ and therefore it was described in Section \ref{sec3}. 
There are exactly $3$ Hopf
algebras with this group of grouplikes: two of them, $H_{B:1}\cong H_{C:1}
^{\ast }$ and $ H_{B:X}\cong H_{E}^{\ast }$, have 
$\mathbf{G}\left( H\right)^{\ast }\cong D_{8}$ and one of them,
$H_{C:\sigma _{1}}$ has  $\mathbf{G}\left( H\right)^{\ast }\cong C_{2}\times 
C_{2}$ and therefore should be selfdual. Comparing the quotients of $H_{B}$ and
$kQ_{8}\#^{\alpha }kC_{2}$ we see that $H_{B:X}\cong kQ_{8}\#^{\alpha }kC_{2}$.
\end{enumerate}
\endproof


\begin{center}
ACKNOWLEDGMENTS
\end{center}
I would like to thank my Ph.D. advisor Professor Susan Montgomery for numerous
discussions, suggestions and comments about this paper.

Part of this paper is contained in my Ph.D. thesis in the University of 
Southern California. The rest of the work was done while I was a postdoctoral
fellow at MSRI. I am grateful to MSRI and the organizers of the 
Noncommutative Algebra program for the support.

I would also like to thank the referee for Proposition \ref{com_quot} and 
for the suggestion to use \cite{Ma7} in the proof of Theorem \ref{com_sub}.

\end{document}